\magnification=1000
\hsize=11.7cm
\vsize=18.9cm
\lineskip2pt \lineskiplimit2pt
\nopagenumbers

\hoffset=-1truein
\voffset=-1truein

\advance\voffset by 4truecm
\advance\hoffset by 4.5truecm

\newif\ifentete

\headline{\ifentete\ifodd	\count0 
      \rlap{\head}\hfill\tenrm\llap{\the\count0}\relax
    \else
        \tenrm\rlap{\the\count0}\hfill\llap{\head} \relax
    \fi\else
\global\entetetrue\fi}

\def\entete#1{\entetefalse\gdef\head{#1}}
\entete{Lluis Puig}

\input amssym.def
\input amssym.tex

\def\-{\hbox{-}}
\def\.{{\cdot}}

\def\K{{\cal K}}
\def\F{{\cal F}}
\def\E{{\cal E}}

\def\N{{\cal N}}

\def\T{{\cal T}}

\def\H{{\cal H}}
\def\s{{\cal S}}

\def\B{{\cal B}}
\def\W{{\cal W}}
\def\X{{\cal X}}

\def\Z{{\cal Z}}

\def\D{{\cal D}}

\def\ad{\frak a\frak c}

\def\Gr{\frak G\frak r}

\def\Fct{\frak F\frak c\frak t}
\def\Nat{\frak N\frak a\frak t}

\def\id{\frak i\frak d}
\def\int{\frak i\frak n\frak t}

\def\qq{\quad{\rm and}\quad}

\def\mod{\frak m\frak o\frak d}
\def\res{\frak r\frak e\frak s}
\def\ind{\frak i\frak n\frak d}

\def\too{\longrightarrow}

\def\Set{\frak S\frak e\frak t}

 3
 2
\font\large=cmr10  scaled \magstep 2
 2
\font\larti=cmti10  scaled \magstep 2
 2
\font\cds=cmr7

\count0=1

\centerline{\large The Hecke algebra}
\medskip
\centerline{\large of a Frobenius {\larti P}-category}
\medskip
\centerline{\bf Lluis Puig }

\medskip
\centerline{\cds CNRS, Institut de Math\'ematiques de Jussieu}
\smallskip
\centerline{\cds 6 Av Bizet, 94340 Joinville-le-Pont, France}
\smallskip
\centerline{\cds puig@math.jussieu.fr}

\bigskip
\noindent
{\bf £1. Introduction}

\bigskip
£1.1. Let $p$ be a prime number and $P$ a finite $p\-$group; as usual, denote by $\F_{\!P}$ the {\it Frobenius
category\/} of $P$ [4,~I\hskip1pt1], namely the category where the objects are the subgroups of $P$ and the morphisms are
the group homomorphisms between them induced by the inner automorphisms of $P\,.$  Recall that  a {\it divi-sible $P\-$category\/}
$\F$ [4,~2.3] is a subcategory of the {\it category of finite groups\/} $\Gr$ such that it contains  $\F_{\!P}$ and has the same objects, that all the morphisms are injective group homomorphisms and that, for any subgroup $Q$ of $P\,,$ the category 
$\F_{\mid Q}$ of {\it $\F\-$morphisms to~$Q$\/} is a {\it full\/} subcategory of~$\Gr_{\mid Q}$ [4,~1.7].

\medskip
£1.2. More precisely, denoting by $\F (Q)$ the group of $\F\-$automorphisms of any subgroup $Q$ of $P\,,$ 
recall that then $\F$ is a  {\it Frobenius $P\-$category\/} [4,~Proposition~2.11] if $\F_{\!P}(P)$ is a Sylow $p\-$subgroup 
of $\F (P)$ and if any $\F\-$morphism 
$\varphi\,\colon Q\to P$ fulfilling
$$\zeta\Big(C_P\big(\varphi (Q)\big)\Big) = C_P \big((\zeta\circ\varphi) (Q)\big)
\eqno £1.2.1\phantom{.}$$
for any $\F\-$morphism $\zeta\,\colon \varphi(Q)\.C_P\big(\varphi (Q)\big)\to P\,,$ can be extended to an $\F\-$mor-phism
from the converse image in $P$ of the intersection $\F_{\!P} (Q)\cap {}^{\varphi^*}\!\F_{\!P}\big(\varphi (Q)\big)\,,$
where we denote by $\varphi^*\,\colon \varphi (Q)\cong Q$ the inverse of the isomorphism induced by $\varphi\,,$
and by ${}^{\varphi^*}\!\F_{\!P}\big(\varphi (Q)\big)$ the corresponding image of $\F_{\!P}\big(\varphi (Q)\big)$ in ${\rm Aut}(Q)\,.$

\medskip
£1.3. On the one hand, in a Frobenius $P\-$category $\F$  still holds a gene-ral form for the {\it Alperin Fusion Theorem\/}
which, in some sense, provides a ``decomposition'' of the $\F\-$morphisms [4,~Corollary~5.14]. On the other hand,
restricting ourselves to the {\it full\/} subcategory $\F^{^{\rm sc}}$ over the subgroups $Q$ of $P$ fulfilling 
$C_P\big(\varphi (Q)\big) \i \varphi (Q)$ for any $\varphi\in \F(P,Q)\,,$ in the quotient category $\tilde\F^{^{\rm sc}}$
of $\F^{^{\rm sc}}$ by the inner automorphisms of the objects  the morphisms admit a ``canonical decomposition'' 
in terms of suitable subsets $\tilde\F^{^{\rm sc}}(Q,R)_{\tilde\iota_R^P}$ of $\tilde\F^{^{\rm sc}}\-$morphisms
as it follows from [4,~Proposition~6.7]. An old aim we had was to find a common framework where we could formulate
both facts.

\medskip
£1.4. More recently, we have proved that the Frobenius $P\-$categories are just an {\it avatar\/} of the {\it basic $P\times P\-$sets\/}
[4,~21.4] --- namely the finite nonempty $P\times P\-$sets $\Omega$ such that $\{1\}\times P$ acts {\it freely\/} on $\Omega\,,$ we have
$$\Omega^\circ \cong \Omega\qq \vert\Omega\vert/\vert P\vert \not\equiv 0 \bmod{p}
\eqno £1.4.1\phantom{.}$$
where we  denote by $\Omega^\circ$ the $P\times P\-$set obtained exchanging both factors
and, for any subgroup $Q$ of $P$ and any group homomorphism $\varphi\,\colon Q\to P$ such that $\Omega$ contains
a $P\times P\-$subset isomorphic to $(P\times P)/\Delta_\varphi (Q)\,,$ we have a $Q\times P\-$isomorphism
$${\rm Res}_{\varphi\times {\rm id}_P} (\Omega)\cong {\rm Res}_{\iota_Q^P\times {\rm id}_P} (\Omega)
\eqno £1.4.2,$$
where we set $\Delta_\varphi (Q)  =\{(\varphi (u),u)\}_{u\in Q}$ and denote by $\iota_Q^P$
the corresponding inclusion map. The last but not the least, their generalization by Kari Ragnarsson and Radu Stancu to {\it virtual\/} 
$P\times P\-$sets [6] has driven our attention upon the so-called {\it double Burnside ring of $P\,.$\/}

\medskip
£1.5. In this paper we introduce a new {\it avatar\/} of a Frobenius $P\-$category~$\F$ in the form of a suitable subring $\H_\F$
of the double Burnside ring of~$P$ --- called the {\it Hecke algebra of $\F$\/} --- where we are able to formulate 
{\it all\/} the results mentioned above, in some sense providing the framework we were looking for. The choice of the name 
comes from the case where $\F$ is the Frobenius category of a finite Chevalley group in the defining characteristic or, 
more generally, of a finite group $G\,,$ where an analogous definition of the  {\it Hecke algebra\/} is usually considered;
 in this case, the {\it Hecke algebra\/} of~$G$ is isomorphic to a quotient of a suitable {\it extension\/} of~$ \H_\F$ (see \S 8 below).

\bigskip
\noindent
{\bf £2. The quoted framework}
\bigskip
£2.1. Let us denote by $P\-\Set$ the category of finite $P\-$sets endowed with the {\it disjoint union\/}
and with the {\it inner\/} direct product mapping any pair of finite $P\-$sets $X$ and $Y$ on the {\it $P\-$set\/}
--- still noted $X \times Y$ --- obtained from the $P\times P\-$set $X\times Y$ restricted through the {\it diagonal\/} map 
$\Delta \,\colon P\to P\times P\,;$
note that, for a third finite $P\-$set $Z\,,$ we have a {\it canonical\/} bijection
$${\rm Hom}_P (Z,X\times Y)\cong {\rm Hom}_P(Z,X)\times {\rm Hom}_P(Z,Y)
\eqno £2.1.1.$$
Then, we denote by $\Fct(P\-\Set)$ the category of functors from $P\-\Set$ to itself preserving the 
disjoint unions, endowed with the {\it  disjoint union\/} induced by the  disjoint union in $P\-\Set$ and with the 
{\it composition\/} of functors.

\medskip
£2.2. More generally, if $P'$ is a second finite $p\-$group, we often have to consider the category 
$\Fct (P\-\Set,P'\-\Set)$ of functors from $P\-\Set$ to~$P'\-\Set$ preserving the 
disjoint union, which is  also endowed with a  {\it  disjoint union\/};
let us say that a functor $\frak f\,\colon P\-\Set\to P'\-\Set$ is {\it indecomposable\/} if it is {\it not\/} 
the disjoint union  of two {\it nonempty\/} functors in~$\Fct(P\-\Set,P'\-\Set)\,.$ If $P''$ is a third finite $p\-$group,
 we clearly have a {\it distributivity\/}, namely 
$$(\frak g \sqcup \frak g')\circ \frak f = (\frak g\circ\frak f)\sqcup (\frak g'\circ\frak f)\qq
\frak g\circ (\frak f \sqcup \frak f') = (\frak g\circ \frak f )\sqcup (\frak g\circ \frak f')
\eqno £2.2.1\phantom{.}$$
for any functors $\frak f$ and $\frak f'$ in $\Fct (P\-\Set,P'\-\Set)$ and  any functors $\frak g$ and $\frak g'$ in  $\Fct(P'\-\Set,P''\-\Set)\,.$ 
We denote by $\Nat (\frak f,\frak f')$ the set of {\it natural maps\/} from $\frak f$ to $\frak f'\,,$ but write 
$\Nat (\frak f)$ instead of~$\Nat (\frak f,\frak f)\,;$ it is clear that for a third functor 
$\frak f''\,\colon P\-\Set\to P'\-\Set$ we have
$$\Nat (\frak f\sqcup\frak f',\frak f'')\cong \Nat (\frak f,\frak f'')\times \Nat (\frak f',\frak f'')
\eqno £2.2.2.$$

\medskip
£2.3. On the other hand, if $\alpha\,\colon P'\to P$ is a group homomorphism,
the {\it restriction\/} defines an evident functor
$$\res_\alpha : P\-\Set\too  P'\-\Set
\eqno £2.3.1\phantom{.}$$
which admits a {\it left adjoin\/} functor
$$\ind_\alpha : P'\-\Set\too  P\-\Set
\eqno £2.3.2;$$
hence,  for any functors $\frak g$ in $\Fct (P''\-\Set,P'\-\Set)$  and $\frak h$ in $\Fct (P''\-\Set,P\-\Set)$ 
we have a {\it canonical\/} bijection
$$\Nat(\ind_\alpha\circ\frak g,\frak h)\cong \Nat(\frak g, \res_\alpha\circ \frak h)
\eqno £2.3.3.$$

\medskip
£2.4. Note that, any pair of functors $\frak f$ and $\frak f'$ in $\Fct(P\-\Set,P'\-\Set)$ defines a functor
$$\frak f\times\frak f' : P\-\Set\times P\-\Set \too P'\-\Set
\eqno £2.4.1\phantom{.}$$
mapping any pair of finite $P\-$sets $X$ and $X'$ on the $P'\-$set $\frak f (X)\times \frak f' (X')\,,$
and any pair of $P\-$set maps $a\,\colon X\to Y$ and $a'\,\colon X'\to Y'$ on the evident $P'\-$set map
$$\frak f(a)\times \frak f' (a') : \frak f (X)\times \frak f' (X')\too \frak f (Y)\times \frak f' (Y')
\eqno £2.4.2;$$
then, for any pair of functors $\frak g$ and $\frak g'$ in $\Fct(P\-\Set,P'\-\Set)$ we have an evident injection
$$\Nat (\frak g,\frak f)\times \Nat (\frak g',\frak f')\too \Nat (\frak g\times \frak g',\frak f\times \frak f')
\eqno £2.4.3.$$
which is bijective if $\frak g$ and $\frak g'$ both map the {\it trivial\/} $P\-\Set$
on the {\it trivial\/} $P'\-\Set\,.$
Indeed, according to the bijection~£2.1.1, any $P\-$set map
$$\frak g (X)\times \frak g'(X')\too \frak f (X)\times \frak f'(X')
\eqno £2.4.4\phantom{.}$$
determines and is determined by the corresponding two $P\-$set maps
$$\frak g (X)\times \frak g'(X')\too \frak f (X)\qq \frak g (X)\times \frak g'(X')\too \frak f'(X')
\eqno £2.4.5;$$
but, it follows from the {\it naturalness\/} that the left-hand map in~£2.4.5 is determined by the map which corresponds to 
the case where $X'$ is the {\it trivial\/} $P\-$set and then it determines an element of $\Nat (\frak g,\frak f)\,;$
similarly,  the right-hand map in~£2.4.5  determines an element of $\Nat (\frak g',\frak f')\,;$ finally, it is easily checked that
this correspondence is the inverse of the injection~£2.4.3.
\eject

\medskip
£2.5. Actually, we are only interested in the {\it full\/} subcategory $\frak B_{P\times P'}$ of~$\Fct(P\-\Set,P'\-\Set)$ over the
functors determined by $P'\times P\-$sets; explicitly, a $P'\times P\-$set $\Omega$ determines a functor
$\frak f_\Omega$ mapping any $P\-$set $X$ on the $P'\-$set
--- noted $\Omega \times_P X $ --- formed by the set of orbits of $P$ on $\Omega\times X$ through the
action sending $(\omega,x)\in \Omega\times X$ to $\big((1,u)\.\omega, u\.x\big)$ for any $u\in P\,,$
endowed with the action of~$P'$ mapping the orbit of $(\omega,x)$ on the orbit of $\big((u',1)\.\omega,x\big)$
for any $u'\in P$ --- we often write $u'\.\omega$ and $\omega\.u^{-1}$ instead of $(u',1)\.\omega$
and $(1,u)\.\omega$ respectively; it is clear that $\frak f_\Omega$ preserves the disjoint unions. For instance, for $P'' = P\,,$  the {\it restriction\/} and the {\it induction\/} functors above
$$\res_\alpha : P\-\Set\too P'\-\Set\qq \ind_\alpha : P'\-\Set\too P\-\Set
\eqno £2.5.1$$
are respectively determined by the $P'\times P\-$ and the $P\times P'\-$sets ${}_\alpha P$ and $P_\alpha$
denoting  the respective $P'\times P\-$ and $P\times P'\-$set structures over $P$ sending
$v\in P$ to $\alpha (u')vu^{-1}$ and to $uv\alpha (u')^{-1}$ for any $u\in P$ and any $u'\in P'\,.$ 
Note that the functor $\frak f_\Omega$ determines the $P\times P\-$set
$\Omega\,;$ indeed, considering the $P\-$set $P$ and its group of automorphisms in $P\-\Set$ --- which is
isomorphic to~$P$ --- we get a $P'\times P\-$set structure on $\frak f_\Omega (P)$ and, with this structure, 
it is clear that $\frak f_\Omega (P)$ is isomorphic to $\Omega\,.$ In particular, for any $P'\times P\-$set $\Omega'\,,$
we have a canonical bijection 
$$\Nat (\frak f_\Omega,\frak f_{\Omega'})\cong {\rm Hom}_{P'\times P}(\Omega,\Omega')
\eqno £2.5.2.$$

\medskip
£2.6.  It is clear that $\frak B_{P\times P'}$ is {\it closed\/} with respect to the
disjoint union,  namely we have $\frak f_\Omega \sqcup \frak f_{\Omega'} = \frak f_{\Omega \sqcup\Omega'}\,,$
 and an {\it indecomposable\/} functor in $\frak B_{P\times P'}$ is given by the {\it transitive\/} $P'\times P\-$set
$(P'\times P)/D$ where $D$ is a subgroup of~$P'\times P\,;$ shortly we set
$$\frak f_D = \frak f_{(P'\times P)/D}
\eqno £2.6.1.$$ 
Then, for any $P'\times P\-$set $\Omega\,,$ it follows from isomorphism~£2.5.2 that a natural
map $\mu\,\colon \frak f_D\to \frak f_\Omega$ is induced by a $P'\times P\-$set map (cf.~£2.5.2)
$${\rm m}_\mu : (P'\times P)/D\too \Omega
\eqno £2.6.2\phantom{.}$$
which is actually determined by the image in $\Omega$ of the class in $(P'\times P)/D$ of~$(1,1)\,,$  and 
this image clearly belongs to the set of $D\-$fixed elements $\Omega^{D}\,;$ conversely, any element $\omega$ 
of~$\Omega^{D}$ determines  a $P'\times P\-$set map from $(P'\times P)/D$ to $\Omega$ sending 
the class of~$(1,1)$ to $\omega\,,$  and therefore it determines  a natural map from $\frak f_{D}$ 
to $\frak f_\Omega\,;$ thus, we have a bijection
$$\Nat (\frak f_{D},\frak f_\Omega)\cong \Omega^{D}
\eqno £2.6.3.$$

\medskip
£2.7. In particular, it is  clear that the group
$$\Z(\frak f_{D}) =   C_{P'\times P}(D)\big/Z(D)
\eqno £2.7.1\phantom{.}$$
 acts on the set of natural maps $\Nat (\frak f_{D},\frak f_{\Omega})\,.$ Actually, all these subgroups 
 $\Z(\frak f_{D})\i \Nat(\frak f_{D})\,,$ where $D$ runs over the set of  subgroups of $P'\times P\,,$\break
 \eject
 \noindent
 determine  an {\it interior structure\/} in  $(\frak B_{P\times P'})^{\!^\circ}$ [4,~1.3] and, coherently, for another  subgroup $E$ 
 of $P'\times P\,,$  we set
 $$\widetilde\Nat (\frak f_D,\frak f_E) =  \Z(\frak f_{D})\big\backslash\Nat (\frak f_D,\frak f_E)
 \eqno £2.7.2.$$
  Note that, for any pair of functors $\frak g$ and $\frak g'$ in $\frak B_{P\times P'}$ we get
$$\Nat (\frak f_{D},\frak g\sqcup\frak g')\cong \Nat (\frak f_{D},\frak g)\sqcup \Nat (\frak f_{D},\frak g')
\eqno £2.7.3.$$

\medskip
£2.8. Moreover, $\frak B_{P\times P}$ is also  {\it closed\/} with respect to the composition, namely we have
$$\frak f_\Omega \circ \frak f_{\Omega'} = \frak f_{\Omega \times_P\Omega'}
\eqno £2.8.1\phantom{.}$$
for any $P\times P\-$sets $\Omega$ and $\Omega'\,.$ On the other hand, denoting by $\Omega^{^{\!\circ}}$ 
the {\it opposite\/} $P\times P\-$set of a $P\times P\-$set $\Omega\,,$ the correspondence mapping 
$\frak f_\Omega$ on $(\frak f_\Omega)^{^{\!\circ}} = \frak f_{\Omega^{^{\circ}}}$ defines an auto-equivalence 
$\frak t$  of $\frak B_{P\times P}$ which preserves the disjoint union and reverses  the composition; then, it follows  from 
bijection~£2.5.2 that, for any pair of functors $\frak f$ and $\frak g$ in $\frak B_{P\times P}\,,$ we still have a canonical bijection
$$\Nat (\frak f^{^\circ},\frak g^{_\circ})\cong \Nat (\frak f,\frak g)
\eqno £2.8.2.$$
Let us call {\it symmetric\/} any $\frak t\-$stable subcategory of $\frak B_{P\times P}\,.$

\medskip
£2.9. A particular kind of functors in $\frak B_{P\times P}$ is provided by the {\it inner\/} direct product in $P\-\Set\,;$
indeed, for any $P\-$set $X\,,$ the correspondence sending any $P\-$set $Y$ to the $P\-$set $X\times Y$ and any $P\-$set map
$g\,\colon Y\to Y'$ to the $P\-$set map
$${\rm id}_X\times g : X\times Y\too X\times Y'
\eqno £2.9.1\phantom{.}$$
defines a functor $\frak m_X$ from $P\-\Set$ to $P\-\Set$ preserving the disjoint unions, and it is easily checked that, for any  subgroup 
$Q$ of $P\,,$ the functor determined by the $P\-$set $P/Q$ coincides with the functor determined by the $P\times P\-$set
$P\times _Q P\,;$ more precisely, we have a {\it functor\/}
$$\frak m^P : P\-\Set\too \frak B_P
\eqno £2.9.2\phantom{.}$$
mapping any $P\-$set $X$ on $\frak m_X$ and  any $P\-$set map $f\,\colon X\to X'$ to the {\it natural map\/}
$\mu_f\,\colon \frak m_X\to \frak m_{X'}$ sending $Y$ to the $P\-$set map 
$$f\times {\rm id}_Y : X\times Y\too X'\times Y
\eqno £2.9.3.$$
Note that $\frak m^P$ preserves the {\it disjoint unions\/}, maps the {\it inner\/} direct product in~$P\-\Set$ 
on the {\it composition\/} in $\frak B_{P\times P}\,,$ and is invariant by the auto-equivalence $\frak t$ of~$\frak B_{P\times P}$ above.

\medskip
£2.10.  Actually, $P$ has an obvious $P\times P\-$set structure and it is clear that $\frak f_P$ is the {\it identity functor\/}.
More generally, for any subgroup $Q$ of $P$ and any group homomorphism $\varphi\,\colon Q\to P\,,$ as in~£1.4 above we set
$$\Delta_\varphi (Q) = \big\{\big(\varphi(u),u\big)\}_{u\in Q}
\eqno £2.10.1\phantom{.}$$
and, coherently with our notation in~£2.5, we still set
$$P_\varphi \times_Q P = (P\times P)/\Delta_\varphi (Q)\qq 
\frak f_{Q,\varphi} = \frak f_{P_\varphi\times_Q P} = \ind_\varphi\circ \res_{\iota_Q^P}
\eqno £2.10.2;$$
it is quite clear that this functor preserves the {\it projective\/} objects in $P\-\Set$ or, equivalently, the $P\-$sets
 where $P$ acts freely; conversely, any {\it indecomposable\/} functor in~$\frak B_P$ which  preserves the {\it projective\/} 
 $P\-$sets has this form.

 \medskip
 £2.11. More explicitly, any functor $\frak f$ in $\frak B_{P\times P}$ preserving the {\it projective\/} $P\-$sets has the form 
 $\frak f   = \bigsqcup_{i\in I} \frak f_{Q_i,\varphi_i}$  for a suitable finite family $\{Q_i,\varphi_i\}_{i\in I}$ where 
 $Q_i$ are subgroups of $P$ and  $\varphi_i\,\colon Q_i\to P$ are group homomorphisms; then, for another functor  $\frak g$ in $\frak B_P$ preserving the {\it projective\/} $P\-$sets, given by the  finite family  $\{R_j,\psi_j\}_{j\in J}\,,$
 there is a natural map $\mu\,\colon \frak g\to \frak f$ if and only~if~there are a map
 $m\,\colon J\to I$ and, for any $j\in J\,,$ a natural map $\mu_j\,\colon 
 \frak f_{R_j,\psi_j}\to \frak f_{Q_{m(j)},\varphi_{m(j)}}$ or, equivalently, a $P\times P\-$set map
$${\rm m_j} : P_{\psi_j}\times_{R_j} P\too  P_{\varphi_{m(j)}}\times_{Q_{m(j)}}  P
\eqno £2.11.1,$$
which amounts to saying that there are $u_j,v_j\in P$ such that 
$$R_j\i (Q_{m(j)})^{u_j}\qq \psi_j = \kappa_{v_j}\circ\varphi_{m(j)}\circ \kappa_{u_j}
\eqno £2.11.2\phantom{.}$$
where $\kappa_{u_j}\,\colon R_j\to Q_{m(j)}$ and $\kappa_{v_j}\,\colon P\to P$ respectively  denote the
conjugation by $u_j$ and by $v_j\,.$

\medskip
£2.12. In particular, for any indecomposable functors $\frak f_{Q,\varphi}$ and $\frak f_{R,\psi}$ in $\frak B_{P\times P}$
such that $(\frak f_{Q,\varphi})^{^\circ}$ and $(\frak f_{R,\psi})^{^\circ}$ still preserve the {\it projective\/} $P\-$sets,
denoting by $\varphi_*\,\colon Q\cong \varphi(Q)$ and $\psi_*\,\colon R\cong\psi(R)$ the corresponding isomorphisms
and setting $\varphi^* =(\varphi_*)^{-1}$ and $\psi^* = (\psi_*)^{-1}\,,$ it follows from condition~£2.11.2 
above that any natural map from $\frak f_{R,\psi}$ to $\frak f_{Q,\varphi}$ determines an element of the intersection
$$\tilde\F_P (Q,R)\cap \big(\tilde\varphi^*\circ\tilde\F_P \big(\varphi (Q),\psi(R)\big)\circ \tilde\psi_*\big)
\eqno £2.12.1,$$
namely an element $\tilde\kappa_u$ in $\tilde\F_P (Q,R)$ such that the composition $\tilde\varphi_*\circ\tilde\kappa_u\circ\tilde\psi^*$ 
belongs to $\tilde\F_P \big(\varphi (Q),\psi(R)\big)\,;$
moreover, it is easily checked that this  cor-respondence determines a bijection
$$\widetilde\Nat (\frak f_{R,\psi},\frak f_{Q,\varphi}) \cong \tilde\F_P (Q,R)\cap \big(\tilde\varphi^*\circ\tilde\F_P \big(\varphi (Q),\psi(R)\big)\circ \tilde\psi_*\big)
\eqno £2.12.2.$$

 \medskip
 £2.13. Moreover, for any pair of subgroups $Q$ and $R$ of $P\,,$ and any pair of group homomorphisms $\varphi$ and $\psi$
 from $Q$ and $R$ to $P$ respectively, the functor  $\frak f_{Q,\varphi}\circ \frak f_{R,\psi}$ is determined by the $P\times P\-$set
$$\eqalign{(P_\varphi \times_Q P)\times_P (P_\psi\times_RP)
&= P_\varphi\times_Q \big( P\times_P (P_\psi\times_R P)\big)\cr
&= P_\varphi\times_Q \big(( P\times_P P_\psi)\times_R P\big)\cr
&=P_\varphi\times_Q ( P_\psi\times_R P)\cr}
\eqno £2.13.1;$$
but, choosing a set of representatives $W\i P$ for the set of double classes $Q\backslash P/\psi(R)\,,$ we have the disjoint decomposition
$$ P = \bigsqcup_{w\in W} Q\.w\.\psi(R)
\eqno £2.13.2;$$
moreover, for any $w\in W\,,$  it is clear that $ Q\.w\.\psi(R)$ has a $Q\times R\-$set structure sending $uw\psi(v)^{-1}$ 
to $u'uw\psi(v'v)^{-1}$ for any $u,u'\in Q$ and any $v,v'\in R\,;$ more precisely, setting 
$U_w = \psi^{-1}\big(Q^w\cap \psi (R)\big)$ and denoting by $\psi_w\,\colon U_w\to Q^w$ the restriction of $\psi$
and by $\kappa_w\,\colon Q^w\to Q$ the corresponding conjugation 
by $w\,,$  it is clear that we have a $Q\times R\-$set isomorphism
$$Q\.w\.\psi(R) \cong Q_{\kappa_w\circ\psi_w}\times_{U_w}  R
\eqno £2.13.3;$$
finally, we get a $P\times P\-$set isomorphism
$$\eqalign{(P_\varphi \times_Q P)\times_P (P_\psi\times_RP)
&\cong \bigsqcup_{w\in W} P_\varphi\times_Q (Q_{\kappa_w\circ\psi_w}\times_{U_w}  R)\times_R P\cr
&\cong \bigsqcup_{w\in W} P_{\varphi\circ\kappa_w\circ\psi_w}\times_{U_w} P\cr}
 \eqno £2.13.4.$$

 \bigskip
 \noindent
 {\bf £3. A version for divisible and Frobenius $P\-$categories} 
\bigskip
£3.1.  For any {\it divisible $P\-$category\/} $\F\,,$ let us denote by $\frak H_\F$ the {\it full\/} sub-category  
of $\frak B_{P\times P}$ over the functors which are isomorphic to disjoint unions of functors $\frak f_{Q,\varphi}$
 where $Q$ is a subgroup of $P$ and $\varphi\,\colon Q\to P$ is an $\F\-$morphism; thus, $\frak H_\F$ is closed
 with respect to the disjoint unions, its elements preserve the {\it projective\/} $P\-$sets, and is clearly symmetric (cf.~£2.8). Moreover,  $\frak H_\F$ is closed with respect to the composition; indeed, for a second subgroup $R$ of~$P$ and any $\F\-$morphism 
 $\psi\,\colon R\to P\,,$ with the notation in~£2.13 above we have
 $$\frak f_{Q,\varphi}\circ \frak f_{R,\psi} = \bigsqcup_{w\in W} \frak f_{U_w,\varphi\circ\kappa_w\circ\psi_w}
 \eqno £3.1.1\phantom{.}$$
 and the divisibility in $\F$ implies that, for any $w\in W\,,$ $\psi_w$ and $\kappa_w$ are $\F\-$morphisms
 and therefore $\varphi\circ\kappa_w\circ\psi_w$ is also an $\F\-$morphism. In particular, since $\F$ contains the  
 Frobenius category $\F_P$ of $P\,,$ $\frak H_\F$ contains $\frak H_{\F_P}$ --- for short, we write $\frak H_P$
  instead of $\frak H_{\F_{\!P}}\,;$ note that $\frak H_P$ coincides with the {\it image\/} of the functor $\frak m^P\,.$
  \eject

 \medskip
 £3.2. More generally, the divisibility of $\F$ implies that 
 \smallskip
 \noindent
 £3.2.1\quad {\it Any functor $\frak f$ in~$\frak B_{P\times P}$ admitting   
 a natural map to a functor in  $\frak H_\F$ is an $\frak H_\F\-$object.\/}
  \smallskip
 \noindent
Indeed, we actually may assume that $\frak f = \frak f_{D}$ where $D$ is a subgroup of $P\times P\,;$ then, our hypothesis and £2.11 above imply that
there are a functor $\frak f_{Q,\varphi}\,,$ where $Q$ is a subgroup of $P$ and $\varphi\,\colon Q\to P$ is an 
$\F\-$morphism, and a natural map $\mu\,\colon  \frak f_{D}\to \frak f_{Q,\varphi}\,,$ which amounts to 
saying that there are $u,v\in P$ such that 
$$D\i \Delta_\varphi (Q)^{(u,v)}
\eqno £3.2.2;$$
since $\F$ is divisible and contains $\F_{\!P}\,,$ it easily follows that $D = \Delta_\psi (R)$ for a subgroup $R$ of $Q^v$ and
an $\F\-$morphism $\psi\,\colon R\to P\,.$ Actually, this condition is equivalent to the divisibility of $\F$ as our next result shows.

 \bigskip
 \noindent
 {\bf Proposition~£3.3.\/} {\it  Let $\frak H$ be a full symmetric subcategory of $\frak B_{P\times P}$ 
 containing~$\frak H_P\,,$ where all the elements preserve the projective $P\-$sets, closed with respect to disjoint unions and to composition, and fulfilling the following condition:
 \smallskip
\noindent
{\rm £3.3.1.}\quad  Any functor $\frak f$ in~$\frak B_{P\times P}$ admitting   
 a natural map to a functor in $\frak H$  is an $\frak H\-$object. 
 \smallskip
\noindent
Then, the family of sets 
$$\F_{\!\frak H} (P,Q) = \{\varphi\in {\rm Hom}(Q,P)\mid \hbox{$\frak f_{Q,\varphi}$ is an object of $\frak H$}\}
\eqno £3.3.2,$$
 where $Q$ runs over the set of subgroups of $P\,,$  determine a  {\it divisible $P\-$cate-gory~$\F_{\!\frak H}\,.$\/}\/}

\medskip
\noindent
{\bf Proof:} Since $\frak H$ contains $\frak H_P\,,$ it is clear that  $\F_{\!\frak H} (P,Q)$ contains $\F_P (P,Q)\,;$
moreover, following [1,~2.4.1], let $Q$ and $R$ be subgroups of~$P$ and
consider a group homomorphism $\theta\,\colon R\to Q$ such that 
$$\F_\frak H (P,R)\cap \big(\F_\frak H (P,Q)\circ \theta\big)\not= \emptyset
\eqno £3.3.3;$$
that is to say, we assume that there is a group homomorphism $\varphi\,\colon Q\to P$ such that
$\frak f_{Q,\varphi}$ and $\frak f_{R,\varphi\circ \theta}$ are $\frak H\-$objects;  hence, the composition
$(\frak f_{Q,\varphi})^{^{\!\circ}} \circ\frak f_{R,\varphi\circ \theta}$ is an $\frak H\-$object too; but, it is clear that
$$(\frak f_{Q,\varphi})^{^{\!\circ}} = \frak f_{\varphi (Q),\iota_Q^P\circ \varphi^*}
\eqno £3.3.4\phantom{.}$$
 where, as in~£1.2 above, $\varphi^*$
denotes the inverse of the isomorphism $Q\cong \varphi(Q)$ determined by~$\varphi\,,$ and then it follows from equality~£3.1.1
that we have a natural map
$$\frak f_{R,\iota_Q^P\circ\theta}\too (\frak f_{Q,\varphi})^{^{\!\circ}} \circ\frak f_{R,\varphi\circ \theta}
\eqno £3.3.5;$$
consequently, $\frak f_{R,\iota_Q^P\circ\theta}$ is also an $\frak H\-$object.

\smallskip
Now, since for any $\varphi'\in \F_\frak H (P,Q)$ the functor $\frak f_{Q,\varphi'}\circ \frak f_{R,\iota_Q^P\circ\theta}$
is an $\frak H\-$object and  it follows again from equality~£3.1.1
that we have a natural map
$$\frak f_{R,\varphi'\circ\theta}\too \frak f_{Q,\varphi'} \circ\frak f_{R,\iota_Q^P\circ \theta}
\eqno £3.3.6,$$
$\frak f_{R,\varphi'\circ\theta}$ is an $\frak H\-$object  too; in conclusion, we get
$$\F_\frak H (P,Q)\circ \theta\i \F_\frak H (P,R)
\eqno £3.3.7,$$
so that the proposition follows from [1,~2.4.1].

\medskip
£3.4. Let us call {\it divisible $\frak H_P\-$category\/} any full symmetric subcategory $\frak H$ of~$\frak B_{P\times P}$ 
 containing~$\frak H_P\,,$ where all the elements preserve the projective $P\-$sets,  closed with respect to disjoint unions and to composition, and fulfilling  condition~£3.3.1; then,
the proposition above guarantees that $\frak H$ determines a {\it divisible $P\-$category\/} $\F_{\!\frak H}$ 
and we already have seen that any {\it divisible $P\-$category\/} $\F$ determines a {\it divisible $\frak H_P\-$category 
$\frak H_\F$\/}. Moreover,
it is quite clear that
$$\frak H_{\F_{\!\frak H}} = \frak H\qq \F_{\frak H_\F} = \F
\eqno £3.4.1,$$
so that these correspondences determine a bijection between the sets of {\it divi-sible $P\-$categories\/} and
 {\it divisible $\frak H_P\-$categories;\/} we discribe below the image of the {\it Frobenius $P\-$categories.\/}
 Recall that we denote by $\tilde\F$  the quotient category of $\F$ by the inner automorphisms of the objects [4,~1.3].

\medskip 
£3.5. First of all,  for any functor $\frak f$ in  $\frak B_{P\times P}$ which preserves the {\it projective\/} $P\-$sets
 we define a positive integer $\ell (\frak f)$ by
$$\vert\frak f (P)\vert = \ell (\frak f)\vert P\vert
\eqno £3.5.1;$$
for another such a functor $\frak g$ in  $\frak B_{P\times P}\,,$ we have
$$\ell (\frak f\sqcup \frak g) = \ell (\frak f) + \ell (\frak g)\qq  \ell (\frak f\circ \frak g) = \ell (\frak f)\,\ell (\frak g)
\eqno £3.5.2.$$
Then, for any  {\it divisible $\frak H_P\-$category $\frak H$\/} denote by 
$$\ell_\frak H : \frak H\too \Bbb N
\eqno £3.5.3\phantom{.}$$
the restriction  to $\frak H$ of this correspondence and by ${\rm Ker}(\ell_\frak H)$ the set of isomorphism classes of functors
$\frak f$ in $\frak H$ fulfilling $\ell (\frak f) = 1\,;$ it is quite clear that the composition in $\frak H$ induces a group 
structure in ${\rm Ker}(\ell_\frak H)$ and we actually have
$${\rm Ker}(\ell_\frak H) \cong \tilde\F_{\!\frak H} (P)
\eqno £3.5.4.$$

\medskip
£3.6. On the other hand, according to~£2.6 above, for any indecomposable functor $\frak f_{D}$ 
in~$\frak B_{P\times P}\,,$ a natural map from $\frak f_D$ to $\frak f_D$ is determined by the
image $\overline{(u,v)}\in (P\times P)/D$ of~$\overline{(1,1)}$ and then it is easily checked that the element
$(u,v)\in P\times P$ normalizes the subgroup $D$ of $P\times P\,;$ conversely, any 
$(u,v)\in N_{P\times P}(D)$ determines a $P\times P\-$set automorphism of
$(P\times P)/D$ and therefore it determines a natural automorphism 
$\frak f_D\cong \frak f_D\,,$ so that we have
$$\Nat (\frak f_D)\cong \bar N_{P\times P}(D)
\eqno £3.6.1$$
\eject
\noindent
and  $\Nat (\frak f_D)$ acts on $(P\times P)/D\,.$
We are specially interested in the case where $\frak f_D$ preserves the {\it projective\/} $P\-$sets or, equivalently, 
where $D = \Delta_\varphi (Q)$ for a subgroup $Q$ of $P$ and a group homomorphism $\varphi\,\colon Q\to P\,;$
then, from bijection~£2.12.2, we get the exact sequence (cf.~£2.7)
$$1\too \Z (\frak f_{Q,\varphi})\too \Nat (\frak f_{Q,\varphi})\too 
\tilde\F_P (Q)\cap {}^{\varphi^*}\!  \tilde\F_P\big(\varphi (Q)\big)\too 1
\eqno £3.6.2.$$

\bigskip
\noindent
{\bf Proposition~£3.7.} {\it Let $\frak H$ be a divisible $\frak H_P\-$category. Then $\F_{\!\frak H}$ is a 
Frobenius $P\-$category if and only if  $\,{\rm Ker}(\ell_\frak H)$ is a $p'\-$group and, for any indecomposable functor 
$\frak f_{Q,\varphi}$ in $\frak H$ such that we have
$$\vert \Z (\frak f_{Q,\varphi'})\vert\le \vert \Z (\frak f_{Q,\varphi})\vert
\eqno £3.7.1\phantom{.}$$
for any indecomposable functor $\frak f_{Q,\varphi'}$ in $\frak H\,,$ there are an indecomposable functor 
$\frak f_{R,\psi}$ in $\frak H$ and  a natural map $\mu$ from 
$\frak f_{Q,\varphi}$ to $\frak f_{R,\psi}$ such that  the stabilizer of $\mu\in \Nat (\frak f_{Q,\varphi}, \frak f_{R,\psi})$ 
in $ \Nat (\frak f_{Q,\varphi})$ is a complement of the image of~$\Z(\frak f_{Q,\varphi})\,.$\/}

\medskip
\noindent
{\bf Proof:} By the very definitions of $\Z (\frak f_{Q,\varphi})$ and of $\F_\frak H\,,$ note that the condition 
on $\frak f_{Q,\varphi}$ is equivalent to
$$\big\vert C_P\big(\varphi' (Q)\big)\big\vert \le \big\vert C_P\big(\varphi (Q)\big)\big\vert
\eqno £3.7.2\phantom{.}$$
for any $\varphi'\in \F_\frak H (P,Q)$ which implies that $\varphi (Q)$ is {\it fully centralized\/} in $\F_{\!\frak H}$ 
[4,~Proposition~2.7]. If $\F_{\!\frak H}$ is a Frobenius $P\-$category then, denoting by $R$ the converse image in $P$ of
$\tilde\F_P (Q)\cap {}^{\varphi^*}\!  \tilde\F_P\big(\varphi (Q)\big)\,,$ we know
that there exists $\psi\in  \F_\frak H (P,R)$ extending $\varphi$ or, equivalently,
such that $\Delta_\psi (R)$ contains $\Delta_\varphi (Q)\,;$ in this case, it is clear that
$\Delta_\varphi (Q)$ is actually normal in $\Delta_\psi (R)\,,$ so that we get
$$N_{P\times P}\big(\Delta_\varphi (Q)\big) = \Big(C_P (Q)\times 
 C_P\big(\varphi (Q)\big)\Big)\. \Delta_\psi (R)
 \eqno £3.7.3,$$
 and therefore the exact sequence~£3.6.2 above is {\it split\/}; moreover, we have an
evident  natural map $\mu\,\colon \frak f_{Q,\varphi}\to \frak f_{R,\psi}$ and 
the stabilizer of $\mu\in \Nat (\frak f_{Q,\varphi}, \frak f_{R,\psi})$ in
$ \Nat (\frak f_{Q,\varphi})$ is a complement of the image of $\Z(\frak f_{Q,\varphi})\,.$
 
 \smallskip
Conversely, assume that  the condition on $\frak f_{Q,\varphi}$ holds
 and that there are an indecomposable functor $\frak f_{T,\eta}$ in $\frak H$ and 
 a natural map $\nu$ from $\frak f_{Q,\varphi}$ to~$\frak f_{T,\eta}$ such~that 
 the stabilizer of $\nu\in \Nat (\frak f_{Q,\varphi}, \frak f_{T,\eta})$ in
$ \Nat (\frak f_{Q,\varphi})$ is a complement of the image of $\Z(\frak f_{Q,\varphi})\,;$ then, with the notation above, we claim that  there exists $\psi\in  \F_\frak H (P,R)$ extending $\varphi\,.$ Indeed, it is clear that $\nu$ is determined
by $(u,v)\in T_{P\times P} \big(\Delta_\varphi (Q),\Delta_\eta (T)\big)$
and therefore the converse image in $N_{P\times P}\big(\Delta_\varphi (Q)\big)$ of 
the stabilizer of $\nu$ coincides with $N_{\Delta_\eta (T)^{(u,v)}}
\big(\Delta_\varphi (Q)\big)$ and covers $\tilde\F_P (Q)\cap {}^{\varphi^*}\!  \tilde\F_P\big(\varphi (Q)\big)\,,$ which forces 
$$N_{\Delta_\eta (T)^{(u,v)}} \big(\Delta_\varphi (Q)\big) = \Delta_\psi (R)\j \Delta_\varphi (Q)
 \eqno £3.7.4\phantom{.}$$
 where $\psi (w) = v\eta (uwu^{-1})v^{-1}$ for any $w\in R\,;$ hence, $\psi$ extends
 $\varphi$ and belongs to $\F_\frak H (P,R)\,,$ proving the claim.
 \eject

 \smallskip
 Now, if $\varphi'\,\colon Q\to P$ is an $\F_\frak H\-$morphism such that $\varphi' (Q)$ is {\it fully centralized\/} in $\F_\frak H$
 [4,~2.6] and $R'$ is a subgroup of $N_P(Q)$ such that 
 $$Q\i R'\qq {}^{\varphi'}\!\F_{R'} (Q)\i \F_P\big(\varphi' (Q)\big)
 \eqno £3.7.5,$$
 we claim that there exists $\psi'\in  \F_\frak H (P,R')$ extending $\varphi'\,.$ Indeed, setting $Q' = \varphi' (Q)\,,$
 it is clear that $\frak f_{Q',\varphi\circ \varphi'^*}$ is an  indecomposable functor  in $\frak H$ which still fulfills the
 condition above; hence, our hypothesis and the argument above proves that, denoting by $T'$ the converse image of 
${}^{\varphi'}\!\F_{R'} (Q)$ in $P\,,$ there exists $\eta'\in \F_\frak H (P,T')$ extending  $\varphi\circ \varphi'^*\,;$
since  $Q'$ is {\it fully centralized\/} in~$\F_\frak H\,,$ in particular we have
$$\eta' \big(C_P (Q')\big) = C_P\big(\varphi (Q)\big)
\eqno £3.7.6;$$
then, since ${}^{\varphi^*}\!\F_{\eta'(T')}\big(\varphi (Q)\big) = \F_{R'}(Q)$ and $C_P\big(\varphi (Q)\big)\i \eta' (T')\,,$ it is clear that we have 
$$R'\i R\qq \psi (R') \i \eta' (T')
\eqno £3.7.7$$ 
and the $\F_\frak H\-$morphism mapping $v\in R'$ on $\eta'^*\big(\psi (v)\big)$ extends $\varphi'\,,$ 
proving our claim. Thus, $\F_\frak H$ is a Frobenius $P\-$category by [4,~Proposition~2.11].

 \medskip
 £3.8. We call {\it Frobenius $\frak H_P\-$category\/} any  divisible $\frak H_P\-$category $\frak H$ such that 
$\F_{\!\frak H}$ is a Frobenius $P\-$category. In this case, we have the following result which --- although we do not 
employ it --- generalizes  [4,~Proposition~6.7]. Recall that, if $\F$ is a {\it divisible $P\-$category\/},
 we call  {\it $\F\-$selfcentralizing\/} any subgroup~$Q$ of $P$ 
such that $C_P \big(\varphi (Q)\big)\i \varphi (Q)$ for any $\varphi\in \F(P,Q)\,;$ note that~[4,~4.8]
\smallskip
\noindent
£3.8.1\quad {\it If $Q$ is $\F\-$selfcentralizing then  any subgroup $R$ of $P$ such that $\F (R,Q)$ is not empty
 is $\F\-$selfcentralizing too.\/}
\smallskip
\noindent
 Moreover, let us say that an 
indecomposable functor $\frak f$ in $\frak H$ is {\it maximal\/} if any natural map from $\frak f$ to
an  indecomposable functor in $\frak H$ is a natural isomorphism.

\bigskip
\noindent
{\bf Proposition~£3.9.} {\it Let $\frak H$ be a Frobenius $\frak H_P\-$category and set $\F = \F_{\!\frak H}\,.$ 
For any indecomposable functor $\frak f_{Q,\varphi}$ in $\frak H$ such that $Q$ is
 {\it $\F\-$selfcentralizing\/} there is a unique isomorphism class of maximal indecomposable functors 
 $\frak f_{R,\psi}$ in $\frak H$  such that $\Nat (\frak f_{Q,\varphi},\frak f_{R,\psi})$ is not empty, and then we have
$$\widetilde\Nat (\frak f_{Q,\varphi},\frak f_{R,\psi})= \{\tilde\mu\}
\eqno £3.9.1.$$
In particular, for any natural map $\mu'\,\colon \frak f_{Q,\varphi}\to \frak f_{R',\psi'}$ there is a unique exterior natural map 
$\tilde\eta\,\colon \frak f_{R',\psi'}\to \frak f_{R,\psi}$ fulfilling $\tilde\eta\circ \tilde\mu' = \tilde\mu\,.$\/}

\medskip\noindent
{\bf Proof:} We argue by induction on $\vert P\,\colon Q\vert$ and may assume that $Q\not= P\,;$ denote by 
$N$ the converse image in $P$ of the intersection $\tilde\F_P (Q)\cap  {}^{\varphi^*}\!  \tilde\F_P\big(\varphi (Q)\big)\,;$ 
then, it follows from [4,~£2.12.1] that $\varphi$ can be extended  to an $\F\-$morphism
$\nu\,\colon N\to P\,,$ so that $\Delta_\nu (N)$ contains $\Delta_\varphi (Q)\,;$ in particular, $N$ is  $\F\-$self-centralizing
and we have a natural map $\alpha$ from $\frak f_{Q,\varphi}$ to~$\frak f_{N,\nu}$~(cf.~£2.6).
It~is\break
\eject
\noindent
 clear that $\frak f_{N,\nu}$ admits a natural map $\beta$ to a maximal  indecomposable functor
$\frak f_{R,\psi}$ in $\frak H\,,$ so that $R$ is also $\F\-$selfcentralizing and we actually may assume that 
$\Delta_\nu (N)\i \Delta_\psi (R)\,;$ denote by $\mu$ the  composition of the natural maps
$$\frak f_{Q,\varphi}\buildrel\alpha\over\too \frak f_{N,\nu}\buildrel\beta\over\too \frak f_{R,\psi}
\eqno £3.9.3.$$

\smallskip
If  $\mu'\,\colon \frak f_{Q,\varphi}\to \frak f_{R',\psi'}$ is a natural map, there are $u,v\in P$
such that (cf.~£2.11.2)
$$Q\i R'^u\qq \varphi = \kappa_v\circ \psi'\circ \iota_{Q^{u^{-1}}}^{R'}\circ \kappa_u
\eqno £3.9.4\phantom{.}$$
where $\kappa_{u}\,\colon Q\to Q^{u^{-1}}$ and $\kappa_{v}\,\colon P\to P$ respectively  denote the
conjugation by~$u$ and by $v\,,$ and therefore we get
$$\tilde\F_{R'^u}(Q)\i \tilde\F_P (Q)\cap {}^{\varphi^*}\!  \tilde\F_P\big(\varphi (Q)\big)
\eqno £3.9.5;$$
moreover, setting $N' = N_{R'}(Q^{u^{-1}})$ and denoting by $\nu'\,\colon N'\to P$ the restriction of $\psi'$
and by $\kappa'_u\,\colon N'^u\to N'$ the conjugation by $u\,,$ we still have $\varphi = \kappa_v\circ\nu'\circ\kappa'_u\circ \iota_Q^{N'^u}$ and $\mu'$ factorizes throughout obvious natural maps
$$\frak f_{Q,\varphi}\buildrel\alpha'\over\too\frak f_{N',\nu'}\buildrel\beta'\over\too \frak f_{R',\psi'}
\eqno £3.9.6.$$

\smallskip
In particular,  both group homomorphisms $\nu\circ\iota_{N'}^N$ and $\kappa_v\circ\nu'\circ\kappa'_u$
from $N'$ to $P$ extend $\varphi\,;$ consequently, it follows from [4,~Proposition~4.6] that there exists $z\in C_P (Q)$
such that $\nu\circ \iota_{N'}^N = \kappa_v\circ\nu'\circ \kappa'_{uz}$ and therefore we get a natural 
map $\gamma\,\colon \frak f_{N',\nu'}\to \frak f_{N,\nu}$ determined by the $P\times P\-$set map
$$P_{\nu'}\times_{N'} P\too P_{\nu}\times_N P
\eqno £3.9.7\phantom{.}$$
sending  the class of $(1,1)$ to the class of $(v^{-1},uz)\,;$ then, since $\alpha$ and $\alpha'$ are respectively determined by
the $P\times P\-$set map
$$P_\varphi \times_Q P\too P_{\nu}\times_N P\qq P_\varphi \times_Q P\too P_{\nu'}\times_{N'} P
\eqno £3.9.8\phantom{.}$$
sending the class of $(1,1)$ to the classes of $(1,1)$ and $(v^{-1},u)\,,$ we finally obtain~$\tilde\gamma\circ \tilde\alpha' 
= \tilde\alpha\,.$

\smallskip
Now, since $Q\i N'^u\,,$ $N'$ is $\F\-$selfcentralizing; thus,
considering the natural maps $\beta\circ\gamma\,\colon \frak f_{N',\nu'}\to \frak f_{R,\psi}$ and 
$\beta'\,\colon \frak f_{N',\nu'}\too \frak f_{R',\psi'}\,,$ it follows from the induction hypothesis that
 there exists a unique exterior natural map 
$$\tilde\eta : \frak f_{R',\psi'}\too \frak f_{R,\psi}
\eqno £3.9.9\phantom{.}$$
 fulfilling $\tilde\eta\circ \tilde\beta' = \tilde\beta\circ\tilde\gamma$ and therefore we finally obtain
 $$\tilde\eta\circ \tilde\mu'  = \tilde\eta\circ \tilde\beta'\circ\tilde\alpha' = \tilde\beta\circ\tilde\gamma\circ\tilde\alpha' 
 = \tilde\beta\circ\tilde\alpha = \tilde\mu
 \eqno £3.9.10.$$
On the other hand, if $\tilde\theta \,\colon\frak f_{R',\psi'}\to \frak f_{R,\psi}$ is another natural map such that
$\tilde\theta\circ \tilde\mu' = \tilde\mu$ then, applying  [4,~Proposition~4.6] as above, we get
$\tilde\theta\circ \tilde\beta' = \tilde\beta\circ\tilde\gamma$ and therefore we still get $\tilde\theta = \tilde\eta\,.$
We are done.
\eject

 \bigskip
 \noindent
 {\bf £4. The Alperin theorem in a Frobenius $\frak H_P\-$category} 
\bigskip
£4.1. Let $\F$ be a {\it divisible $P\-$category\/} and set $\frak H = \frak H_\F\,;$ in [4, Lemma~5.4] we state
an {\it additive\/} version of Alperin's original formulation of his celebrated  theorem [1], which in $\F$ is expressed {\it via\/}
equalities [4,~5.4.2]; but, as a matter of fact, the corresponding equalities in $\tilde\F$ admit an equivalent formulation 
in the category $\frak H$  as our next result shows.

\bigskip
\noindent
{\bf Proposition~£4.2.} {\it Let $\frak f$ be the composition of a family $\frak F = \{\frak f_{Q_i,\varphi_i}\}_{i\in \Delta_n}$~of indecomposable functors in $\frak H\,.$ For any indecomposable functor $\frak f_{Q,\varphi}$ in~$\frak H$ 
there is a natural map  $\frak f_{Q,\varphi}\to \frak f$ if and only if for any $i\in \Delta_n$ there is an $\tilde\F\-$morphism 
$\tilde\psi_i\,\colon Q\to Q_i$ fulfilling
$$\eqalign{\tilde\iota_Q^P &= \tilde\iota_{Q_0}^P\circ \tilde\psi_0\cr
\tilde\varphi_{i-1}\circ \tilde\psi_{i-1} &=  \tilde\iota_{Q_i}^P\circ \tilde\psi_i \hbox{ for any $1\le i\le n$}\cr
\tilde\varphi_n\circ \tilde\psi_n &= \tilde\varphi\cr}
\eqno £4.2.1.$$\/}
\medskip
\noindent
{\bf Proof:} By the composition of $\frak F$ we mean that we inductively define a new family 
$\{\frak g_i\}_{i\in \Delta_n}$ of functors setting $\frak g_0 =\frak f_{Q_0,\varphi_0}$ and 
$\frak g_i = \frak f_{Q_i,\varphi_i}\circ \frak g_{i-1}$ for any~$1\le i\le n\,;$ in order to compute $\frak f =\frak g_n\,,$ the first step is the choice of a set of representatives $W_{R',R''}^{R}\i R$ for the set of double classes  $R'\backslash R/R''$ when $R$ runs over the set of subgroups of $P$ and $(R',R'')$
over the set of pairs of subgroups of $R\,.$

\smallskip
 Then, we consider the set $\W(\frak F)$ of finite double sequences   of elements  of~$P$  
 $$w = \{w_{j +1,i}\}_{j,i\in \Delta_n,\, j +i  < n}
 \eqno £4.2.2\phantom{.}$$
such that, setting $w(0,i) = \emptyset$  for any $i\in \Delta_n$ and considering
$$w(j +1,i) =  \{w_{\ell +1,i +k}\}_{\ell, k\in \Delta_j,\, \ell +k\le j}
\eqno £4.2.3\phantom{.}$$
 for any $j,i\in \Delta_n$ fulfilling $j +i < n\,,$ and inductively defining $R_{w(0,i)} = Q_i$ for any $i\in \Delta_n$ and
$$R_{w(j +1,i)} = (\varphi_i)^{-1}\big((R_{w(j,i+1)})^{w_{j +1,i }}
\cap  \varphi_i (R_{w(j,i)})\big)
\eqno £4.2.4\phantom{.}$$
 for any $j,i\in \Delta_n$ fulfilling $j +i < n\,,$ the element $(w_{j +1,i})^{-1}$ belongs to the chosen set of representatives 
 $$W_{R_{w(j,i+1)},\varphi_i (R_{w(j,i)})}^{R_{w(j-1,i +1)}}\cong R_{w(j,i+1)} \backslash R_{w(j-1,i +1)}/\varphi_i (R_{w(j,i)})
\eqno £4.2.5\phantom{.}$$
 for any $j,i\in \Delta_n$ fulfilling $j +i < n\,,$ where we still set $R_{w(-1,i +1)} = P\,.$  For any $j,i\in \Delta_n$ 
 fulfilling $j +i < n\,,$ we denote by $\kappa_{w_{j +1,i}}\,\colon P\cong P$ the conjugation by~$w_{j +1,i}\,.$

\smallskip
For such a sequence, note that $R_{w(j +1,i)}$ is contained in $R_{w(j,i)}$ and that
$\kappa_{w_{j +1,i}}\circ \varphi_i$ induces a group homomorphism  
$\varphi_{j,i}\,\colon R_{w(j +1,i)}\to R_{w(j,i+1)}\,;$ then, we have an evident commutative $\tilde\F\-$diagram
$$\matrix{&&R_{w(j-1,i +1)}\cr
&{\tilde\varphi{j-1,i}\atop}\hskip-5pt\nearrow\hskip-5pt&&\hskip-15pt\nwarrow\cr
R_{w(j,i)}\hskip-35pt&&&&\hskip-10pt R_{w(j,i+1)}\cr
&\nwarrow\hskip-25pt&&\hskip-25pt{\tilde\varphi{j,i}\atop}\hskip-5pt\nearrow\cr
&&R_{w(j +1,i)}\cr}
\eqno £4.2.6,$$
where the left-hand oriented $\tilde\F\-$morphisms are induced by the corresponding inclusions;
moreover, we inductively define 
$$\omega_{w(j,i)} : R_{w(j ,i)}\too Q_{j +i}
\eqno £4.2.7\phantom{.}$$
setting $\omega_{w(0,i)} = {\rm id}_{Q_i}$ and $\omega_{w(j +1,i)} 
= \omega_{w(j,i+1)}\circ \varphi_{j,i}\,,$ and from the commutativity of diagram~£4.2.6
above it is not difficult to get the following commutative $\tilde\F\-$diagram
$$\matrix{&&&\hskip-25pt P\cr
&&\hskip-15pt{\tilde\varphi_{j +i}\atop}\hskip-5pt\nearrow\hskip-5pt&&\hskip-50pt\nwarrow\cr
&&\hskip-10ptQ_{j +i}&&Q_{j +i +1}\cr
&{\tilde\omega {j,i}\atop}\hskip-5pt\nearrow\hskip-5pt\cr
R_{w(j,i)}\hskip-35pt&&&\hskip-20pt {\tilde\omega_{j +1,i}\atop}\hskip-5pt \nearrow\cr
&\nwarrow\hskip-25pt\cr
&&R_{w(j +1,i)}\cr}
\eqno £4.2.8.$$

\smallskip
Now, we claim that
$$\frak f = \frak g_n = \bigsqcup_{w\in \W(\frak F)} 
\frak f_{R_{w(n,0)},\,\varphi_n\circ\omega_{w(n,0)}}
\eqno £4.2.9;$$
 indeed, arguing by induction on $n\,,$ we may assume that $n\ge 1\,;$
then, by the induction hypothesis, we already have
$$\frak g_{n-1} = \bigsqcup_{w'\in \W(\frak F')} \frak f_{R_{w'(n-1,0)},\,\varphi_{n-1}\circ\omega_{w'(n-1,0)}}
\eqno £4.2.10\phantom{.}$$
where we set $\frak F' = \{\frak f_{Q_i,\varphi_i}\}_{i\in \Delta_{n-1}}\,;$ hence, we get
$$\frak g_n =  \bigsqcup_{w'\in \W(\frak F')} \frak f_{Q_n,\varphi_n}\circ \frak f_{R_{w'(n-1,0)},\,\varphi_{n-1}
\circ\omega_{w'(n-1,0)}}
\eqno £4.2.11.$$ 
But, it is easily checked that we have a surjective map $\pi\,\colon \W(\frak F)\to \W(\frak F')$ sending
$w\in \W(\frak F)$ to the double sequence $w(n -1,0)$ and that, for any $w'\in \W(\frak F')\,,$ the set of ordered products
$$W =\Big\{\prod_{j\in \Delta_{n-1}} (w_{j +1,n -j -1})^{-1}\Big\}_{w\in \pi^{-1}(w')}
\eqno £4.2.12\phantom{.}$$
\eject
\noindent
is a set of representatives for the set of double classes 
$$Q_n\backslash P/(\varphi_{n-1}\circ \omega_{w' (n-1,0)})(R_{w'(n-1,0)})
\eqno £4.2.13.$$
Consequently, from equality~£3.1.1 we have
$$\frak f_{Q_n,\varphi_n}\circ \frak f_{R_{w'(n-1,0)},\,\varphi_{n-1}\circ\omega_{w'(n-1,0)}} = \bigsqcup_{w\in \pi^{-1}(w')} \frak f_{R_{w(n,0)},\,\varphi_n\circ\omega_{w(n,0)}}
\eqno £4.2.14,$$
which proves our claim.

\smallskip
In particular, if there is a natural map $\frak f_{Q,\varphi}\to \frak f$ from an indecomposable functor  $\frak f_{Q,\varphi}$
in $\frak H\,,$ then there are $w\in \W(\frak F)$ and a natural map
$$\frak f_{Q,\varphi}\too  \frak f_{R_{w(n,0)},\, \varphi_n\circ\omega_{w(n,0)}}
\eqno £4.2.15,$$ 
which amounts to saying that there are $u,v\in P$ such that (cf.~£2.11)
$$Q\i (R_{w(n,0)})^{u}\qq \varphi = \kappa_{v}\circ \varphi_n\circ
\omega_{w(n,0)}\circ \kappa_{u}
\eqno £4.2.16\phantom{.}$$
where $\kappa_u\,\colon Q\to R_{w(n,0)}$ and $\kappa_v\,\colon P\to P$ respectively denote the corresponding 
conjugation by $u$ and by $v\,;$ then, setting
$$\psi_j = \omega_{w(j,0)}\circ \iota_{R_{w(n,0)}}^{R_{w(j,0)}}\circ \kappa_u
\eqno £4.2.17\phantom{.}$$
 for any $j\in \Delta_n\,,$ we get the top and the bottom equalities in~£4.2.1, whereas if $j\not= n$ then from the commutativity 
of diagram~£4.2.8 we obtain
$$\eqalign{\tilde\varphi_{j -1}\circ \tilde\psi_{j-1} &= \tilde\varphi_{j -1}\circ
\tilde\omega_{w(j -1,0)}\circ\tilde\iota_{R_{w(j,0)}}^{R_{w(j-1,0)}}\circ
\tilde\iota_{R_{w(n,0)}}^{R_{w(j,0)}}\circ \tilde\kappa_u\cr
&= \tilde\iota_{Q_j}^P\circ \tilde\omega_{w(j,0)}\circ 
\tilde\iota_{R_{w(n,0)}}^{R_{w(j,0)}}\circ \tilde\kappa_u = 
\tilde\iota_{Q_j}^P\circ \tilde\psi_j\cr}
\eqno £4.2.18,$$
which proves the middle equalities in~£4.2.1.

\smallskip
Conversely, assume that the equalities~£4.2.1 hold; arguing by induction on~$n\,,$ we may assume that $n\not= 0\,;$ 
then, setting $\varphi' = \iota_{Q_n}\circ\psi_n\,,$ the induction hypothesis guarantees the existence of a
natural map $\mu'\,\colon \frak f_{Q,\varphi'}\to\frak g_{n-1}$ and therefore we get the natural map
$$\frak f_{Q_n,\varphi_n} *\mu'  : \frak f_{Q_n,\varphi_n}\circ \frak f_{Q,\varphi'}\too
\frak f_{Q_n,\varphi_n}\circ \frak g_{n-1} = \frak f
\eqno £4.2.19;$$
but, according to equality~£3.1.1, for a suitable functor $\frak f'$ in $\frak H$ we have
$$\frak f_{Q_n,\varphi_n}\circ \frak f_{Q,\varphi'} = \frak f_{Q,\varphi_n\circ\psi_n}  \sqcup \frak f'
\eqno £4.2.20;$$
moreover, since $\tilde\varphi_n\circ \tilde\psi_n = \tilde\varphi\,,$ we have 
$ \frak f_{Q,\varphi_n\circ\psi_n}\cong \frak f_{Q,\varphi}$ and therefore we still have a natural map
$\nu\,\colon \frak f_{Q,\varphi}\to \frak f_{Q_n,\varphi_n}\circ \frak f_{Q,\varphi'}\,;$ finally, we get the natural map
$$\mu = (\frak f_{Q_n,\varphi_n} *\mu' )\circ \nu : \frak f_{Q,\varphi}\too \frak f
\eqno £4.2.21.$$

\medskip
£4.3. Translating here the terminology of [4,~Ch.~5], we say that an indecomposable functor $\frak f_{Q,\varphi}$ in $\frak H$
is {\it $\frak H\-$reducible\/} if either $Q = P$ or there is a natural map from $\frak f_{Q,\varphi}$ to the composition of a family $\frak F = \{\frak f_{Q_i,\varphi_i}\}_{i\in \Delta_n}$~of indecomposable functors in~$\frak H$ such that we have $\vert Q\vert < \vert Q_i\vert$ for any $i\in \Delta_n\,.$ From equality~£3.1.1 it is easily checked that 
the composition of $\frak H\-$reducible indecomposable functors is a disjoint union of $\frak H\-$reducible indecomposable functors.
It follows easily from Proposition~£4.2 that the indecomposable functor $\frak f_{Q,\varphi}$
is {\it $\frak H\-$reducible\/} if and only if,  with the notation and the terminology in [4,~Ch.~5], the {\it Alperin $\F\-$fusion
$\varphi - \iota_Q^P$ in the $\Bbb Z\-$free module over $\F (P,Q)$ is $\F\-$reducible.\/}

\medskip
£4.4. Coherently, we call {\it $\frak H\-$irreducible\/} any indecomposable
functor in $\frak H$ which is {\it not\/} $\frak H\-$reducible. As in [3,~Ch.~III], we say that two {\it $\frak H\-$irreducible\/} 
functors $\frak f_{Q,\varphi}$ and $\frak f_{Q',\varphi'}$ are {\it exchangeable\/} if there are two natural maps
$$\frak f_{Q',\varphi'}\too \frak f\circ \frak f_{Q,\varphi}\circ \frak g\qq
\frak f_{Q,\varphi}\too \frak f'\circ \frak f_{Q',\varphi'}\circ \frak g'
\eqno £4.3.1\phantom{.}$$
for suitable $\frak H\-$reducible indecomposable functors $\frak f\,,$ $\frak g\,,$ $\frak f'$ and $\frak g'\,;$
it is clear that this relation defines an {\it equivalence\/} in the set of $\frak H\-$irreducible indecomposable functors.
Actually,  the existence of any natural map in £4.3.1 forces the existence of the other one, as our next result shows.

\bigskip
\noindent
{\bf Proposition~£4.5.\/} {\it  Two $\frak H\-$irreducible  functors $\frak f_{Q,\varphi}$ and $\frak f_{Q',\varphi'}$ 
are {\it exchangeable\/} if and only if there exists a group isomorphism $\theta\,\colon Q\cong Q'$ such that
the functors $\frak f_{\varphi (Q), \varphi'\circ\theta\circ\varphi^*}$ and $\frak f_{Q',\iota_Q^P\circ\theta^{-1}}$
are $\frak H\-$reducible.\/}

\medskip
\noindent
{\bf Proof:} If we have a natural map from $\frak f_{Q',\varphi'}$ to $\frak f_{R,\psi}\circ \frak f_{Q,\varphi}\circ \frak f_{R',\psi'}$ for suitable  $\frak H\-$reducible indecomposable functors $\frak f_{R,\psi}$ and $\frak f_{R',\psi'}\,,$ the 
{\it $\frak H\-$irreducibility\/} of~$\frak f_{Q',\varphi'}$ forces $\vert Q'\vert = \vert Q\vert\,;$ moreover,
 it follows from equality~£3.1.1 applied twice that there are $w,w'\in P$ such that, setting
$$\eqalign{Q'' &= \psi'^{-1}\bigg(\Big(\varphi^{-1}\big(R^w\cap \varphi(Q)\big)\Big)^{w'}\cap \psi'(R')\bigg)\cr
\varphi'' &= \psi\circ \kappa_w\circ \varphi_w\circ \kappa_{w'}\circ \psi'_{w'}\cr}
\eqno £4.5.1\phantom{.}$$
where $\varphi_w$ and $\psi'_{w'}$ are suitable restrictions of $\varphi$ and $\psi'\,,$ and~$\kappa_w$
and $\kappa_{w'}$ are induced by the conjugation by $w$ and by $w'$ respectively, 
we still have a natural map from $\frak f_{Q',\varphi'}$ to $\frak f_{Q'',\varphi''}\,,$ which actually has to be a natural isomorphism
since $\frak f_{Q',\varphi'}$ is $\frak H\-$irreducible.

\smallskip
Consequently, according to £2.11,  there are $v,v'\in P$ fufilling  $Q' = (Q'')^v$ and 
$\varphi'  = \kappa_{v'}\circ\varphi''\circ\kappa_v\,;$ in particular, we get 
$$\varphi (Q)\i R^w\quad,\quad Q^{w'}\i \psi' (R')\qq Q'' = \psi'^{-1}(Q^{w'})
\eqno £4.5.2\phantom{.}$$
and we denote by $\theta\,\colon Q\cong Q'$ the group isomorphism fulfilling 
$$\psi'\big(\theta (u)^{v^{-1}}\big) = u^{w'}
\eqno £4.5.3\phantom{.}$$ 
for any $u\in Q\,.$ 
\eject

\smallskip
Then, we have natural maps
$$\frak f_{\varphi (Q), \varphi'\circ\theta\circ\varphi^*}\too \frak f_{R,\psi}\qq 
\frak f_{Q',\iota_Q^P\circ\theta^{-1}}\too \frak f_{R',\psi'}\eqno £4.5.4;$$
indeed, we already have $\varphi (Q)\i R^w$ and, moreover, we can write
$$\eqalign{\varphi'\circ\theta\circ\varphi^* &= (\kappa_{v'}\circ\varphi''\circ\kappa_v)\circ\theta\circ\varphi^*\cr
& = \kappa_{v'}\circ (\psi\circ \kappa_w\circ \varphi_w\circ \kappa_{w'}\circ \psi'_{w'})\circ\kappa_v\circ\theta\circ\varphi^*\cr
& = \kappa_{v'}\circ \psi\circ \kappa_w\cr}
\eqno £4.5.5\phantom{.}$$
which proves the existence of the left-hand natural map; similarly, we still have $Q'\i R'^v$ and we can write
$$\theta^{-1} = \kappa_{w'}\circ \psi'_{w'}\circ \kappa_v
\eqno £4.5.6\phantom{.}$$
which proves the existence of the right-hand natural map. Hence, the functors $\frak f_{\varphi (Q), \varphi'\circ\theta\circ\varphi^*}$ and $\frak f_{Q',\iota_Q^P\circ\theta^{-1}}$
are $\frak H\-$reducible too.

\smallskip
Conversely, assume that there is a group isomorphism $\theta\,\colon Q\cong Q'$ such that
the functors $\frak f_{\varphi (Q), \varphi'\circ\theta\circ\varphi^*}$ and $\frak f_{Q',\iota_Q^P\circ\theta^{-1}}$
are $\frak H\-$reducible; on the one hand, from equality~£3.1.1 it is easily checked that there exists a natural map
$$\frak f_{Q',\varphi'}\too \frak f_{\varphi (Q), \varphi'\circ\theta\circ\varphi^*}\circ \frak f_{Q,\varphi}\circ
\frak f_{Q',\iota_Q^P\circ\theta^{-1}}
\eqno £4.5.7;$$
on the other hand, it is quite clear that we have 
$$(\frak f_{\varphi (Q), \varphi'\circ\theta\circ\varphi^*})^\circ = \frak f_{\varphi' (Q'), \varphi\circ\theta^{-1}\circ\varphi'^*}
\qq (\frak f_{Q',\iota_Q^P\circ\theta^{-1}})^\circ = \frak f_{Q,\iota_{Q'}^P\circ\theta}
\eqno £4.5.8,$$
so that $\frak f_{\varphi' (Q'), \varphi\circ\theta^{-1}\circ\varphi'^*}$ and $\frak f_{Q,\iota_{Q'}^P\circ\theta}$
are $\frak H\-$reducible too, and once again we have a natural map
$$\frak f_{Q,\varphi}\too \frak f_{\varphi' (Q'), \varphi\circ\theta^{-1}\circ\varphi'^*}\circ \frak f_{Q',\varphi'}\circ
\frak f_{Q,\iota_{Q'}^P\circ\theta}
\eqno £4.5.9.$$
We are done.

\bigskip
\noindent
{\bf Proposition~£4.6.} {\it Let $\E$ be  a set of representatives for the exchangeability classes
of $\frak H\-$irreducible indecomposable functors. Then, any indecomposable functor $\frak f_{Q,\varphi}$ in $\frak H$ 
admits a natural map to the composition of $\,\frak f_{P,\sigma}$ for a suitable $\sigma\in \F (P)\,,$ with the composition 
of a family $\frak F = \{\frak f_{Q_i,\varphi_i}\}_{i\in \Delta_n}$ of indecomposable functors in  $\E\,.$\/}

\medskip
\noindent
{\bf Proof:} According to Lemma~£4.7 below, it suffices to prove that $\frak f_{Q,\varphi}$ in $\frak H$ admits a natural map
to the composition of a finite family $\frak F = \{\frak f_{Q_i,\varphi_i}\}_{i\in \Delta_n}$  of indecomposable functors in
$$\E' = \{\frak f_{P,\sigma}\}_{\sigma\in \F (P)}\cup\E
\eqno £4.6.1;$$
\eject
\noindent
we argue by induction on $\vert P\,\colon Q\vert$ and may assume that $Q\not= P\,;$ we also may 	assume that 
 $\frak f_{Q,\varphi}$ does not belong to~$\E\,.$ If $\frak f_{Q,\varphi}$ is $\frak H\-$reducible then it admits
 a natural map to the composition of a family $\{\frak f_{R_j,\psi_j}\}_{j\in \Delta_m}$ of indecomposable functors 
 in $\frak H$ fulfilling $\vert Q\vert < \vert R_j\vert$ for any $j\in \Delta_m$ and it suffices to apply the induction hypothesis
 to any $\frak f_{R_j,\psi_j}\,.$ If $\frak f_{Q,\varphi}$ is $\frak H\-$irreducible then we have a natural map
 $$\frak f_{Q,\varphi}\too \frak f\circ \frak f_{Q',\varphi'}\circ \frak g
 \eqno £4.6.2\phantom{.}$$
where $\frak f$ and $\frak g$ are $\frak H\-$reducible indecomposable functors and $\frak f_{Q',\varphi'}$ belongs to~$\E\,;$ 
then, it suffices to apply the precedent argument to both $\frak f$ and $\frak g\,.$

\bigskip
\noindent
{\bf Lemma~£4.7.} {\it For any $\sigma\in {\rm Aut} (P)$ and any indecomposable functor $\frak f_{Q,\varphi}$
 we have
 $$\frak f_{P,\sigma}\circ \frak f_{Q,\varphi} = \frak f_{\sigma (Q),\sigma\circ\varphi\circ\sigma_Q^{-1}}\circ \frak f_{P,\sigma}
 \eqno £4.7.1\phantom{.}$$
 where $\sigma_Q\,\colon Q\cong \sigma (Q)$ denotes the isomorphism induced by $\sigma\,.$\/}

 \medskip
 \noindent
 {\bf Proof:} Straightforward.

 \medskip
 £4.8. As stated in~£4.3 above, an indecomposable functor $\frak f_{Q,\varphi}$ is {\it $\frak H\-$irredu-cible\/} if and only if,  with the notation and the terminology in [4,~5.7], the {\it Alperin $\F\-$fusion\/} $\varphi - \iota_Q^P$ in the $\Bbb Z\-$free 
 module over $\F (P,Q)$ is {\it $\F\-$irreducible,\/} so that $Q$ is  {\it $\F\-$essential\/}. Recall that [4,~Theorem~5.11]
\smallskip
\noindent
£4.8.1\quad {\it  If $\F$ is a Frobenius $P\-$category, $Q$ is  $\F\-$essential if and only if it is 
 $\F\-$selfcentralizing and $\tilde\F (Q)$ admits a proper subgroup $\tilde M_\F$ such that $p$ divides $\vert \tilde M_\F\vert$ 
 and does not divide $\vert \tilde M_\F\cap\tilde M_\F^{\tilde\tau}\vert$ for any $\tilde\tau\in \tilde\F (Q) -\tilde M_\F\,.$\/}
 \smallskip
\noindent
In this case, $\tilde M_\F$ contains a Sylow $p\-$subgroup of $\tilde\F (Q)$ and, according to [4,~5.12.1], these subgroups are exactly the proper subgroups of $\tilde\F (Q)$ containing a conjugate of one of them $\tilde M_\F(P,Q)$ which contains
$\tilde\F_P (Q)\,.$ The following result completes our equivalent formulation.

 \bigskip
 \noindent
 {\bf Theorem~£4.9.} {\it Assume that  $\F$ is a Frobenius $P\-$category. An indecomposable functor 
 $\frak f_{Q,\varphi}$ is $\frak H\-$irreducible if and only if $Q$ is $\F\-$essential and $\tilde M_\F(P,Q)$ 
 does not contain ${}^{\varphi^*}\!\tilde\F_{\!P}\big(\varphi (Q)\big)\,.$ In this case,  for any 
 $\tilde\tau\in \tilde\F (Q) -\tilde M_\F (P,Q)\,,$ the $\frak H\-$irreducible functors $\frak f_{Q,\varphi}$ 
 and~$\frak f_{Q,\iota_Q^P\circ\tau}$ are exchangeable if and only if $\tilde M_\F(P,Q)^{\tilde\tau}$ contains  
 ${}^{\varphi^*}\!\tilde\F_{\!P}\big(\varphi (Q)\big)\,.$\/}

\medskip
\noindent
{\bf Proof:} First of all, we claim that if $Q$ is not $\F\-$selfcentralizing then  $\frak f_{Q,\varphi}$ is $\frak H\-$reducible.
Indeed, let $\theta\,\colon \varphi (Q)\to P$ be an $\F\-$morphism such that $\theta\big(\varphi(Q)\big)$
is fully normalized in $\F$ [4,~Proposition~2.7]; then, it follows from [4,~Proposition~2.11] that, up to modifying
our choice of $\theta\,,$ there is an $\F\-$morphism $\zeta\,\colon N_P (Q)\to P$ extending $\theta\circ\varphi_*$ where
we denote by $\varphi_*\,\colon Q\cong \varphi(Q) = Q'$ the isomorphism determined by $\varphi\,;$
similarly, there exists an $\F\-$morphism $\xi\,\colon N_P (Q')\to P$ which extends $\theta\circ\tau'$ for a suitable
$\F\-$automorphism $\tau'$ of~$Q'\,;$ finally, since $Q''= \theta (Q')$ is also fully centralized in $\F$ by  [4,~Proposition~2.11], 
we have an $\F\-$automorphism $\sigma''$ of  $C'' = Q''\.C_P (Q'')$ extending $\theta_*\circ\tau'\circ\theta_*^{-1}\,.$
\eject

\smallskip
 At this point, since we clearly have the equalities 
$$\eqalign{\iota_Q^P &= \iota_{N_P(Q)}^P\circ \iota_Q^{N_P(Q)}\cr
\zeta\circ\iota_{N_P(Q)}^P &= \iota_{C''}^P\circ \iota_{Q''}^{C''}\circ\theta_*\circ\varphi_*\cr
\iota_{C''}^P\circ\sigma''\circ \iota_{Q''}^{C''}\circ\theta_*\circ\varphi_* &= \xi\circ\iota_{Q'}^{N_P(Q')}\circ\varphi_*\cr
\iota_{N_P(Q')}^P\circ\iota_{Q'}^{N_P(Q')}\circ\varphi_* &= \varphi\cr}
\eqno £4.9.1,$$
it follows from Proposition~£4.2 that we have a natural map
$$\frak f_{Q,\varphi}\too (\frak f_{N_P(Q'),\xi})^\circ \circ \frak f_{C'',\iota_{C''}^P\circ\sigma''}\circ \frak f_{N_P(Q),\zeta}
\eqno £4.9.2,$$
which proves that $\frak f_{Q,\varphi}$ is $\frak H\-$reducible since $Q''\not= C''\,.$

\smallskip
Thus, we may assume that $Q$ is $\F\-$selfcentralizing. Now, we claim that if $\frak f_{Q,\varphi}$ is $\frak H\-$reducible
 then $\tilde M_\F(P,Q)$ contains ${}^{\varphi^*}\!\tilde\F_{\!P}\big(\varphi (Q)\big)\,;$ we may assume that
 $Q\not= P$ and then it follows from Proposition~£4.2 that there are a family $\frak F 
 = \{\frak f_{Q_i,\varphi_i}\}_{i\in \Delta_n}$ of indecomposable functors in $\frak H$ and, for any $i\in \Delta_n\,,$
  a  non-bijective $\tilde\F\-$morphism
$\tilde\psi_i\,\colon Q\to Q_i$ fulfilling
$$\eqalign{\tilde\iota_Q^P &= \tilde\iota_{Q_0}^P\circ \tilde\psi_0\cr
\tilde\varphi_{i-1}\circ \tilde\psi_{i-1} &=  \tilde\iota_{Q_i}^P\circ \tilde\psi_i \hbox{ for any $1\le i\le n$}\cr
\tilde\varphi_n\circ \tilde\psi_n &= \tilde\varphi\cr}
\eqno £4.9.3;$$
 for any 
$i\in \Delta_n\,,$ consider the subgroups of $\tilde\F (Q)$ 
$$\tilde R_i = {}^{(\psi_i)^*}\!\tilde\F_{Q_i} \big(\psi_i (Q)\big)\qq 
\tilde T_i = {}^{(\psi_i)^*}\!\tilde\F_{P} \big(\psi_i (Q)\big)
\eqno £4.9.4\phantom{.}$$ 
where $\psi_i$ is a representative of $\tilde\psi_i\,,$ and $(\psi_i)^*$ denotes the inverse of the isomorphism $Q\cong \psi_i (Q)$ determined by $\psi_i\,.$

\smallskip
Note that $\tilde R_i\i \tilde T_i$ and, since  $Q$ is $\F\-$selfcentralizing, it is clear that we have $\tilde R_i\not= \{1\}\,;$
moreover, it follows from equalities~£4.9.4 that we have $\tilde\F_P (Q) =\tilde T_0$ and
that $\tilde T_{i+1}$ still contains $\tilde R_i$ for any $i\in \Delta_{n}\,,$ where we set 
$$\tilde T_{n+1} = {}^{\varphi^*}\!\tilde\F_{\!P}\big(\varphi (Q)\big)
\eqno £4.9.5.$$
Now, arguing by induction on $i\,,$ we claim that  $\tilde M_\F(P,Q)$ contains $\tilde T_i\,;$ it is clear that  
$\tilde M_\F(P,Q)$ contains a Sylow $p\-$subgroup of $\tilde\F (Q)$ containing $\tilde T_0\,;$ moreover,
if $\tilde M_\F(P,Q)$  contains $\tilde T_i$ and $\tilde M_\F(P,Q)^{\tilde \tau}$ contains $\tilde T_{i+1}$
then $\tilde R_i$ is contained in the intersection $\tilde M_\F(P,Q)\cap \tilde M_\F(P,Q)^{\tilde \tau}\,,$
which implies that $\tilde\tau$ belongs to $\tilde M_\F(P,Q)\,,$ proving our claim.
\eject

\smallskip
Conversely, we claim that if $\frak f_{Q,\varphi}$ is $\frak H\-$irreducible  then $\tilde M_\F(P,Q)$ does not contain
 ${}^{\varphi^*}\!\tilde\F_{\!P}\big(\varphi (Q)\big)\,;$ indeed,  if $\frak f_{Q,\varphi}$ is $\frak H\-$irreducible
 then, with the notation and the terminology in [4,~5.7], the {\it Alperin $\F\-$fusion\/} $\varphi - \iota_Q^P$ in the 
 $\Bbb Z\-$free  module over $\F (P,Q)$ is {\it $\F\-$irreducible,\/} so that the respective images $\bar\varphi$ and 
 $\bar\iota_Q^P$ of $\varphi$ and $\iota_Q^P$  in $\overline{\F(P,Q)}$ do not coincide; then, it follows from [4,~Theorem~5.11] that there is $\tilde\tau\in \tilde\F (Q) - \tilde M_\F(P,Q)$ such that 
 $$\bar\varphi = \bar\iota_Q^P\.\tilde\tau = \overline{\iota_Q^P\circ \tau}
 \eqno £4.9.6\phantom{.}$$
 where $\tau$ is a representative of $\tilde\tau\,,$ which amounts to saying that the {\it Alperin $\F\-$fusion\/}
 $\varphi\circ\tau^{-1} - \iota_Q^P$ is {\it $\F\-$reducible\/} or, equivalently, that the functor $\frak f_{Q,\varphi\circ\tau^{-1}}$
 is $\frak H\-$reducible. Thus, by the argument above, $\tilde M_\F(P,Q)$ contains 
 $${}^{(\varphi\circ\tau^{-1})_*}\tilde\F_P \big(\varphi (Q)\big) = {}^{\tau\circ\varphi_*}\tilde\F_P \big(\varphi (Q)\big)
 \eqno £4.9.7,$$
so that $\tilde M_\F(P,Q)^{\tilde\tau}$ contains ${}^{\varphi_*}\tilde\F_P \big(\varphi (Q)\big)$ which is a nontrivial
$p\-$group, and therefore 
$\tilde M_\F(P,Q)$ can not contain ${}^{\varphi_*}\tilde\F_P \big(\varphi (Q)\big)$ since $\tilde\tau$ does not belong 
to~$\tilde M_\F(P,Q)\,;$ this proves our claim.

\smallskip
More generally, for any  $\tilde\sigma\in \tilde\F (Q) - \tilde M_\F(P,Q)$ such that $\tilde M_\F(P,Q)^{\tilde\sigma}$ contains ${}^{\varphi_*}\tilde\F_P \big(\varphi (Q)\big)\,,$ on the one hand $\frak f_{Q,\iota_Q^P\circ\sigma}$ is $\frak H\-$irreducible since
$\tilde M_\F(P,Q)$ cannot contain ${}^{\sigma_*}\tilde\F_P (Q)\,,$ and on the other hand, since ${}^{(\varphi\circ\sigma^{-1})_*}\tilde\F_P \big(\varphi (Q)\big)$ is contained in $\tilde M_\F(P,Q)\,,$ $\frak f_{Q,\varphi\circ\sigma^{-1}}$ is $\frak H\-$reducible;
but, it is clear that we have a natural map
$$\frak f_{Q,\varphi} \too \frak f_{Q,\varphi\circ\sigma^{-1}}\circ \frak f_{Q,\iota_Q^P\circ\sigma}\circ\frak f_{P,{\rm id}_P}
\eqno £4.9.8;$$
hence, $\frak f_{Q,\varphi}$ and $\frak f_{Q,\iota_Q^P\circ\sigma}$ are exchangeable.

\smallskip
Conversely assume that, for an element $\tilde\sigma$ in $\tilde\F (Q) - \tilde M_\F(P,Q)\,,$ the $\frak H\-$irreducible
functors $\frak f_{Q,\iota_Q^P\circ\sigma}$ and $\frak f_{Q,\varphi}$  are exchangeable; then, it follows from Proposition~£4.5 that there is an automorphism
$\theta$ of $Q$ such that the functors $\frak f_{Q, \varphi\circ\theta\circ\sigma^{-1}}$ and 
$\frak f_{Q,\iota_Q^P\circ\theta^{-1}}$ are $\frak H\-$reducibles; hence, $\theta$ belongs to $\F (Q)$
and $\tilde M_\F(P,Q)$ contains the groups 
$${}^{\theta}\F_P (Q)\qq {}^{(\varphi\circ\theta\circ\sigma^{-1})_*} \F\big(\varphi(Q)\big) = {}^{\sigma\circ\theta^{-1}\circ\varphi_*} \F\big(\varphi(Q)\big)
\eqno £4.9.9,$$
so that $\theta$ belongs to $\tilde M_\F(P,Q)$ and then $\tilde M_\F(P,Q)^{\tilde\sigma}$ contains 
${}^{\varphi_*} \F\big(\varphi(Q)\big)\,.$ We are done.

 \bigskip
 \noindent
{\bf £5. Associate Burnside algebras} 
\bigskip
£5.1. Let us denote by $\B_P$ the free $\Bbb Z\-$module over the set of isomorphism classes 
of indecomposable $P\-$sets; it is clear that any object in $P\-\Set$ determines a {\it positive\/} 
element of $\B_P$ and that the {\it disjoint union\/} of objects in~$P\-\Set$ corresponds to the sum in $\B_P\,;$
similarly, the {\it inner\/} product of  objects in~$P\-\Set$ induces an {\it associative and distributive product\/} 
in $\B_P\,,$ so that $\B_P$ becomes a $\Bbb Z\-$algebra --- called the {\it Burnside $\Bbb Z\-$algebra\/}
of $P\,.$
\eject

\medskip
£5.2. {\it Mutatis mutandis\/}, we denote by $\B_{P\times P}$ the free $\Bbb Z\-$module over the set 
of isomorphism classes of indecomposable functors in $\frak B_{P\times P}$ --- called {\it indecomposable elements
of $\B_{P\times P}\,.$\/} Once again, the {\it disjoint union\/} and the {\it composition of functors\/} induce on 
$\B_{P\times P}$ an {\it associative $\Bbb Z\-$algebra\/}
structure--- called the {\it double Burnside $\Bbb Z\-$algebra\/} of $P\,.$ Actually, we are mostly interested in the
$\Bbb Z\-$subalgebra $\B_{P\times P}^{^{\rm pp}}$ of $\B_{P\times P}$ determined by the functors preserving 
the {\it projective\/} $P\-$sets; that is to say, according to~£2.10 and~£2.11 above, $P\times P$ acts on the set of
pairs $(Q,\varphi)\,,$ where $Q$ is a subgroup of $P$ and $\varphi\,\colon Q\to P$ a group homomorphism, 
sending $(Q,\varphi)$ to $(Q^u,\kappa_v\circ\varphi\circ\kappa_u)$ for any $u,v\in P\,,$ where 
$\kappa_{u}\,\colon Q^u\to Q$ and $\kappa_{v}\,\colon P\to P$ respectively denote the  conjugation by $u$ and by $v\,,$ 
and, denoting by $\X$ a set of representatives for the set of $P\times P\-$orbits  under this action and by $f_{Q,\varphi}$
the isomorphism class of the functor $\frak f_{Q,\varphi}\,,$ it is clear that
the family $\{f_{Q,\varphi}\}_{(Q,\varphi)\in \X}$ is the set 
of isomorphism classes of indecomposable functors in $\frak B_{P\times P}^{^{\rm pp}}\,.$

\medskip
£5.3. It is clear that $\B_{P\times P}$ acts on the  $\Bbb Z\-$module~$\B_P$ or, equivalently, that $\B_P$ 
becomes a $\B_{P\times P}\-$module; moreover, the functor $\frak m^P\,\colon P\-\Set\to \frak B_{P\times P}$
in~£2.9 above clearly induces a $\Bbb Z\-$algebra homomorphism
$$m^P : \B_P\too \B_{P\times P}
\eqno £5.3.1\phantom{.}$$
mapping any $x\in \B_P$ on the element $m_x$ of $\B_{P\times P}$ acting on $\B_P$ {\it via\/} the
multiplication by $x\,;$  further, the auto-equivalence $\frak t$ of $\frak B_{P\times P}$ in~£2.8 above induces an 
{\it involutive anti-isomorphism\/} $t$ on the  $\Bbb Z\-$algebra $\B_{P\times P}\,,$ acting trivially on the image of~$m^P\,.$
Similarly, the correspondence $\ell$ in~£2.10 above induces a $\Bbb Z\-$algebra homomorphism
$$\ell : \B_{P\times P}^{^{\rm pp}}\too \Bbb Z
\eqno £5.3.2;$$
we call {\it length of $f$\/} the image of any $f\in \B_{P\times P}^{^{\rm pp}}\,.$

\medskip
£5.4. More generally, if $P'$ is another finite $p\-$group, we denote by $\B_{P\times P'}$ the free $\Bbb Z\-$module 
over the set  of isomorphism classes of {\it indecomposable\/} functors in $\frak B_{P\times P'}$ --- called 
{\it indecomposable elements of $\B_{P\times P'}\,;$\/} then, the composition of functors clearly 
induces a $\B_{P\times P}\-$module structure on $\B_{P\times P'}$ and, for a third  finite $p\-$group $P''\,,$ a bilinear map
$$\B_{P\times P'}\times \B_{P'\times P''}\too \B_{P\times P''}
\eqno £5.4.1;$$
moreover, for any group homomorphism
$\alpha\,\colon P'\to P\,,$ it is clear that the functors $\res_\alpha$ and $\ind_\alpha$ in £2.5 above determine
elements ${\rm res}_\alpha$ in~$\B_{P\times P'}$ and~${\rm ind}_\alpha$ in $\B_{P'\times P}\,.$ Let us denote by
$\B_{P\times P'}^+$ the elements of $\B_{P\times P'}$ determined by the functors in $\frak B_{P\times P'}$ or, equivalently,
the elements with {\it non-negative\/} coefficients in the canonical $\Bbb Z\-$basis.

\medskip
£5.5. We are interested in  the following {\it partially defined $\Bbb Q\-$valued scalar product\/} in~$\B_{P\times P'}\,,$ 
defined by the natural maps in  $\frak B_{P\times P'}\,;$  if $f,g\in \B_{P\times P'}$ come from two functors 
$\frak f$ and $\frak g$  in~$\frak B_{P\times P'}\,,$ we set
$$\vert f,g\vert = \vert \Nat(\frak f,\frak g)\vert
\eqno £5.5.1;$$
note that, if  $g'\in \B_{P\times P'}$ still  comes from a functor $\frak g'$  in~$\frak B_{P\times P'}$  
and $\frak f$ is {\it indecomposable\/}, form isomorphism~£2.7.3 we get
$$\vert f,g + g'\vert = \vert f,g\vert + \vert f,g'\vert
\eqno £5.5.2,$$
which shows that, setting $\vert f,- g\vert = -\vert f,g\vert\,,$ our definition can be extended to all the pairs 
$(f,g)$ such that $f$ is an {\it indecomposable element\/} and $g\in \B_{P\times P'}\,.$ Moreover, with the notation above, 
if  $f'\in \B_{P\times P'}$ also  comes from a functor $\frak f'$   in~$\frak B_{P\times P'}\,,$ from isomorphism~£2.2.2 we still get 
$$\vert f +f',g\vert = \vert f,g\vert\vert f',g\vert
\eqno £5.5.3\phantom{.}$$
and therefore our definition can be extended to all the pairs $(f,g)$ such that $f$ comes from a functor $\frak f$ and 
$g$ belongs to $\B_{P\times P'}\,.$ Finally, we define $\vert f,0\vert = 0$ for any $f\in \B_{P\times P'}\,,$  and
$$\vert f - f', g\vert = \vert f,g\vert/\vert f',g\vert
\eqno £5.5.4\phantom{.}$$
for any $g\in \B_{P\times P'}- \{0\}$ and any $f,f'\in \B_{P\times P'}^+$ such that $\vert f',g\vert \not= 0\,.$

\medskip
£5.6. Note that if $g\not= 0$ then $\vert 0,g\vert = 1$ and clearly there is an indecomposable element 
$f$ of $\B_{P\times P'}$ such that  $\vert f,g\vert \not= 0\,.$ Moreover, for any $f,f',g\in \B_{P\times P'}$ 
equality~£5.5.3 holds provided that all the terms are defined; indeed, this is clear if $g = 0\,;$ otherwise, 
there are $h,\bar h,h',\bar h'\in \B_{P\times P'}^+$ fulfilling 
$$f = h - \bar h\quad ,\quad f' = h' - \bar h'\quad ,\quad \vert \bar h,g\vert\not= 0\qq \vert \bar h',g\vert\not= 0
\eqno £5.6.1;$$
then, since $\vert \bar h +\bar h',g\vert = \vert\bar h,g\vert\vert\bar h',g\vert\not= 0\,,$ according to~£5.5.4, we have 
$$\eqalign{\vert f +f',g\vert &= \vert (h + h') - (\bar h + \bar h'),g\vert = \vert h + h',g\vert\big/\vert \bar h + \bar h',g\vert\cr
&= \vert h ,g\vert\vert h' ,g\vert\big/\vert \bar h ,g\vert\vert \bar h' ,g\vert = \big(\vert h ,g\vert\big/\vert \bar h ,g\vert\big)
\big(\vert h' ,g\vert\big/\vert \bar h' ,g\vert\big)\cr
&=\vert f,g\vert\vert f',g\vert\cr}
\eqno £5.6.2.$$

\medskip
£5.7. On the other hand, it follows from isomorphism~£2.3.3 that,  for any group homomorphism $\alpha\,\colon P'\to P\,,$ 
any $f\in \B_{P''\times P'}$ and any $h\in \B_{P''\times P}\,,$ we get
$$\vert {\rm ind}_\alpha\. f,h\vert = \vert f,{\rm res}_\alpha\. h\vert
\eqno £5.7.1\phantom{.}$$
provided that both members are defined; indeed, if $f$ is {\it indecomposable\/} and comes from $\frak f_D$ for some subgroup $D$ of $P'\times P''$ then ${\rm ind}_\alpha\. f$ comes from $\frak f_{(\alpha\times {\rm id}_{P''})(D)}$ and thus it is {\it indecomposable\/} too, so that the corresponding equality~£5.7.1 holds for any $h\in \B_{P''\times P}$ 
(cf.~equality~£5.5.2); moreover, if $f = f' - f''$ for suitable $f',f''\in \B_{P''\times P'}^+\,,$ we have 
$$\vert {\rm ind}_\alpha\. f'',h\vert = \vert f'',{\rm res}_\alpha\. h\vert
\eqno £5.7.2\phantom{.}$$
and whenever both members are not zero, we get
$$\eqalign{\vert {\rm ind}_\alpha\. f,h\vert &= \vert {\rm ind}_\alpha\. f',h\vert/\vert {\rm ind}_\alpha\. f'',h\vert\cr 
&= \vert f',{\rm res}_\alpha\. h\vert/\vert f'',{\rm res}_\alpha\. h\vert = \vert f,{\rm res}_\alpha\. h\vert\cr}
\eqno £5.7.3.$$

\medskip
£5.8. Any {\it divisible $\frak H_P\-$category\/} $\frak H$ introduced in~£3.4 clearly determines
a $\Bbb Z\-$subalgebra $\H$ of $\B_{P\times P}$ which is $t\-$stable --- here the word ``symmetric''
would be misleading --- and contains the $\Bbb Z\-$subalgebra determined by $\frak H_P$ --- noted~$\H_P\,;$
it is clear that $\H$ is a free $\Bbb Z\-$module over the set of {\it indecomposable elements\/} in~$\H\,.$
Let us say that an {\it indecomposable\/} element $f$ is {\it maximal\/}~in~$\H$ if it comes from a {\it maximal\/}
indecomposable functor in $\frak H$ or, equivalently, we have $\vert f,g\vert = 0$ for any
{\it indecomposable element\/} $g$ of~$\H$ different from~$f\,.$ We already know that any {\it divisible 
$\frak H_P\-$category\/} $\frak H$ comes from a  {\it divisible $P\-$category\/}~$\F = \F_\frak H$ (cf.~Proposition~£3.3) 
and we denote by~$\H_\F$ the corresponding $\Bbb Z\-$subalgebra --- called {\it Hecke $\Bbb Z\-$algebra\/} of $\F\, ;$ 
actually, these $\Bbb Z\-$subalgebras can be directly described by rewriting  Proposition~£3.3 in our new context.

  \bigskip
 \noindent
 {\bf Proposition~£5.9.\/} {\it A $\Bbb Z\-$subalgebra $\H$ of $\B_{P\times P}^{^{\rm pp}} $ comes from a 
divisible $P\-$cate-gory if and only if it is $t\-$stable,  contains $\H_P\,,$ and fulfills the following condition
 \smallskip
\noindent
{\rm £5.9.1.}\quad  Any indecomposable element $f\in \B_{P\times P}$ such that $\vert f,\H\vert \not= \{0\}$ belongs to $\H\,.$\/}

\medskip
\noindent
{\bf Proof:} If $\F$ is a divisible $P\-$category then it is clear that $\H_\F$ contains $\H_P$ and is $t\-$stable.
Let $f$ be an indecomposable element of $\B_{P\times P}$ such that  $\vert f,\H_\F\vert \not= \{0\}\,;$
then, $f$ comes from an indecomposable functor $\frak f$ in $\frak H_\F$ and, by the very definition of 
the partial scalar product in $\B_{P\times P}\,,$ there exists a natural map from $\frak f$ to some functor
 in $\frak H_\F$ and therefore, according to statement~£3.2.1, 
$\frak f$ is an object in $\frak H_\F\,,$ so that $f$ belongs to~$\H_\F\,.$

\smallskip
Conversely, let us denote by $\frak H$ the {\it full\/} subcategory of $\frak B_{P\times P}$ over the objects having
their isomorphism class in $\H\,;$ it is clear that $\frak H$ is symmetric, contains~$\frak H_P$ and is closed with respect
to disjoint unions and  composition; moreover, since we assume that $\H\i \B_{P\times P}^{^{\rm pp}}\,,$ all the functors
in $\frak H$ preserve the {\it projective\/} $P\-$sets. Thus, according to Proposition~3.3, it suffices to prove that $\frak H$
fulfills condition~£3.3.1; let $\frak f$ be a functor in $\frak B_{P\times P}$ admitting a natural map to some functor 
in~$\frak H\,;$ once again according to the very definition of the scalar product in~$\B_{P\times P}\,,$ 
denoting by $f$ the isomorphism class of $\frak f\,,$ we have $\vert f,\H\vert \not= \{0\}\,;$ in particular,
for any {\it proper\/} decomposition $f = f' + f''$ such that $f',f''\in \B_{P\times P}^+\,,$ we still have 
$\vert f',\H\vert \not= \{0\} \not= \vert f'',\H\vert\,;$ thus, arguing by induction on $\ell (f)\,,$ condition~£5.9.1
implies that $f$ belongs to $\H$ and therefore $\frak f$ is indeed a functor in $\frak H\,,$
We are done.

\medskip
£5.10. Similarly, most of the statements in section~£4 can be easily translated to the new context; namely,
if $\H$ is a {\it Hecke $\Bbb Z\-$algebra\/}, we say that an indecomposable element $f\in \H$ is {\it $\H\-$irreducible\/} if $\ell (f)\not= 1$ and, for any family 
$\{f_i\}_{i\in \Delta_n}$ of indecomposable elements of $\H$ such that 
$$\vert f,\prod_{i\in \Delta_n} f_i\vert \not= 0
\eqno £5.10.1,$$
we have $\ell (f) = \ell (f_i)$ for some $i\in \Delta_n\,;$ then, two  {\it $\H\-$irreducible\/} elements $f$ and $f'$
are said {\it exchangeable\/} if we have 
$$\vert f',gfh\vert\not= 0\qq \vert f,g'f'h'\vert\not= 0
\eqno £5.10.2\phantom{.}$$
for suitable indecomposable elements $g,h,g',h'\in \H$ {\it not $\H\-$irreducible\/}; at this point , it follows from 
Proposition~£4.6 that
\smallskip
\noindent
£5.10.3\quad  {\it Let $\E$ be a set of representatives for the exchangeability classes
of $\H\-$ irreducible elements. For any indecomposable element $f\in \H$ there is a family $\{f_i\}_{i\in \Delta_n}$ 
of elements in $\E$ and an indecomposable element $g\in \H$ of length one fulfilling
$$\vert f,g\prod_{i\in \Delta_n}f_i\vert\not= 0
\eqno .$$\/}

\medskip
£5.11. We see below that the fact that a $\Bbb Z\-$subalgebra $\H$ of $\B_{P\times P}$ is the {\it Hecke algebra\/}
of some  {\it Frobenius $P\-$category\/}  can be decided from the own structure of $\H\,.$ First of all, we state a new
version of Proposition~£3.9; as above, let us call {\it positive\/} any element of $\H_\F\cap \B_{P\times P}^+\,;$
moreover, if $\F$ is a {\it divisible $P\-$category\/}, let us denote by  $\N_\F$ the $\Bbb Z\-$submodule generated by the elements of $\H_\F$ coming from 
the indecomposable functors $\frak f_{Q,\varphi}$ such that $Q$ is {\it not\/} $\F\-$selfcentralizing.

\bigskip
\noindent
{\bf Lemma~£5.12.} {\it $\N_\F$  is a two-sided ideal of~$\H_\F\,.$\/}

\medskip
\noindent
{\bf Proof:} Since $\H_\F$ and $\N_\F$ are $t\-$stable, it suffices to prove that, for any indecomposable functors
$\frak f_{Q,\varphi}$ and $\frak f_{R,\psi}$ in $\frak H_\F$  such that $Q$ is not $\F\-$selfcentralizing, the element of $\H_\F$ determined by the composition $\frak f_{Q,\varphi}\circ \frak f_{R,\psi}$ belongs to $\N_\F\,;$
 but, choosing a set of representatives $W$ for~$Q\backslash P/\psi (R)\,,$ with the notation in~£2.13 above we have (cf.~£3.1.1)
 $$\frak f_{Q,\varphi}\circ \frak f_{R,\psi} = \bigsqcup_{w\in W} \frak f_{U_w,\varphi\circ\kappa_w\circ\psi_w}
 \eqno £5.12.1\phantom{.}$$
and therefore if $U_w$ is $\F\-$selfcentralizing then $Q$ is $\F\-$selfcentralizing too  [4,~4.8]. 
\eject

\bigskip
\noindent
{\bf Proposition~£5.13.} {\it  If $\F$ is a Frobenius $P\-$category, for any indecomposable element  $f$ 
in $\H_\F -\N_\F$ there is a  unique {\it maximal\/} indecomposable element $\hat f$ in $\H_\F$ fulfilling 
$\vert f,\hat f\vert\not= 0\,.$ Denote by
$$e_\F : \H_\F\too \H_\F
\eqno £5.13.1\phantom{.}$$
the $\Bbb Z\-$module endomorphism sending $\N_\F$ to $\{0\}$ and mapping  any indecomposable element  $f\in \H_\F -\N_\F$ on  $\displaystyle{\ell (f)\over\ell (\hat f)}\.\hat f\,.$ For any positive elements $g$ and $h$ of~$\H_\F\,,$ the difference 
$ e_\F \big(e_\F (g)e_\F (h)\big) - e_\F (gh)$ is also positive.\/} 

\medskip
\noindent
{\bf Proof:} It follows from Proposition~£3.9 that, if $f$ comes from an indecom-posable functor $\frak f_{Q,\varphi}$
 in~$\frak H_\F$ such that $Q$ is $\F\-$selfcentralizing, there is a unique isomorphism class $\hat f$ of {\it maximal\/} indecomposable functors 
 $\frak f_{\hat Q,\hat\varphi}$ admitting a natural map $\mu\,\colon \frak f_{Q,\varphi}\to \frak f_{\hat Q,\hat\varphi}$
 or, equivalently, fulfilling $\vert f,\hat f\vert\not= 0\,.$

 \smallskip
 Consider a second {\it indecomposable\/} element  $g$ in $\H_\F -\N_\F\,,$ coming from an indecomposable functor 
 $\frak f_{R,\psi}$  in~$\frak H_\F\,,$ and the corresponding isomorphism class $\hat g$ of {\it maximal\/} 
 indecomposable functors  $\frak f_{\hat R,\hat\psi}$ admitting a natural map $\nu\,\colon \frak f_{R,\psi}
 \to \frak f_{\hat R,\hat\psi}\,;$ then, the following commutative diagram
  $$\matrix{&&\frak f_{\hat Q,\hat\varphi}\circ\frak f_{R,\psi}\cr
&\hskip-30pt{\mu *\frak f_{R,\psi}\atop}\hskip-5pt\nearrow &
&\searrow\hskip-4pt{\frak f_{\hat Q,\hat\varphi}*\nu\atop}\hskip-30pt\cr
 \frak f_{Q,\varphi}\circ \frak f_{R,\psi}&&\hbox to 40pt{\rightarrowfill} &
 &\frak f_{\hat Q,\hat\varphi}\circ \frak f_{\hat R,\hat\psi}\cr$$
&\hskip-30pt{\atop\frak f_{Q,\varphi}*\nu}\hskip-5pt\searrow&
&\nearrow\hskip-5pt{\atop\mu *\frak f_{\hat R,\hat\psi}}\hskip-30pt\cr
&&\frak f_{Q,\varphi}\circ\frak f_{\hat R,\hat\psi}\cr}
\eqno £5.13.3\phantom{.}$$
supplies a natural map from $\frak f_{Q,\varphi}\circ \frak f_{R,\psi}$ to $\frak f_{\hat Q,\hat\varphi}
\circ \frak f_{\hat R,\hat\psi}\,;$ actually,  according to~£2.11, we may assume that $Q\i \hat Q\,,$ 
that $\psi (R)\i \hat\psi (\hat R)\,,$ and that
this natural map comes from the canonical map 
$$Q\backslash P /\psi (R)\too \hat Q\backslash P /\hat\psi (\hat R)
\eqno £5.13.4\phantom{.}$$
 and  from a family of natural maps
$$\eta_w : \frak f_{U_w,\varphi\circ\kappa_w\circ\psi_w}\too 
\frak f_{\hat U_{\hat w},\hat\varphi\circ\hat\kappa_{\hat w}\circ\hat\psi_{\hat w}}
\eqno £5.13.5\phantom{.}$$
where, borrowing notation  from~£2.13 above, we have chosen  sets of representatives $\hat W$ for 
 $\hat Q\backslash P /\hat\psi (\hat R)\,,$ $X$ for $\hat Q/Q$ and $Y$ for $\hat\psi (\hat R)/\psi (R)\,,$
and $w$ runs over $W = X^{-1}\.\hat W\.Y$ determining $\hat w\in \hat W\,.$

\smallskip
Moreover, if $U_w$ is $\F\-$selfcentralizing then it follows again from~£2.11 that $\hat U_{\hat w}$ is $\F\-$selfcentralizing too;
 in this case, the isomorphism class $h_{\hat w}$ of {\it maximal\/} 
 indecomposable functors  admitting a natural map from $\frak f_{\hat U_{\hat w},\hat\varphi\circ\hat\kappa_{\hat w}\circ\hat\psi_{\hat w}}$  coincides with the isomorphism class of {\it maximal\/} 
 indecomposable functors  admitting a natural map from $\frak f_{ U_{w},\varphi\circ\kappa_{w}\circ\psi_{w}}\,.$
 \eject

 \smallskip
 Consequently, denoting by $W^{^{\rm sc}}\i W$ and by $\hat W^{^{\rm sc}}\i \hat W$ the corresponding subsets fulfilling these conditions, by the very definition of $e_\F$ we get
 $$e_\F (fg) = \sum_{w\in W^{^{\rm sc}}} {\vert P \colon U_w\vert\over\ell (h_{\hat w})}\.h_{\hat w}\qq 
 e_\F (\hat f\hat g) = \sum_{\hat w\in \hat W^{^{\rm sc}}} {\vert P \colon \hat U_{\hat w}\vert\over\ell (h_{\hat w})}\.h_{\hat w}
 \eqno £5.13.6;$$
but note that, if the subset $X^{-1}\.\hat w\.Y\i W$ has an element in $\hat W^{^{\rm sc}}\,,$ then this subset is actually contained in
$\hat W^{^{\rm sc}}\,;$ hence, we still get
$$\eqalign{e_\F (fg) &= \sum_{\hat w\in W^{^{\rm sc}}\cap \hat W}\Big(\sum_{w\in X^{-1}\.\hat w\.Y} 
{\vert P \colon U_w\vert\over\ell (h_{\hat w})}\Big)\.h_{\hat w}\cr 
 e_\F \big(e_\F ( f)e_\F ( g)\big) &= \vert \hat Q\colon Q\vert\vert \hat R\colon R\vert
 \sum_{\hat w\in \hat W^{^{\rm sc}}} {\vert P \colon \hat U_{\hat w}\vert\over\ell (h_{\hat w})}\.h_{\hat w}\cr}
\eqno £5.13.7;$$

moreover, for any $\hat w\in W^{^{\rm sc}}\cap \hat W\,,$ it is clear that
$$\eqalign{\sum_{w\in X^{-1}\.\hat w\.Y} \vert P \colon U_w\vert 
&= \sum_{w\in X^{-1}\.\hat w\.Y} {\vert P\vert\vert Q\.w\.\psi(R)\vert\over \vert Q\vert \vert R\vert}\cr
& = {\vert P\vert\vert\hat Q\.\hat w\.\hat\psi(\hat R)\vert\over \vert Q\vert \vert R\vert} 
=  \vert \hat Q\colon Q\vert\vert \hat R\colon R\vert \vert P \colon \hat U_{\hat w}\vert\cr}
\eqno £5.13.8;$$
since $\hat W^{^{\rm sc}}$ contains $W^{^{\rm sc}}\cap \hat W\,,$ the difference $ e_\F \big(e_\F (g)e_\F (h)\big) - e_\F (gh)$ is indeed positive. We are done.

\bigskip
\noindent
{\bf £6. The Hecke $\Bbb Z_{(p)}\-$algebras of Frobenius $P\-$categories}
\bigskip

 £6.1. Let $\H$ be a $t\-$stable $\Bbb Z\-$subalgebra of $\B_{P\times P}^{^{\rm pp}}$ containing $\H_P$ 
and fulfilling condition~£5.9.1; for any $f\in \H\,,$ shortly we set $f^{^\circ} = t (f)\,.$ Consider the {\it divisible $P\-$category\/} 
$\F= \F_\H$ determined by $\H$  (cf.~Proposition~£5.9); we see below that whenever $\F$ is a Frobenius $P\-$category, 
the extended $\Bbb Z_{(p)}\-$algebra~$\H_{(p)}$ has a quite particular structure and that, conversely, some of 
these particularities of $\H_{(p)}$  imply that  $\F$ is a Frobenius $P\-$category. For this purpose, let us consider the 
{\it evaluation\/}   $\Bbb Z\-$module homomorphism
$$v_\H: \H\too \B_P
\eqno £6.1.1$$
mapping $f\in \H$ on $f(1)$ where $1$ denotes the isomorphism class of the {\it trivial\/} $P\-$set; then, denote
by $\K_\H$ the kernel of this homomorphism, which is clearly a left-hand ideal.

\medskip
£6.2. Moreover, consider the {\it contravariant\/}  functors from $\F$ to the category $\Bbb Z\-\mod$ of 
$\Bbb Z\-$modules respectively mapping any subgroup $Q$ of $P$  on~$\B_{Q\times P}$ and $\B_{P\times Q}\,,$ 
and any $\F\-$morphism $\varphi\,\colon R\to Q$ on the corresponding ``multiplication'' by ${\rm res}_\varphi$
(cf.~£5.4.1); then, identifying the {\it inverse limits\/} of these\break
\eject
\noindent
 functors with their image in $\B_{P\times P}\,,$
we denote by $(\B_{P\times P})^\F$ their intersection and set
$$\H^\F = \H\cap (\B_{P\times P})^\F\qq (\H^\F)^+ = \H^\F\cap \B_{P\times P}^+
\eqno £6.2.1.$$
Analogously, we respectively denote by $(\B_P)^\F\i \B_P$
and by $\B_P/\F$ the {\it inverse\/} and the {\it direct limits\/}  of the  functors from $\F^{^\circ}$ and from 
$\F$ to $\Bbb Z\-\mod$  sending any subgroup $Q$ of $P$  to $\B_{Q}$ and any $\F\-$morphism 
$\varphi\,\colon R\to Q$ to the actions of ${\rm res}_\varphi$ and of ${\rm ind}_\varphi\,;$ moreover, we denote by
$\bar v_\H$ the composition of $v_\H$ with the canonical map $\B_P\to \B_P/\F$ 
and by  $\bar s_Q\in \B_P/\F$  the image of the isomorphism class $s_Q$ of the $P\-$set $P/Q$ for  any subgroup~$Q$ of $P\,.$

\medskip
£6.3. More explicitly, $f\in \H$ belongs to $\H^\F$ if and only if we have
$${\rm res}_\varphi\.f ={\rm res}_{\iota_Q^P}\.f\qq {\rm res}_\varphi\.f^{^\circ} ={\rm res}_{\iota_Q^P}\.f^{^\circ}
\eqno £6.3.1\phantom{.}$$
 for any subgroup $Q$ of $P$ and any $\varphi\in \F (P,Q)\,;$ in order to relate these equalities with the 
 {\it partial scalar product\/} in $\B_{P\times P}\,,$ let us explicit the set of indecomposable elements of $\H\,;$ 
 for any~$n\in \Bbb N\,,$ let $\s_n$ be  a set of representatives  for the set of $\F\-$isomorphism classes  of subgroups $Q$ of $P$ of index~$p^n\,,$  in such a way that any element of~$\s_n$ is {\it fully normalized\/} in $\F$ [4,~2.6]; moreover, for any $\varphi,\varphi'\in \F (P,Q)\,,$ let us set
$$f_{\varphi,\varphi'} = {\rm ind}_\varphi \. {\rm res}_{\varphi'}\qq \Delta_{\varphi,\varphi'}(Q)
= \{\big(\varphi (u),\varphi' (u)\big)\}_{u\in Q}
\eqno £6.3.2;$$
note that $\bar v_\H (f_{\varphi,\varphi'}) = \bar s_Q\,;$ finally, for any $Q\in \s_n\,,$ let $\D_Q$ be
a set of re-presentatives in the product $\F (P,Q)\times \F (P,Q)$ 
for the set of $\tilde\F(Q)\-$orbits in the product $\tilde\F(P,Q)\times \tilde\F (P,Q)$ through the {\it diagonal\/} action;
we may assume that  $(\iota_Q^P,\iota_Q^P)\in \D_Q\,.$ Then, setting
$$\s = \bigsqcup_{n\in \Bbb N} \s_n\qq \D = \bigsqcup_{Q\in \s} \D_Q
\eqno £6.3.4,$$
the family $\big\{f_{\varphi,\varphi'}\big\}_{(\varphi,\varphi')\in \D}$
is the $\Bbb Z\-$basis of indecomposable elements of $\H\,.$

\bigskip
\noindent
{\bf Lemma~£6.4.} {\it An element $f$ of $\H$  belongs to $\H^\F$ if and only if  we have $\vert f_{\varphi,\varphi'},f\vert 
 = \vert m_{s_Q},f\vert$ for any $Q\in \s$ and any $(\varphi,\varphi')\in \D_Q\,.$
 In particular, for any $f\in \H^\F$ we have $f^{^\circ} = f\,.$\/}

 \medskip
 \noindent
 {\bf Proof:} If $f$ belongs to $\H^\F$ then it follows from isomorphism~£2.8.2 and equalities~£5.7.1 and~£6.3.1 that
  $$\vert f_{\iota_Q^P,\varphi},f\vert = \vert m_{s_Q},f\vert\qq \vert f_{\varphi,\iota_Q^P},f\vert = \vert m_{s_Q},f\vert
 \eqno £6.4.1\phantom{.}$$
  and therefore, denoting by $\varphi^*\,\colon \varphi (Q)\cong Q$ the inverse of the isomorphism determined 
 by $\varphi\,,$ we still get
 $$\eqalign{\vert f_{\varphi,\varphi'},f\vert  &= \vert f_{\iota_{\varphi (Q)}^P,\varphi'\circ\varphi^*},f\vert
 = \vert m_{s_{\varphi (Q)}},f\vert\cr
 & = \vert f_{\iota_Q^P\circ \varphi^*,\iota_{\varphi (Q)}^P},f\vert = \vert f_{\iota_Q^P,\varphi},f\vert 
 = \vert m_{s_Q},f\vert\cr}\eqno £6.4.2.$$
 \eject

 \smallskip
On the other hand, it follows from £5.6 that equalities~£6.3.1 are equivalent to the equalities
 $$\vert g,{\rm res}_\varphi\.f \vert= \vert g,{\rm res}_{\iota_Q^P}\.f\vert \qq 
 \vert g,{\rm res}_\varphi\.f^{^\circ}\vert= \vert g,{\rm res}_{\iota_Q^P}\.f^{^\circ}\vert
\eqno £6.4.3\phantom{.}$$
for any indecomposable element $g$ of $\B_{Q\times P}\,;$ but, it is easily checked that all the members of these equalities
vanish unless $g = {\rm ind}_\theta\.{\rm res}_{\iota_R^P}$ for some $\F\-$morphism $\theta\,\colon R\to Q$
(cf.~£2.6.3), and then we get (cf. isomorphism~£2.8.2 and equality £5.7.1)
$$\eqalign{\vert g,{\rm res}_\varphi\.f \vert = \vert f_{\varphi\circ\theta,\iota_R^P},f \vert&\qq
\vert g,{\rm res}_{\iota_R^P}\.f \vert = \vert f_{\iota_R^P\circ\theta,\iota_R^P},f \vert\cr
\vert g,{\rm res}_\varphi\.f^{^\circ} \vert = \vert f_{\iota_R^P,\varphi\circ\theta},f \vert&\qq
\vert g,{\rm res}_{\iota_R^P}\.f^{^\circ} \vert = \vert f_{\iota_R^P,\iota_R^P\circ\theta},f \vert\cr}
\eqno £6.4.4\,;$$
consequently,  equalities~£6.3.1 follow from equalities~£6.4.2.

\smallskip
In particular, equalities~£6.4.2 force (cf.~isomorphism~£2.8.2)
$$\vert f_{\varphi,\varphi'},f\vert = \vert f_{\varphi',\varphi},f\vert = \vert f_{\varphi,\varphi'},f^{^\circ}\vert
\eqno £6.4.5\phantom{.}$$
and therefore, according to~£5.6, we get $f = f^{^\circ}\,.$
We are done.

\bigskip
\noindent
{\bf Proposition~£6.5.} {\it With the notation above, $\,(\B_P)^\F$ is a unitary $\Bbb Z\-$subalgebra of~$\B_P$
of $\,\Bbb Z\-$rank equal to ${\rm rank}_\Bbb Z (\B_P/\F)\,,$ $\H^\F$ is a $\Bbb Z\-$subalgebra of $\H$ fulfilling
$\,\H^\F\.\B_P\i   (\B_P)^\F\,,$ and we have
$$\K_\H\.\H^\F = \{0\} = \K_\H\cap \H^\F \qq \H^\F = (\H^\F)^+ - (\H^\F)^+ 
\eqno £6.5.1.$$
Moreover, any $f\in \H^\F$ commutes with any $g\in \K_\H + \H^\F$ such that $g^{_\circ} = g\,.$\/}

\medskip
\noindent
{\bf Proof:} It is clear that $s\in \B_P$ belongs to $(\B_P)^\F$ if and only if we have
$${\rm res}_\varphi (s) ={\rm res}_{\iota_Q^P}(s)
\eqno £6.5.2\phantom{.}$$
 for any subgroup $Q$ of $P$ and any $\varphi\in \F (P,Q)\,;$ since ${\rm res}_\varphi$ and ${\rm res}_{\iota_Q^P}$
are both unitary $\Bbb Z\-$algebra homomorphisms,  $(\B_P)^\F$ is indeed a unitary $\Bbb Z\-$subalgebra of~$\B_P\,.$
Moreover, equality~£6.3.1 and the  associativity of the ``multiplication'' show that  $\H^\F$ is a $\Bbb Z\-$subalgebra of $\H$
and that, for any $s\in \B_P\,,$  we have
 $${\rm res}_\varphi\big(f(s)\big) ={\rm res}_{\iota_Q^P}\big(f(s)\big)
\eqno £6.5.3\phantom{.}$$
which proves that $\,\H^\F\.\B_P\i   (\B_P)^\F\,.$

\smallskip
On the other hand, it is clear that the family $\{f_{\iota_Q^P,\varphi} - m_{s_Q}\}_{(\iota_Q^P,\varphi)\in \D}$
is a $\Bbb Z\-$basis of $\K_\H$ and therefore the  associativity of the ``multiplication'' and equalities~£6.3.1 
force $\K_\H\.\H^\F = \{0\}\,.$ In particular,  for any $f\in \K_\H\cap \H^\F $ we get $f^2 = 0$ and therefore, setting 
$$f = \sum_{(\varphi,\varphi')\in \D} z_{\varphi,\varphi'}\.f_{\varphi,\varphi'}
\eqno £6.5.4\phantom{.}$$
\eject
\noindent
 for suitable $z_{\varphi,\varphi'}\in \Bbb Z$ and choosing sets of representatives $W^{\varphi,\varphi'}_{\psi,\psi'}$ 
 for the sets of double classes $\varphi'(Q)\backslash P /\psi(R)$ for any $(\varphi,\varphi')\in \D_Q$ and any 
 $(\psi,\psi')\in \D_R$ where $Q,R\in \s\,,$ we still have (cf.~£3.1.1)
$$f^2 = \sum_{(\varphi,\varphi'),(\psi,\psi')\in \D}\,\sum_{w\in W^{\varphi,\varphi'}_{\psi,\psi'}} 
z_{\varphi,\varphi'}z_{\psi,\psi'}\.f_{\varphi\circ\varphi'^*_w\circ \kappa_w, \psi'\circ \psi^*_w}
\eqno £6.5.5\phantom{.}$$
where $\varphi'^*\,\colon\varphi' (Q)\cong Q$ and $\psi^*\,\colon \psi (R)\cong R$ denote the inverse of the
respective isomorphisms determined by $\varphi'$ and $\psi\,,$ and we denote by
$$\varphi'^*_w : \varphi' (Q)\cap {}^w \psi (R)\too Q\qq \psi^*_w : \varphi' (Q)^w\cap \psi (R)\too R
\eqno £6.5.6\phantom{.}$$
the corresponding restrictions.

\smallskip
Thus, if we assume that $f\not= 0$ and choose $(\hat\eta,\hat\eta') \in \D_T$ fulfilling $z_{\hat\eta,\hat\eta'}\not= 0$ 
for a $T\in \s$ such that $\vert T\vert$ is maximal then, setting 
$$f_T = \sum_{(\eta,\eta')\in \D_T} z_{\eta,\eta'}\.f_{\eta,\eta'}
\eqno £6.5.7,$$ 
it is clear that $f^2 = 0$ forces $(f_T)^2 = 0\,;$ but, since $f$ belongs to $\H^\F\,,$ it follows from Lemma~£6.4 
that for any $(\eta,\eta')\in \D_T$ we have
$$z_{\eta,\eta'}\vert f_{\eta,\eta'},f_{\eta,\eta'}\vert = \vert f_{\eta,\eta'},f_T\vert  =\vert f_{\eta,\eta'},f\vert  = \vert m_{s_T},f\vert  = \vert m_{s_T},f_T\vert
\eqno £6.5.8,$$
so that either $f_T$ or $-f_T$ is {\it positive\/} and therefore we get $(f_T)^2 \not= 0\,,$ a contradiction;
consequently, we finally get $K_\H\cap \H^\F = \{0\}\,.$

\smallskip
 Moreover, for any $g\in \K_\H + \H^\F$ such that $g^{_\circ} = g\,,$ setting  $g = g' + g''$ with 
 $g'\in \K_\H$ and $g''\in \H^\F\,,$
we still have $g'^{_\circ} = g'\,;$ then, the equality $\K_\H\.\H^\F = \{0\}$ and Lemma~£6.4 above imply that
$$(fg')^{\!^\circ} = g'^{_\circ}f^{^\circ} = g'f = 0\qq   fg'' = f^{^\circ}g''^{_\circ} = (g''f)^{\!^\circ} = g''f
 \eqno £6.5.9,$$
 so that $fg = gf\,.$

\smallskip
Now, we claim that, for any $f\in \H^\F\,,$ there are {\it positive\/} (cf.~£5.11) elements $f'$ and $f''$ of $\H^\F$ such that we have $f = f' - f''\,;$ that is to say, we claim that, setting 
$$f = \sum_{(\varphi,\varphi')\in \D} z_{\varphi,\varphi'}\.f_{\varphi,\varphi'}
\eqno £6.5.10,$$
there are $z'_{\varphi,\varphi'},z''_{\varphi,\varphi'} \in \Bbb N\,,$
where $(\varphi,\varphi')$ runs over $\D\,,$ such that we have $z_{\varphi,\varphi'} 
= z'_{\varphi,\varphi'} - z''_{\varphi,\varphi'}$
and that, setting
$$f' = \sum_{(\varphi,\varphi')\in \D} z'_{\varphi,\varphi'}\.f_{\varphi,\varphi'}\qq
f'' = \sum_{(\varphi,\varphi')\in \D} z''_{\varphi,\varphi'}\.f_{\varphi,\varphi'}
\eqno £6.5.11,$$ 
for any $Q\in \s$ and any $(\varphi,\varphi')\in \D_Q$
we still have (cf.~Lemma~£6.4)
$$\vert f_{\varphi,\varphi'},f' \vert = \vert m_{s_Q},f' \vert \qq 
\vert f_{\varphi,\varphi'}f''\vert = \vert m_{s_Q},f''\vert
\eqno £6.5.12.$$
\eject

\smallskip
We argue by induction on $\vert P\colon Q\vert $ and assume that, for $m\in \Bbb N\,,$ we have found 
$z'_{\varphi,\varphi'},z''_{\varphi,\varphi'} \in \Bbb N\,,$ where $(\varphi,\varphi')$ runs over $\D_Q$ for any
$Q\in \s$ fulfilling $\vert P\colon Q\vert < p^m\,,$ such that we have
 $z_{\varphi,\varphi'} = z'_{\varphi,\varphi'} - z''_{\varphi,\varphi'}$ and that,
 setting
$$\eqalign{f'_m = \sum_Q\,\sum_{(\varphi,\varphi')\in \D_Q} z'_{\varphi,\varphi'}\.f_{\varphi,\varphi'}\cr
f''_m = \sum_Q\,\sum_{(\varphi,\varphi')\in \D_Q} z''_{\varphi,\varphi'}\.f_{\varphi,\varphi'}\cr}
\eqno £6.5.13\phantom{.}$$ 
where $Q$ runs over the elements of $\s$ fulfilling $\vert P\colon Q\vert < p^m\,,$
we still have
$$\vert f_{\varphi,\varphi'},f'_m \vert = \vert m_{s_Q},f'_m \vert \qq 
\vert f_{\varphi,\varphi'}f''_m\vert = \vert m_{s_Q},f'_m\vert
\eqno £6.5.14\phantom{.}$$
for any $Q\in \s$ fulfilling $\vert P\colon Q\vert < p^m\,,$ and any $(\varphi,\varphi')\in \D_Q\,.$

\smallskip
Thus, since $f$ belongs to~$\H^\F\,,$  it follows from Lemma~£6.4 that  for any $R\in \s_m$ and any $(\psi,\psi')\in \D_R$ we have
$$\eqalign{\vert f_{\iota_R^P,\iota_R^P},f'_m - f''_m \vert + z_{\iota_R^P,\iota_R^P} &\vert f_{\iota_R^P,\iota_R^P},f_{\iota_R^P,\iota_R^P} \vert\cr  
&= \vert f_{\psi,\psi'},f'_m - f''_m  \vert + z_{\psi,\psi'} \vert f_{\psi,\psi'},f_{\psi,\psi'} \vert\cr}
\eqno £6.5.15;$$
moreover, we know that (cf.~£3.6.1 and~£5.5.1)
$$\vert f_{\psi,\psi'},f_{\psi,\psi'} \vert = \big\vert \bar N_{P\times P}\big(\Delta_{\psi,\psi'} (R)\big)\big\vert
\eqno £6.5.16\phantom{.}$$
and therefore equalities~£6.5.15 considered in $\Bbb Q$ are obviously equivalent to
$$z_{\psi,\psi'} =  {\vert f_{\iota_R^P,\iota_R^P},f_m \vert - \vert f_{\psi,\psi'},f_m\vert 
\over \big\vert \bar N_{P\times P} \big(\Delta_{\psi,\psi'} (R)\big)\big\vert} 
+ z_{\iota_R^P,\iota_R^P} {\big\vert\bar N_{P\times P}
\big(\Delta (R)\big)\big\vert\over \big\vert \bar N_{P\times P}\big(\Delta_{\psi,\psi'} (R)\big)\big\vert}
\eqno £6.5.17;$$
but, up to a suitable modification of our choice of $\s_m$ in~£6.3 above, we may assume that the order of 
$\bar N_{P\times P}\big(\Delta_{\psi,\psi'} (R)\big)$ divides the order of  $\bar N_{P\times P}\big(\Delta (R)\big)\,;$ 
then,  the order of $\bar N_{P\times P}\big(\Delta_{\psi,\psi'} (R)\big)$ also divides $\vert f_{\iota_R^P,\iota_R^P},f_m \vert - \vert f_{\psi,\psi'},f_m\vert \,.$

\smallskip
At this point, we clearly can choose  $z'_{\iota_R^P,\iota_R^P},z''_{\iota_R^P,\iota_R^P} \in \Bbb N$
in such a way that we have $z_{\iota_R^P,\iota_R^P} = z'_{\iota_R^P,\iota_R^P}- z''_{\iota_R^P,\iota_R^P}$
and that the integers 
$$\eqalign{z'_{\psi,\psi'} =  {\vert f_{\iota_R^P,\iota_R^P},f_m \vert - \vert f_{\psi,\psi'},f_m\vert 
\over \big\vert \bar N_{P\times P} \big(\Delta_{\psi,\psi'} (R)\big)\big\vert} 
+ z'_{\iota_R^P,\iota_R^P} {\big\vert\bar N_{P\times P}
\big(\Delta (R)\big)\big\vert\over \big\vert \bar N_{P\times P}\big(\Delta_{\psi,\psi'} (R)\big)\big\vert}\cr
z''_{\psi,\psi'} =  {\vert f_{\iota_R^P,\iota_R^P},f_m \vert - \vert f_{\psi,\psi'},f_m\vert 
\over \big\vert \bar N_{P\times P} \big(\Delta_{\psi,\psi'} (R)\big)\big\vert} 
+ z''_{\iota_R^P,\iota_R^P} {\big\vert\bar N_{P\times P}
\big(\Delta (R)\big)\big\vert\over \big\vert \bar N_{P\times P}\big(\Delta_{\psi,\psi'} (R)\big)\big\vert}\cr}
\eqno £6.5.18\phantom{.}$$
are all {\it positive\/};  then equalities~£6.5.16 imply that $z_{\psi,\psi'} = z'_{\psi,\psi'}- z''_{\psi,\psi'}\,.$
Hence, we have found
$$\eqalign{g'_m &= \sum_{R\in \s_m}\sum_{(\psi,\psi')\in \D_R} z'_{\psi,\psi'}\.f_{\psi,\psi'}\cr
g''_m &= \sum_{R\in \s_m}\sum_{(\psi,\psi')\in \D_R} z''_{\psi,\psi'}\.f_{\psi,\psi'}\cr}
\eqno £6.6.19\phantom{.}$$
such that, setting $f'_{m+1} = f'_m+ g'_m$ and  $f''_{m+1} = f''_m+ g''_m\,,$ we still have
$$\vert f_{\varphi,\varphi'},f'_{m +1} \vert = \vert m_{s_Q},f'_{m +1} \vert \qq 
\vert f_{\varphi,\varphi'}f''_{m +1}\vert = \vert m_{s_Q},f'_{m +1}\vert
\eqno £6.5.20\phantom{.}$$
for any $Q\in \s$ fulfilling $\vert P\colon Q\vert \le p^m$ and any 
$(\varphi,\varphi')\in \D_Q\,.$ This proves our claim above.

\smallskip
Finally, it remains to prove that the $\Bbb Z\-$ranks of $(\B_P)^\F$ and of $\B_P/\F$ coincide.
For any subgroup $Q$ of $P\,,$ it is quite clear that the bilinear map 
$$\Bbb Q\otimes_\Bbb Z \B_Q\times \Bbb Q\otimes_\Bbb Z \B_Q\too \Bbb Q
\eqno £6.5.21\phantom{.}$$
sending $(1\otimes t,1\otimes t')$ to $\vert {\rm Hom}_Q (t,t')\vert$ for any $t,t'\in \B_Q$ determines 
a $\Bbb Q\-$linear isomorphism between $\Bbb Q\otimes_\Bbb Z \B_Q$ and its {\it dual\/} $(\Bbb Q\otimes_\Bbb Z \B_Q)^*\,;$
but, for any $\F\-$mor-phism $\varphi\,\colon R\to Q\,,$ it is easily checked that we have a commutative diagram
$$\matrix{(\Bbb Q\otimes_\Bbb Z \B_R)^*&\buildrel ({\rm res}_\varphi)^*\over\too &(\Bbb Q\otimes_\Bbb Z \B_Q)^*\cr
\wr\Vert&&\wr\Vert\cr
\Bbb Q\otimes_\Bbb Z \B_R&\buildrel {\rm ind}_\varphi \over\too &\Bbb Q\otimes_\Bbb Z \B_Q\cr}
\eqno £6.5.22;$$
hence, since the {\it dual\/} of an {\it inverse limit\/} of some functor to $\Bbb Q\-\mod$ is the {\it direct limit\/} 
of the {\it dual\/} of this functor, we get a $\Bbb Q\-$linear isomorphism and we are done
$$\Bbb Q\otimes_\Bbb Z \B_P/\F\cong \big(\Bbb Q\otimes_\Bbb Z (\B_P)^\F\big)^*
\eqno £6.5.23.$$
We are done.

\bigskip
\noindent
{\bf Theorem~£6.6.} {\it With the notation above, if $\F$ is a Frobenius $P\-$category then $(\K_\H + \H^\F)_{(p)}$ is a
unitary $\Bbb Z_{(p)}\-$subalgebra of $\H_{(p)}$ containing $m^P\big((\B_P)^\F\big)$ and we have
$$(\K_\H + \H^\F)_{(p)} \cong (\K_\H)_{(p)}\rtimes \H^\F_{(p)}
\eqno £6.6.1.$$
Moreover, $\H^\F$ centralizes $m^P\big((\B_P)^\F\big)\,.$ In particular, $v_\H$ determines a $\Bbb Z\-$al-gebra isomorphism 
$$\H^\F_{(p)}\cong (\B_P)^\F_{(p)}
\eqno £6.6.2.$$\/}

\par
\noindent
{\bf Remark~£6.7.} For any $f\in \H^\F\,,$ since $f(s)$ belongs to $(\B_P)^\F$ for any $s\in \B_P\,,$  if $s'\in \B_P$ then we have
$$f\big(f(s)\,s'\big) = (f\,m_{f(s)})(s') = (m_{f(s)} f)(s') = f(s)f(s')
\eqno £6.7.1$$
which amounts to saying that $f$ fulfills the {\it Frobenius condition\/} introduced in~[6]. Moreover, the isomorphism
$(\H^\F)_{(p)}\cong (\B_P)^\F_{(p)}$ shows that $(\H^\F)_{(p)}$ contains a unique nonzero idempotent
which coincides with the idempotent exhibit in~[5].

\medskip
\noindent
{\bf Proof:} First of all we claim that, for any $\bar s\in \B_P/\F$ there is $f\in  \H^\F_{(p)}$ such 
that $\bar v_\H(f) = \bar s\,;$ with the notation in~£6.3 above, for suitable $z_R\in \Bbb Z\,,$  we have 
$$\bar s = \sum_{n\in \Bbb N}\sum_{R\in \s_n} z_R\.\bar s_R
\eqno £6.6.3\phantom{.}$$
and we set $\bar s_n = \sum_{m\le n}\sum_{R\in \s_m} z_R\.\bar s_R\,.$

\smallskip
Let $\{f_n\}_{n\in \Bbb N}$ be a family of elements of $\H_{(p)}$  inductively defined as follows. We set
$$f_0 = {z_P\over \vert\tilde\F (P)\vert}\. \sum_{\tilde\sigma \in \tilde\F (P)} f_{\sigma,{\rm id}_P}
\eqno £6.6.4,$$
 so that we have $\bar v_\H (f_0) = \bar s_0$ and 
$$\vert f_{{\rm id}_P,{\rm id}_P},f_0 \vert = \vert f_{\sigma,\sigma'},f_0 \vert 
\eqno £6.6.5\phantom{.}$$
for any $\sigma,\sigma'\in \F(P)$ (cf.~£6.3). If $n\ge 1$ and $m < n\,,$ we may assume that $f_m$ belongs to 
$\H_{(p)}\,,$ that $\bar v_\H (f_m) = \bar x_m$ and that we have
$$\vert f_{\iota_Q^P,\iota_Q^P},f_m \vert = \vert f_{\varphi,\varphi'},f_m \vert 
\eqno £6.6.6\phantom{.}$$
 for any $Q\in  \s_m$ and any $\varphi,\varphi'\in\D_Q\,.$ Then we claim that for any $R\in \s_n$ and any $(\psi,\psi')\in \D_R\,,$ we can  choose $k_{\psi,\psi'}\in \Bbb Z_{(p)}$ fulfilling 
the following equalities
$$\eqalign{\vert f_{\iota_R^P,\iota_R^P},f_{n-1} \vert + k_{\iota_R^P,\iota_R^P} &\vert f_{\iota_R^P,\iota_R^P},f_{\iota_R^P,\iota_R^P} \vert\cr  
&= \vert f_{\psi,\psi'},f_{n-1} \vert + k_{\psi,\psi'} \vert f_{\psi,\psi'},f_{\psi,\psi'} \vert\cr
\sum_{(\psi,\psi')\in \D_R} k_{\psi,\psi'} = z_R\cr}
\eqno £6.6.7.$$

\smallskip
Since we have (cf.~£3.6.1 and~£5.5.1)
$$\vert f_{\psi,\psi'},f_{\psi,\psi'} \vert = \big\vert \bar N_{P\times P}\big(\Delta_{\psi,\psi'} (R)\big)\big\vert
\eqno £6.6.8,$$
equalities~£6.6.7 considered in $\Bbb Q$ are obviously equivalent to
$$\eqalign{k_{\psi,\psi'} = & {\vert f_{\iota_R^P,\iota_R^P},f_{n-1} \vert - \vert f_{\psi,\psi'},f_{n-1} 
\vert\over \big\vert \bar N_{P\times P} \big(\Delta_{\psi,\psi'} (R)\big)\big\vert} 
+ k_{\iota_R^P,\iota_R^P} {\big\vert\bar N_{P\times P}
\big(\Delta (R)\big)\big\vert\over \big\vert \bar N_{P\times P}\big(\Delta_{\psi,\psi'} (R)\big)\big\vert}\cr
k_{\iota_R^P,\iota_R^P} & \sum_{(\psi,\psi')\in \D_R}  {\big\vert\bar N_{P\times P}\big(\Delta (R)\big)\big\vert 
\over \big\vert \bar N_{P\times P}\big(\Delta_{\psi,\psi'} (R)\big)\big\vert}\cr
&= z_R - \sum_{(\psi,\psi')\in \D_R}{\vert f_{\iota_R^P,\iota_R^P},f_{n-1} \vert - \vert f_{\psi,\psi'},f_{n-1} \vert
\over \big\vert \bar N_{P\times P} \big(\Delta_{\psi,\psi'} (R)\big)\big\vert}\cr}
\eqno £6.6.9.$$

\smallskip
But,  it follows from Lemmas~£6.8 and~£6.9 below that 
$\big\vert \bar N_{P\times P}\big(\Delta_{\psi,\psi'} (R)\big)\big\vert$
 divides  $\big\vert \bar N_{P\times P}\big(\Delta (R)\big)\big\vert$ for any pair
 $(\psi,\psi')\in \D_R$ and that the sum 
 $$q_R = \sum_{(\psi,\psi')\in \D_R}  {\big\vert\bar N_{P\times P}\big(\Delta (R)\big)\big\vert\over \big\vert \bar N_{P\times P}\big(\Delta_{\psi,\psi'} (R)\big)\big\vert}
 \eqno £6.6.10\phantom{.}$$
 is an integer prime to $p\,.$ 
 Moreover, since $\F$ is a Frobenius $P\-$category, it follows from~[4,~condition~2.8.2] and again from Lemma~£6.8
 below that there is a group homomorphism 
$$\zeta\times \zeta' : N_{P\times P}\big(\Delta_{\psi,\psi'} (R)\big)\too P\times P
\eqno £6.6.11\phantom{.}$$
which maps $\Delta_{\psi,\psi'} (R)$ onto~$\Delta (R)$ and therefore it fulfills
$$(\zeta\times\zeta')\Big(N_{P\times P}\big(\Delta_{\psi,\psi'} (R)\big)\Big) \i N_{P\times P}\big(\Delta (R)\big)
\eqno £6.6.12;$$

\smallskip
At this point, since we have $k\.f_{n-1}\in \H$ for a suitable $k\in \Bbb Z- p\Bbb Z\,,$ we still have 
$k\.f_{n-1} = g - g'$ where $g$ and $g'$ are the isomorphism classes of functors $\frak f_\Omega$ and 
$\frak f_{\Omega'}$ in $\frak H_\F$ for some $P\times P\-$sets $\Omega$ and $\Omega'\,.$ Consider the 
$N_{P\times P}\big(\Delta (R)\big)\-$sets $\Omega^{\Delta(R)}$ and $\Omega'^{\Delta(R)}\,,$ which become
$N_{P\times P}\big(\Delta_{\psi,\psi'} (R)\big)\-$sets {\it via\/} the group homomorphism $\zeta\times \zeta'\,,$
and the $N_{P\times P}\big(\Delta_{\psi,\psi'} (R)\big)\-$sets $\Omega^{\Delta_{\psi,\psi'}(R)}$ and 
$\Omega'^{\Delta_{\psi,\psi'}(R)}\,;$ then,  for any  nontrivial subgroup $\bar D$ of 
$\bar N_{P\times P}\big(\Delta_{\psi,\psi'} (R)\big)\,,$ denoting by $D$ the converse image of $\bar D$ in 
$N_{P\times P}\big(\Delta_{\psi,\psi'} (R)\big)\,,$ the fixed points of $\bar D$ in 
$$\Omega^{\Delta(R)}\quad ,\quad\Omega'^{\Delta(R)}\quad ,\quad\Omega^{\Delta_{\psi,\psi'}(R)}
\qq \Omega'^{\Delta_{\psi,\psi'}(R)}
\eqno £6.6.13\phantom{.}$$
respectively coincide with 
$$\Omega^{(\zeta\times \zeta')(D)}\quad ,\quad\Omega'^{(\zeta\times \zeta')(D)}\quad ,\quad\Omega^D
\qq \Omega'^D
\eqno £6.6.14,$$
 and therefore it follows from~equalities~£6.6.6 and from the fact that 
$g$ and $g'$ belong to $\H$ that we have (cf.~£2.6.3 and~£5.5)
$$\eqalign{\vert\Omega^{(\zeta\times \zeta')(D)}\vert - \vert\Omega'^{(\zeta\times \zeta')(D)}\vert 
&= \vert f_{(\zeta\times \zeta')(D)},g\vert - \vert f_{(\zeta\times \zeta')(D)},g'\vert\cr
& =\vert f_{(\zeta\times \zeta')(D)},k\.f_{n-1}\vert\cr
& = \vert f_D,k\.f_{n-1} \vert = \vert\Omega^D\vert - \vert\Omega'^D\vert \cr}
\eqno £6.6.15\phantom{.}$$
where $f_{(\zeta\times \zeta')(D)}$ and $ f_D$ respectively denote the isomorphism classes of the functors
$\frak f_{(\zeta\times \zeta')(D)}$ and $\frak  f_D\,.$

\smallskip
This is clearly equivalent to saying that the multiplicities of any non-regular indecomposable 
$\bar N_{P\times P}\big(\Delta_{\psi,\psi'} (R)\big)\-$set in  
$$\Omega^{\Delta(R)}\sqcup\Omega'^{\Delta_{\psi,\psi'}(R)}\qq \Omega^{\Delta_{\psi,\psi'}(R)}\sqcup \Omega'^{\Delta(R)}
\eqno £6.6.16\phantom{.}$$ 
coincide with each other; hence, one of these $\bar N_{P\times P}\big(\Delta_{\psi,\psi'} (R)\big)\-$sets is isomorphic to the disjoint union of the other with a multiple of the regular
$\bar N_{P\times P}\big(\Delta_{\psi,\psi'} (R)\big)\-$set; in particular, the order of $ \bar N_{P\times P}
\big(\Delta_{\psi,\psi'} (R)\big)$ divides the following difference (cf.~£2.6.3 and~£5.5)
$$\eqalign{\vert\Omega^{\Delta(R)}\vert + &\vert\Omega'^{\Delta_{\psi,\psi'}(R)}\vert - \vert\Omega'^{\Delta(R)}\vert - 
\vert\Omega^{\Delta_{\psi,\psi'}(R)}\vert \cr
& = \vert f_{\iota_R^P,\iota_R^P},g \vert  + \vert f_{\psi,\psi'},g' \vert -\vert f_{\iota_R^P,\iota_R^P},g' \vert 
- \vert f_{\psi,\psi'},g \vert\cr
& = \vert f_{\iota_R^P,\iota_R^P},k\.f_{n-1} \vert - \vert f_{\psi,\psi'},k\.f_{n-1} \vert\cr}
\eqno £6.6.17,$$
proving our claim.

\smallskip
 Now, it suffices to choose 
$$f_n = f_{n-1} + \sum_{R\in\s_n}
\,\sum_{(\psi,\psi')\in \D_R}k_{\psi,\psi'}\. f_{\psi,\psi'}
\eqno £6.6.18;$$
indeed, it is clear that $f_n$ belongs to $\H_{(p)}\,;$ moreover, since 
$\bar v_\H (f_{\psi,\psi'}) = \bar s_R\,,$ the last equality in~£6.6.7 above shows that $\bar v_\H (f_n) = \bar s_n\,;$
finally, it is easily checked that, for any $R,\hat R\in \s_n\,,$ any $(\psi,\psi')\in \D_R$ and any 
$(\hat\psi,\hat\psi')\in \D_{\hat R}$
the inequality $\vert f_{\psi,\psi'},f_{\hat\psi,\hat\psi'}\vert \not= 0$ forces $(R,\psi,\psi') = (\hat R,\hat\psi,\hat\psi')\,,$
and therefore from equalities~£6.6.7 we obtain
$$\eqalign{\vert f_{\iota_R^P,\iota_R^P},f_n\vert &= \vert f_{\iota_R^P,\iota_R^P},f_{n-1} \vert + k_{\iota_R^P,\iota_R^P} \vert f_{\iota_R^P,\iota_R^P},f_{\iota_R^P,\iota_R^P} \vert\cr
& = \vert f_{\psi,\psi'},f_{n-1} \vert + k_{\psi,\psi'} \vert f_{\psi,\psi'},f_{\psi,\psi'} \vert = \vert f_{\psi,\psi'},f_n\vert\cr}
\eqno £6.6.19.$$
In conclusion, setting $f = f_n$ for $n$ big enough, we obtain $\bar v_\H (f) = \bar s$ and 
$$\vert f_{\varphi,\varphi'},f\vert = \vert f_{\iota_Q^P,\iota_Q^P},f\vert
\eqno £6.6.20\phantom{.}$$
for any $Q\in \s$ and any $(\varphi,\varphi')\in \D_Q$ which, according to Lemma~£6.4,  proves that $f$ belongs  to~$(\H^\F)_{(p)}\,.$
\eject

\smallskip
On the other hand, since $\K_\H$ is a left-hand ideal in $\H\,,$ so that $\K_\H$ contains $(\K_\H + \H^\F)\.\K_\H\,,$ 
and since it follows from Proposition~£6.5 that
$$\H^\F\.\H^\F\i \H^\F\qq \K^\F\.\H^\F = \{0\}
\eqno £6.6.21,$$
the sum $\K_\H + \H^\F$ is a $\Bbb Z\-$subalgebra of $\H\,;$ moreover, according again to  Proposition~£6.5,
we have
 $$v_\H (\H^\F)\i (\B_P)^\F\qq {\rm rank}_\Bbb Z \big( (\B_P)^\F\big) = {\rm rank}_\Bbb Z  (\B_P/\F)
 \eqno £6.6.22.$$
 But, the argument above proves that
 $$\bar v_\H(\H^\F_{(p)}) = (\B_P/\F)_{(p)}
 \eqno £6.6.23.$$
 Consequently, since $\K_H\cap \H^\F = \{0\}\,,$ the canonical map $\B_P\to \B_P/\F$ and the {\it evaluation\/} map
  $v_\H$ induce  $\Bbb Z_{(p)}\-$module isomorphisms 
 $$(\B_P)^\F_{(p)}\cong (\B_P/\F)_{(p)}\qq \H^\F_{(p)}\cong (\B_P)^\F_{(p)}
 \eqno £6.6.24\phantom{.}$$ 
 and $m^P\big((\B_P)^\F\big)$ is contained in $ (\K_\H + \H^\F)_{(p)}\,;$
 in particular,  $(\K_\H + \H^\F)_{(p)}$ is a unitary $\Bbb Z_{(p)}\-$subalgebra of $\H_{(p)}$
 and we have
$$(\K_\H + \H^\F)_{(p)}\cong (\K_\H)_{(p)}\rtimes \H^\F_{(p)}
 \eqno £6.6.25\phantom{.}$$
 where the symbol $\rtimes$ means that it is a {\it direct product\/} of $\Bbb Z_{(p)}\-$modules
 with a product defined by the product in each term and the (left-hand) $\H^\F_{(p)}\-$module structure of $(\K_\H)_{(p)}\,.$

 \smallskip
 Finally, for any $s\in (\B_P)^\F\,,$ since $m_s$ belongs to $ (\K_\H + \H^\F)_{(p)}$ and we have 
 $(m_s)^{^\circ} = m_s\,,$ it follows from proposition~£6.5 that $m_s$ centralizes $\H^\F\,.$
 In particular, for any $f,g\in \H^\F\,,$ since $g(1)$ belongs to~$(\B^P)^\F\,,$ we obtain
 $$(fg)(1) = f\big(g(1)\big) = (fm_{g(1)})(1) = (m_{g(1)}f)(1) = g(1)f(1)
 \eqno £6.6.26\phantom{.}$$
 which, together with  isomorphisms~£6.6.24, proves the last
 $\Bbb Z_{(p)}\-$algebra isomorphism. We are done.

\bigskip
\noindent
{\bf Lemma~£6.8.} {\it Assume that $\F$ is a Frobenius $P\-$category and let $Q$ be a subgroup of $P$ 
fully normalized in $\F\,.$ Then $\Delta (Q)$ is fully normalized in the Frobenius $P\times P\-$category $\F\times \F\,.$ \/}
\medskip
\noindent
{\bf Proof:} Any $\F\times \F\-$morphism from $\Delta (Q)$ to $P\times P$ is determined by two $\F\-$morphisms $\varphi,\varphi' \in \F(P,Q)\,;$ since we assume that $Q$ is {\it fully normalized\/} in~$\F$ [4~2.6], considering the inverse of the
isomorphism $\varphi\,\colon Q\cong \varphi (Q)$ determined by $\varphi\,,$ it follows from
condition [4,~£2.12.2] that there are an $\F\-$morphism $\zeta\, \,\colon N_P\big(\varphi (Q)\big)\to P$ and 
an element $\chi\in \F (Q)$~such that $\zeta\big(\psi (u)\big) = \chi (u)$ for~any $u\in Q\,;$ moreover, according
to [4,~Proposition~£2.13], $Q$ is also fully centralized in~$\F$ and therefore, considering
the group homomorphism
$$\eta : \varphi' (Q)\too P
\eqno £6.8.1\phantom{.}$$
 mapping~$\varphi' (u)$ on $\chi (u)$~for~any $u\in Q\,,$ and the image
$$R = \pi'\Big(N_{P\times P}\big(\Delta_{\varphi,\varphi'}(Q)\big)\Big)\i N_P\big(\varphi' (Q)\big)
\eqno £6.8.2\phantom{.}$$ 
where $\pi'\,\colon P\times P\to P$ is the second projection, there is $\zeta'\in \F (P,R)$
extending $\eta$ and, in particular, fulfilling $\zeta'\big(\psi' (u)\big) = \chi (u)$ for any $u\in Q$ 
[4,~statement~2.10.1]; thus, $\zeta\times\zeta'$ determines a group
homomorphism 
$$N_{P\times P}\big(\Delta_{\varphi,\varphi'} (Q)\big)\too P\times P
\eqno £6.8.3\phantom{.}$$
which maps $\Delta_{\varphi,\varphi'} (Q)$ onto~$\Delta (Q)$ and therefore it fulfills
$$(\zeta\times\zeta')\Big(N_{P\times P}
\big(\Delta_{\varphi,\varphi'} (Q)\big)\Big)\i N_{P\times P}\big(\Delta (Q)\big)
\eqno £6.8.4.$$
We are done.

\bigskip
\noindent
{\bf Lemma~£6.9.}  {\it Assume that $\F$ is a Frobenius $P\-$category, let $Q$ be a subgroup of~$P$ and~set
$$\eqalign{\tilde\F (P,Q)_{\rm fc} &= \{\tilde\varphi\in \tilde\F (P,Q)\mid \hbox{$\varphi (Q)$ is fully centralized in $\F$}\}\cr
\tilde\F (P,Q)_{\rm fn} &= \{\tilde\varphi\in \tilde\F (P,Q)\mid \hbox{$\varphi (Q)$ is fully normalized in $\F$}\}\cr}
\eqno £6.9.1.$$
Then, $\vert\tilde\F (P,Q)_{\rm fc}\vert$ and $\vert\tilde\F (P,Q)_{\rm fn}\vert$ are both prime to $p\,.$ In particular, if $Q$ is fully normalized in $\F\,,$ the sum 
$$\sum_{\tilde\varphi\in \tilde\F (P,Q)}{\vert \bar N_P(Q)\vert\over \vert \bar N_P\big(\varphi(Q)\big)\vert}
\eqno £6.9.2\phantom{.}$$
is an integer prime to $p\,.$\/}

\medskip
\noindent
{\bf Proof:} We actually may assume that $Q$ is fully normalized in $\F$ [4,~2.6] and therefore fully centralized too [4,~Proposition~2.11].  First of all, we prove that $p$ does not divide  $\vert\tilde\F (P,Q)_{\rm fc}\vert\,;$
 if $Q$ is $\F\-$selfcentralizing then we simply have $\tilde\F (P,Q)_{\rm fc} =\tilde\F (P,Q)$ and the claim 
 follows from [4,~6.7.2]. Assuming that $Q$ is {\it not\/}
$\F\-$selfcentralizing, it follows from [4,~statement~2.10.1] that any $\varphi\in \F (P,Q)_{\rm fc}$ can be extended to
an $\F\-$morphism $Q\.C_P (Q)\to P\,;$ that is to say, the inclusion $Q\i Q\.C_P (Q)$ induces a surjective map
$$\tilde\F\big(P,Q\.C_P(Q)\big)\too \tilde\F (P,Q)_{\rm fc}
\eqno £6.9.3\phantom{.}$$
and note that, since $Q\.C_P(Q)$ is $\F\-$selfcentralizing, it follows again from [4,~6.7.2] that $p$ does not divide
 $\vert\tilde\F\big(P,Q\.C_P(Q)\big)\vert\,.$
 \eject

 \smallskip
 Moreover, if $\psi,\psi'\in \F\big(P,Q\.C_P(Q)\big)$ have the same image in $\F (P,Q)_{\rm fc}$ then we have
 $\psi (Q) = \psi' (Q)$ and, in particular, we still have 
 $$\psi\big(Q\.C_P(Q)\big) = \psi'\big(Q\.C_P(Q)\big)
 \eqno £6.9.4$$
  so that
 $\psi' = \psi\circ \sigma$ for a suitable $\sigma\in \F\big(Q\.C_P(Q)\big)$ which stabilizes and acts trivially on $Q\,;$
 that is to say, $\sigma$ determines an element in $\big(C_\F (Q)\big)\big(C_P(Q)\big)$ [4,~2.14]. It easily follows that the $p'\-$group $\big(\widetilde{C_\F (Q)}\big)\big(C_P(Q)\big)$ [4,~Proposition~2.16] acts regularly on the ``fibers''
 of the map~£6.9.3 above and therefore we get
 $${\big\vert\tilde\F\big(P,Q\.C_P(Q)\big)\big\vert\over \big\vert\big(\widetilde{C_\F (Q)}\big)\big(C_P(Q)\big)\big\vert}
 = \vert \tilde\F (P,Q)_{\rm fc} \vert 
 \eqno £6.9.5\phantom{.}$$
 which proves our claim.

 \smallskip
 Secondly, it follows from [4,~Proposition~2.11] that $\varphi\in \F (P,Q)_{\rm fc}$ belongs to 
 $\F (P,Q)_{\rm fn}$ if and~only if $\F_P\big(\varphi (Q)\big)$ is a Sylow $p\-$subgroup of $\F\big(\varphi (Q)\big)\,;$ 
 thus, for any $\varphi\in \F (P,Q)_{\rm fc} -\F (P,Q)_{\rm fn}\,,$ setting $R = \varphi(Q)\,,$ the $p\-$group
$\F_P(R)$ is {\it not\/} a Sylow $p\-$subgroup of $\F(R)$ and therefore $p$ divides 
$\big\vert\bar N_{\F (R)}\big(\F_P(R)\big)\big\vert\,.$ But, since $R$ is fully centralized in $\F\,,$ 
it follows from [4,~statement~2.10.1]   that any $\sigma\in N_{\F (R)}\big(\F_P(R)\big)$ can be lifted to some 
$\tau\in \F\big(P,N_P(R)\big)$ such that $\tau (R) = R$ and therefore we have $\tau\big(N_P(R)\big)= N_P(R)\,;$ 
moreover, it is clear that if the image $\bar\sigma$ of $\sigma$ in   $\bar N_{\F (R)}\big(\F_P(R)\big)$ is not trivial 
then, denoting by $\varphi'\,\colon Q\to P$ the $\F\-$morphism mapping $u\in Q$ on $\tau\big(\varphi (u)\big)$ 
we have $\tilde\varphi'\not= \tilde\varphi\,.$
  Consequently, for any image $R$ of an element of $\F (P,Q)_{\rm fc} -\F (P,Q)_{\rm fn}\,,$ the group 
  $\bar N_{\F (R)}\big(\F_P(R)\big)$ acts freely on the set
  $$\{\tilde\varphi\in \tilde\F (P,Q)_{\rm fc} - \tilde\F (P,Q)_{\rm fn}\mid \varphi(Q) = R\}
  \eqno £6.9.6;$$
  hence, $p$ divides $\vert \tilde\F (P,Q)_{\rm fc} - \tilde\F (P,Q)_{\rm fn}\vert$ and therefore $p$ does not divide
  $\vert \tilde\F (P,Q)_{\rm fn}\vert\,.$ We are done.

  \bigskip
\noindent
{\bf £7. Basic elements in the double Burnside $\Bbb Z\-$algebra of $P$ }

\bigskip
£7.1. As a matter of fact, if $\H$ is the {\it Hecke $\Bbb Z\-$algebra\/} of a {\it Frobenius $P\-$category\/} $\F\,,$
the elements $f$ of $\H^\F$ with {\it length\/} (cf.~£5.3) prime to $p$ determine $\F$ and $\H\,,$ and they can be
described {\it a priori\/} in $\B_{P\times P}$ either {\it via\/} the {\it basic $P\times P\-$sets\/}
in [4,~Ch.~21] or by the {\it Frobenius condition\/} introduced in~[6]. In this section, we consider these
elements in our new context.

\medskip
£7.2. Let $f$ be an element of $\B_{P\times P}^{^{\rm pp}}\cap t(\B_{P\times P}^{^{\rm pp}})\,;$ according 
to Proposition~£5.9, it is quite clear  that $f$ belongs to a {\it Hecke $\Bbb Z\-$algebra\/} ---
namely to  the intersection $\B_{P\times P}^{^{\rm pp}}\cap t(\B_{P\times P}^{^{\rm pp}})$ --- and that the
intersection of two  {\it Hecke $\Bbb Z\-$algebras\/} is also a  {\it Hecke $\Bbb Z\-$algebra\/}; thus, it makes sense 
to consider the unique minimal  {\it Hecke $\Bbb Z\-$algebra\/} $\H_f$ containing $f\,.$ 
\eject

\bigskip
\noindent
{\bf Lemma~£7.3.} {\it With the notation above, $\H_f$ coincides with the minimal Hecke $\Bbb Z\-$algebra containing
the set of indecomposable elements  $g$ of $\B_{P\times P}^{^{\rm pp}}\cap t(\B_{P\times P}^{^{\rm pp}})$ such that 
$\vert g,f\vert \not= 0\,.$\/}
\medskip
\noindent
{\bf Proof:} According to condition~£5.9.1, all these elements $g$ belong to $\H_f\,.$
Conversely, let $\H$ be a Hecke $\Bbb Z\-$algebra which does not contain $f\,;$ thus, if we have $f = \sum_g z_g\.g$
where $g$ runs over the set of indecomposable elements  $g$ of~$\B_{P\times P}^{^{\rm pp}}\cap t(\B_{P\times P}^{^{\rm pp}})\,,$
there is such a $g$ such that $z_g\not= 0$ and $g\not\in \H\,;$ moreover, according to Proposition~£5.9, we can choose this indecomposable element in such a way that if $z_{g'}\not= 0$ then $\vert g,g'\vert = 0$ and therefore,  in this case, we get (cf.~£2.6.3)
$$\vert g,f\vert = z_g\,\vert g,g\vert \not= 0
\eqno £7.3.1.$$
We are done.

\medskip
£7.4. On the other hand, as in Proposition~£3.3, denoting by $\F^f (P,Q)$ the set of injective group homomorphisms 
$\varphi\,\colon Q\to P$ fulfilling
$${\rm res}_\varphi \.f = {\rm res}_{\iota_Q^P} \.f\qq {\rm res}_\varphi \.f^{^\circ} = {\rm res}_{\iota_Q^P} \.f^{^\circ}
\eqno £7.4.1,$$
the family of sets $\F^f (P,Q)$
 where $Q$ runs over the set of subgroups of $P$  determine a  {\it divisible $P\-$category
 $\F^f\,;$\/}
 indeed, if $\theta\,\colon R\to Q$ is a group homomorphism such that
$$\F^f (P,R)\cap \big(\F^f (P,Q)\circ \theta\big)\not= \emptyset
\eqno £7.4.2,$$
there is $\varphi\in \F^f (P,Q)$ such that $\varphi\circ\theta$ belongs to $\F^f (P,R)$ and therefore we~get 
$${\rm res}_{\iota_Q^P\circ\theta} \.f = {\rm res}_{\iota_R^P} \.f\qq {\rm res}_{\iota_Q^P\circ\theta} \.f^{^\circ} 
= {\rm res}_{\iota_R^P} \.f^{^\circ}
\eqno £7.4.3;$$ 
hence, for any $\varphi'\in \F^f (P,Q)$
we still get
$${\rm res}_{\varphi'\circ\theta} \.f = {\rm res}_{\iota_Q^P\circ\theta} \.f = {\rm res}_{\iota_R^P} \.f\qq
{\rm res}_{\varphi'\circ\theta} \.f^{^\circ}  = {\rm res}_{\iota_R^P} \.f^{^\circ}
\eqno £7.4.4.$$

\bigskip
\noindent
{\bf Lemma~£7.5.} {\it For any element $f$ of $\B_{P\times P}^{^{\rm pp}}\cap 
t(\B_{P\times P}^{^{\rm pp}})$ and any subgroup $Q$ of~$P$ we have
$$\sum_{\tilde\varphi\in \tilde\F_{\H_f}(Q,P)} \big\vert f_{\iota_Q^P,\varphi},f\big\vert 
\big/\big\vert C_P\big(\varphi (Q)\big)\big\vert \equiv \ell (f) \bmod{p} 
\eqno £7.5.1.$$
In particular, if $\,\ell (f)$ is prime to $P$ and $\F_{\H_f}$ is contained in $\F^f\,,$ 
 we have the equality $\F_{\H_f}=\F^f\,.$\/}

\medskip
\noindent
{\bf Proof:} Since we have $\vert f_{\iota_Q^P,\varphi},f\vert = \vert {\rm res}_\varphi,{\rm res}_{\iota_Q^P}\.f\vert$ 
(cf.~equality~£5.7.1), let us consider the decomposition   ${\rm res}_{\iota_Q^P}\.f = \sum_g z_g\. g$  
of ${\rm res}_{\iota_Q^P}\.f$  in the $\Bbb Z\-$basis of {\it indecomposable\/} elements~$g$ of~$\B_{Q\times P}\,;$  
for any $\varphi\in {\rm Hom}(Q,P)\,,$ it is easily checked that  $\vert {\rm res}_\varphi,g\vert = 0$ unless 
$g = {\rm res}_\varphi$ (cf.~£2.6.3) and then we have 
$$\vert {\rm res}_\varphi,{\rm res}_\varphi\vert = \big\vert C_P\big(\varphi (Q)\big)\big\vert
\eqno £7.5.2;$$
\eject
\noindent
hence, since $p$ divides $\ell (g)$ unless $g = {\rm res}_\varphi$ for some $\varphi\in {\rm Hom}(Q,P)$ and then we have 
$\ell ({\rm res}_\varphi) = 1\,,$ it is quite clear that
$$\sum_{\tilde\varphi\in \tilde\F_{\H_f}(Q,P)} \big\vert f_{\iota_Q^P,\varphi},f\big\vert 
\big/\big\vert C_P\big(\varphi (Q)\big)\big\vert = \!\!\!\sum_{\tilde\varphi\in \tilde\F_{\H_f}(Q,P)} z_{{\rm res}_\varphi} \equiv \ell (f)\bmod{p} 
\eqno £7.5.3.$$

\smallskip
In particular,  if $\ell (f)$ is prime to $P\,,$  for  any subgroup $Q$ of~$P$  there is $\varphi\in \F_{\H_f}(P,Q)$ such that 
$\vert f_{\iota_Q^P,\varphi},f\vert\not= 0\,;$ hence, if $\F_{\H_f}(P,Q)\i \F^f(P,Q)$ then for any 
$\psi\in \F^f(P,Q)$ it follows from isomorphism~£5.7.1 and equalities~£6.3.1 that we have
$$\vert f_{\iota_Q^P,\psi},f\vert =\vert f_{\iota_Q^P,\varphi},f\vert\not= 0
\eqno £7.5.4;$$
now, it follows from statement~£3.2.1 that $ f_{\iota_Q^P,\psi}$ is also an {\it indecomposable\/} element of $\H_f$
and therefore $\psi$ also belongs to $\F_{\H_f}(P,Q)\,.$ We are done.

\bigskip
\noindent
{\bf Proposition~£7.6.} {\it Let $\F$ be a divisible $P\-$category and set $\H = \H_\F\,.$ If $\H^\F$ contains
an element $f$ such that $\ell (f)$ is prime to~$p\,,$ then $\F$ is a
Frobenius $P\-$category and we have $\F^f =\F$ and $\H_f = \H\,.$\/}

\medskip
\noindent
{\bf Proof:} Since $f$ belongs to $\H\,,$ $\H$ contains $\H_f$ and therefore $\F_{\H_f}$ is contained in $\F\,;$
similarly, it is clear that $\F$ is contained in $\F^f\,;$ since  $\ell (f)$ is prime to~$p\,,$ it follows from~Lemma~£7.5
that  we have the equalities
$$\F_{\H_f} = \F = \F^f
\eqno £7.6.1\phantom{.}$$
and therefore we still have $H_f = \H\,.$

\smallskip
Moreover, it follows from Proposition~£6.5 that we have $f = f' - f''$ for suitable {\it positive\/} elements $f'$ and $f''$
 of~$\H^\F\,;$ in particular, $\H$ contains $\H_{f'}$ and~$\H_{f''}\,,$ and $\F$ is contained in $\F^{f'}$
 and $\F^{f''}\,.$ We actually may assume that  $\ell (f')$ is prime to $p\,;$ then, since $f'$ is {\it positive\/}, 
 it is the isomorphism class of a functor $\frak f'$ in $\frak H_\F\,;$ we already know that this functor comes from a 
 suitable $P\times P\-$set  $\Omega'$ where $\{1\}\times P$ and $P\times \{1\}$ act freely and that $p$ does not divide~$\ell (f') = 
\vert \Omega'\vert/\vert P\vert\,;$ moreover, since $f'$ belongs to~$\H^\F\,,$ it follows from Lemma~£6.4 that
we have $f'^{^\circ} = f'$ and therefore we still have $\Omega'^\circ \cong \Omega'\,.$

\smallskip
That is to say, the  
$P\times P\-$set $\Omega'$ fulfills all the conditions in [4,~21.2]; moreover, for any subgroup $Q$ of $P\,,$
in [4,~21.3] we denote by ${\rm Hom}^{\Omega'} (Q,P)$ the set of injective group homomorphisms 
$\varphi\,\colon Q\to P$ fulfilling (cf.~£2.6.3)
$$\Nat (\frak f_{\iota_Q^P,\varphi},\frak f') \not= \emptyset
\eqno £7.6.3;$$
then, since $\frak f'$ is a functor in $\frak H_\F\,,$ it follows from statement~£3.2.1 that $\frak f_{\iota_Q^P,\varphi}$
is also a functor in $\frak H_\F$ and, in particular, that $\varphi$ belongs to $\F (P,Q)\,;$
consequently, since $\F\i \F^{f'} = \F^{^{\Omega'}}\,,$ the inclusions in~[4,~21.5.1] hold and therefore 
$\Omega'$ is a {\it basic $P\times P\-$set\/}, we have $\F = \F^{f'} $ and it follows from [4,~Proposition~21.9]
that $\F^{f'}$ is a {\it Frobenius $P\-$category\/}. We are done.
\eject

\medskip
£7.7. Now, we say that $f\in \B_{P\times P}$ is {\it basic\/} if $f$ belongs to $\B_{P\times P}^{^{\rm pp}}\,,$ 
$\ell (f)$ is prime to $p$ and we have $\H_f = \H_{\F^f}$ or, equivalently, $\F^f = \F_{\!\H_f}\,;$ thus, according
to Lemma~£7.5, $f\in \B_{P\times P}^{^{\rm pp}}$  is {\it basic\/} if $\ell (f)$ is prime to $p$ and $\F^f$ contains
$\F_{\H_f}\,.$ Note that, if $f$ is {\it basic\/} then
it belongs to $(\H_f)^{\F_{\!\H_f}}$ and therefore it follows from Lemma~£6.4 that we have $f^{^\circ} = f\,;$
moreover, according to Proposition~£7.6, in this case $\F^f$ is a {\it Frobenius $P\-$category\/}. Conversely,
it follows from isomorphism~£6.6.2 and Proposition~£7.6 that it $\F$ is a {\it Frobenius $P\-$category\/} then
$(\H_\F)^\F$ contains {\it basic\/} elements.

\medskip
£7.8. In [6], Kari Ragnarsson and Radu Stancu find an equivalent definition of a {\it basic\/} element {\it via\/} the
{\it Frobenius condition\/}; we say that $f\in \B_{P\times P}$ fulfills  the {\it Frobenius condition\/} if, for any
$s, s'\in \B_P\,,$ we have
$$ f\big(f(s)\,s'\big) = f(s)f(s')
\eqno £7.8.1\phantom{.}$$
or, equivalently, if $f$  centralizes $m^P\big(f(\B_P)\big)\,;$
as we mention in Remark~£6.7 above, it follows from Theorem~£6.6 that any {\it basic\/} element fulfills the
{\it Frobenius condition\/}. Following Kari Ragnarsson and Radu Stancu, in order to discuss  this condition,
for any pair of functors $\frak f$ and $\frak g$ in $\frak B_{P\times P}\,,$ we consider the functors  (cf.~£2.3)
$$\frak f\times \frak g :  \frak B_P \times \frak B_P\too \frak B_P\qq 
\frak g\circ (\frak f\times \id) : \frak B_P \times \frak B_P\too \frak B_P
\eqno £7.8.2;$$
moreover, for any subgroup $Q$ of $P$ and any injective group homomorphisms $\varphi,\varphi'\in {\rm Hom} (P,Q)\,,$ we consider the {\it indecomposable\/} functor
$$\frak f_{\iota_Q^P,\varphi,\varphi'} = \ind_Q^P\circ (\res_\varphi\times \res_{\varphi'})
 : \frak B_P \times \frak B_P\too \frak B_P
 \eqno £7.8.3;$$
 with this notation, in~[6~Lemma~7.5] they prove the following result.

\bigskip
\noindent
{\bf Proposition~£7.9.} {\it With the notation above, for any subgroup $Q$ of $P$ and  any injective group homomorphisms 
$\varphi,\varphi'\in {\rm Hom} (P,Q)\,,$ we have natural bijections
$$\eqalign{\Nat (\frak f_{\iota_Q^P,\varphi,\varphi'}, \frak f\times\frak g)&\cong \Nat (\frak f_{\iota_Q^P,\varphi} , \frak f)
\times \Nat (\frak f_{\iota_Q^P,\varphi'} , \frak g)\cr
\Nat \big(\frak f_{\iota_Q^P,\varphi,\varphi'}, \frak g\circ (\frak f\times \id)\big)&\cong \Nat (\frak f_{\varphi',\varphi} , \frak f)
\times \Nat (\frak f_{\iota_Q^P,\varphi'} , \frak g)\cr}
\eqno £7.9.1$$\/}

\par
\noindent
{\bf Proof:} Since the restriction preserves the {\it trivial\/} actions, it follows from bijection~£2.3.3 and from~£2.4  that we have
$$\eqalign{\Nat (\frak f_{\iota_Q^P,\varphi,\varphi'}, \frak f\times \frak g) &\cong 
\Nat \big(\res_\varphi\times \res_{\varphi'},(\res_{Q}^{P}\circ\frak f)\times 
(\res_{Q}^{P}\circ\frak g)\big)\cr
&\cong \Nat (\res_\varphi,\res_{Q}^{P}\circ\frak f)\times \Nat (\res_{\varphi'}, \res_{Q}^{P}\circ\frak g)\cr
&\cong \Nat (\frak f_{\iota_Q^P,\varphi} , \frak f)
\times \Nat (\frak f_{\iota_Q^P,\varphi'} , \frak g)}
\eqno £7.9.2\,.$$

\smallskip
Once again, it follows from bijection~£2.3.3 that we have
 $$\Nat \big(\frak f_{\iota_Q^P,\varphi,\varphi'}, \frak g\circ (\frak f\times \id)\big) \cong
\Nat \big(\res_\varphi\times \res_{\varphi'}, \res_Q^P\circ \frak g\circ (\frak f\times \id)\big)
\eqno £7.9.3.$$
But, in order to prove the bottom bijection in~£7.9.1, it suffices to consider the case
where $\frak g = \frak f_{\iota_R^P,\psi}$ for a subgroup $R$ of $P$ and an injective group homomorphism 
$\psi\,\colon R\to P\,;$ then, choosing a set of representatives $W$ for $R\backslash P/Q$ and, for any $w\in W\,,$ setting 
$U_w= Q\cap R^w$ and $\psi_w = \psi\circ\kappa_w$ where  $\kappa_w\,\colon U_w\to R$ denotes the group homomorphism determined by the conjugation 
by $w\,,$ we get
$$\eqalign{\res_Q^P\circ \frak f_{\iota_R^P,\psi}\circ (\frak f\times \id)
&= (\res_Q^P\circ\ind_R^P)\circ \big((\res_\psi\circ\frak f)\times \res_\psi\big)\cr
&= \big(\bigsqcup_{w\in W} \ind_{U_w}^Q\circ \res_{\kappa_w}\big)\circ 
\big((\res_\psi\circ\frak f)\times \res_\psi\big)\cr
&= \bigsqcup_{w\in W} \ind_{U_w}^Q\circ \big((\res_{\psi_w}\circ\frak f)\times \res_{\psi_w}\big)\cr}
\eqno £7.9.4\,;$$
moreover, it is easily checked that for any $w\in W$ such that $U_w\not= Q$ we have
$$\Nat \Big(\res_\varphi\times \res_{\varphi'}, \ind_{U_w}^Q\circ \big((\res_{\psi_w}\circ\frak f)\times \res_{\psi_w}\big)\Big)
= \emptyset
\eqno £7.9.5.$$

\smallskip
Consequently, setting $\hat W = W\cap T_P(R,Q)\,,$ once again from~£2.4 we get
$$\eqalign{\Nat \big(&\frak f_{\iota_Q^P,\varphi,\varphi'}, \frak f_{\iota_R^P,\psi}\circ (\frak f\times \id)\big)\cr 
&\cong \bigsqcup_{w\in \hat W} \Nat \big(\res_{\varphi}\times \res_{\varphi'}, 
(\res_{\psi_w}\circ\frak f)\times \res_{\psi_w}\big)\cr
&\cong \bigsqcup_{w\in \hat W} \Nat (\res_{\varphi},\res_{\psi_w}\circ\frak f)\times  
\Nat (\res_{\varphi'}, \res_{\psi_w})\cr
&\cong \bigsqcup_{w\in \hat W} \Nat (\frak f_{\psi_w,\varphi},\frak f)\times  
\Nat (\frak f_{\psi_w,\varphi'}, \id)\cr}
\eqno £7.9.6;$$
moreover, the set $\Nat (\frak f_{\psi_w,\varphi'}, \id)$ is empty unless $\psi_w = \varphi'$ which at most happens once.
On the other hand, we also have
$$\eqalign{\Nat (\frak f_{\iota_Q^P,\varphi'},\frak f_{\iota_R^P,\psi})
&\cong \Nat (\res_{\varphi'}, \res_Q^P\circ \ind_R^P\circ\res_\psi)\cr
&\cong  \bigsqcup_{w\in W} \Nat (\res_{\varphi'},\ind_{U_w}^Q\circ \res_{\psi_w})\cr
&\cong \bigsqcup_{w\in \hat W} \Nat (\res_{\varphi'}, \res_{\psi_w})\cong 
\bigsqcup_{w\in \hat W} \Nat (\frak f_{\psi_w,\varphi'}, \id)\cr}
\eqno £7.9.7.$$
Now, the bottom bijection in~£7.9.1 follows from bijections~£7.9.6 and~£7.9.7.
\eject

\bigskip
\noindent
{\bf Corollary~£7.10.} {\it With the notation above, any element $f$ in $\B_{P\times P}^{^{\rm pp}}\cap 
t(\B_{P\times P}^{^{\rm pp}})$ fulfills the Frobenius condition only if, for any subgroup $Q$ of $P$ and  any injective group homomorphisms $\varphi,\varphi'\in {\rm Hom} (P,Q)\,,$ we have
$$\vert f_{\iota_Q^P,\varphi} ,  f\vert\,\vert f_{\iota_Q^P,\varphi'} ,  f\vert = 
\vert f_{\varphi',\varphi} ,  f\vert\, \vert f_{\iota_Q^P,\varphi'} , f\vert
\eqno £7.10.1.$$\/}

\par
\noindent
{\bf Proof:} With the notation in the proposition above, respectively denote by $f\,,$ $g$ and $f_{\iota_Q^P,\varphi,\varphi'}$
the isomorphism classes of $\frak f\,,$ $\frak g$ and $\frak f_{\iota_Q^P,\varphi,\varphi'}\,;$ it is clear that 
$f_{\iota_Q^P,\varphi,\varphi'}$ becomes an element of $\B_{P\times P\times P}$ (cf.~£5.4) and that $\frak f\times \frak g$
 and $\frak f\times \id$ also determine elements in  $\B_{P\times P\times P}\,,$ respectively denoted by $f\times g$ and $f\times 1\,;$
 thus, the bijections in~£7.9.1 imply the equalities  (cf.~£5.4)
$$\eqalign{\vert f_{\iota_Q^P,\varphi,\varphi'},  f\times g\vert& =  \vert f_{\iota_Q^P,\varphi} ,  f\vert
\, \vert f_{\iota_Q^P,\varphi'} ,  g\vert\cr
\vert f_{\iota_Q^P,\varphi,\varphi'},  g\. ( f\times 1)\vert& = \vert f_{\varphi',\varphi} , f\vert
\, \vert f_{\iota_Q^P,\varphi'} ,  g\vert\cr}
\eqno £7.10.2;$$
since all the members of these equalities are {\it bi-additive\/} on $f$ and $g\,,$ it follows from our definition of the
{\it partial scalar product\/} in~£5.5 that these equalities remain true  for any $f$ and any $g$ in $\B_{P\times P}\,.$

\smallskip
On the other hand,  an element $f$ in $\B_{P\times P}^{^{\rm pp}}\cap  t(\B_{P\times P}^{^{\rm pp}})$ fulfilling the Frobenius condition clearly fulfills the equality $f\.(f\times 1) = f\times f\,;$ then, equality~£7.10.1 follows from equalities~£7.10.2.

\bigskip
\noindent
{\bf Theorem~£7.11.}  {\it Assume that $f\in \B_{P\times P}^{^{\rm pp}}$ fulfills the Frobenius condition, that we have $f^{^\circ} = f$ 
and that $\ell (f)$ is prime to $p\,.$ Then $f$ is  basic.\/}
\medskip
\noindent
{\bf Proof:} According to Lemma~£7.5 and to~£7.7, it suffices to prove that $\F^f$ contains $\F_{\H_f}$ or equivalently
that, for any subgroup $Q$ of $P\,,$  any  injective group homomorphism $\varphi\,\colon Q\to P$ fulfilling 
$\vert f_{\iota_Q^P,\varphi},f\vert \not= 0$ belongs to~$\F^f (P,Q)\,.$ Thus, let us denote by ${\rm Hom}^f (Q,P)$
the set of  injective group homomorphisms $\varphi\,\colon Q\to P$ such that $\vert f_{\iota_Q^P,\varphi},f\vert \not= 0\,,$ 
and by ${\rm Aut}^f (Q)$ the stabilizer of ${\rm Hom}^f (Q,P)$ in ${\rm Aut} (Q)\,;$ that is to say, $\sigma \in {\rm Aut} (Q)$
belongs to ${\rm Aut}^f (Q)$ if and only if we have $\vert f_{\iota_Q^P,\varphi\circ \sigma},f\vert \not= 0$ for any 
$\varphi\in {\rm Hom}^f (Q,P)$ and note that $ f_{\iota_Q^P,\varphi\circ \sigma} =  f_{\iota_Q^P,\varphi}$
is equivalent to $\sigma\in {}^{\varphi*}\!\F_P \big(\varphi (Q)\big)\,.$ In particular, since  we are assuming that
 $f^{^\circ} = f$ and then from equality~£7.10.1 we get
$$\vert f_{\iota_Q^P,\varphi} ,  f\vert = \vert f_{\varphi',\varphi} ,  f\vert = \vert f_{\varphi,\varphi'} ,  f\vert
= \vert f_{\iota_Q^P,\varphi'} ,  f\vert
\eqno £7.11.1\phantom{.}$$
for any $\varphi,\varphi'\in {\rm Hom}^f (Q,P)\,,$ it follows from Lemma~£7.5 above that, if $\iota_Q^P$ belongs to 
${\rm Hom}^f (Q,P)$ and $\vert N_P(Q)\vert$ is maximal in the set of integers $\big\vert N_P\big(\zeta (Q)\big)\big\vert$
when $\zeta$ runs over ${\rm Hom}^f (Q,P)\,,$ we have
$$\vert {\rm Aut}^f (Q)\,\colon \F_P (Q)\vert\,\vert f_{\iota_Q^P,\zeta} ,  f\vert\big/\vert C_P (Q)\vert\not\equiv 0 \bmod{p} 
\eqno £7.11.2.$$

\smallskip
First of all we claim that, for any subgroup $Q$ of $P$ such that $\iota_Q^P$ belongs to ${\rm Hom}^f (Q,P)\,,$
and any $\sigma\in {\rm Aut}^f (Q)\,,$ we have
$${\rm res}_{\iota_Q^P\circ\sigma}\. f =  {\rm res}_{\iota_Q^P}\.f
\eqno £7.11.3;$$
since we are assuming that  $f^{^\circ} = f \,,$ this equality is equivalent to the equalities (cf.~bijection~£2.6.3 and~£5.6)
$$\vert f_{\sigma_R,\psi},f\vert = \vert f_{\iota_R^P,\psi},f\vert
\eqno £7.11.4\phantom{.}$$ 
where $R$ runs over the set of subgroups of $Q\,,$  $\sigma_R\,\colon R\to P$ denotes the group homomorphism 
determined by $\sigma\,,$ and $\psi$ runs over the set of  injective group homomorphism from $R$ to $P\,;$ hence, 
according to the corresponding equality~£7.10.1, it suffices to prove that $\vert f_{\iota_R^P,\sigma_R},f\vert\not= 0$ and,
since $\iota_Q^P\circ\sigma$ belongs to ${\rm Hom}^f (Q,P)\,,$  we may assume that
$R\not= Q\,.$

\smallskip
As above,  choose $\psi\in {\rm Hom}^f (R,P)$ such that $\big\vert N_P\big(\psi (R)\big)\big\vert$ is maximal;
it follows from~£7.11.2 above that $\F_P \big(\psi (R)\big)$ is a Sylow $p\-$subgroup of 
${\rm Aut}^f \big(\psi (R)\big)$ and that we have 
$$\vert f_{\iota_{\psi (R)}^P,\iota_{\psi (R)}^P},f\vert \big/\vert C_P \big(\psi (R)\big) \vert\not\equiv 0 \bmod{p} 
\eqno £7.11.5;$$
then, it follows again from equality~£7.10.1 that, up to modifying  our choice of $\psi$ with  a suitable 
$\tau\in {\rm Aut}^f \big(\psi (R)\big)\,,$ we may assume that
$${}^\psi \F_P (R)\i \F_P\big(\psi (R)\big)
\eqno £7.11.6;$$
in particular, denoting by $\psi^*$ the inverse of the isomorphism $R\cong \psi (R)$ determined by $\psi\,,$
 ${}^{\psi^*}\!\tilde\F_P\big(\psi (R)\big)$ contains  the group
$$\tilde\F_Q (R) =  {}^{\sigma^{\!-1}}\!\tilde\F_Q\big(\sigma (R)\big)\cong \bar N_Q (R)/\bar C_Q (R)
\eqno £7.11.7.$$

\smallskip
The {\it indecomposable\/} elements of $\B_{Q\times P}^{^{\rm pp}}$ are clearly of the form (cf.~£5.4)
$$g_{\theta,\eta} = {\rm ind}_\theta\.{\rm res}_\eta
\eqno £7.11.8\phantom{.}$$
where  $\theta$ is an injective group homomorphism to $Q$ from a subgroup $T$ of $Q$ and $\eta\in {\rm Hom} (T,P)\,;$
let us call {\it length\/} of $g_{\theta,\eta}$ the index $\vert Q\colon T\vert\,.$
Then, in the decomposition of $g = {\rm res}_{\iota_Q^P}\.f$ in this $\Bbb Z\-$basis denote by $g'$ the sum
of all the terms with  {\it length\/} strictly smaller than $\vert Q\colon R\vert\,,$ and set $g'' = g - g'\,.$

\smallskip
Now, in order to prove equalities~£7.11.4, we argue by induction on $\vert Q\colon R\vert\,;$ that is to say, for
any subgroup $T$ of $Q$ such that $\vert Q\colon T\vert < \vert Q\colon R\vert$ and any injective group
homomorphism $\eta\,\colon T\to P\,,$ we may assume that
$$\vert f_{\sigma_T,\eta},f\vert = \vert f_{\iota_T^P,\eta},f\vert
\eqno £7.11.9;$$
in this case, denoting by $\sigma_T^Q\,\colon T\to Q$ the restriction of $\sigma_T\,,$ it is easily checked that we still have 
(cf. bijection~£2.6.3)
$$\vert g_{\sigma_T^Q,\eta},g'\vert = \vert g_{\iota_T^Q,\eta},g'\vert
\eqno £7.11.10\phantom{.}$$
for  any  injective group homomorphism $\eta\,\colon T\to P\,.$

\smallskip
Then, according to~£5.6, we get ${\rm res}_\sigma\. g' = g'$
which implies that
$$\vert g_{\sigma_R^Q,\psi},g'\vert = \vert g_{\iota_R^Q,\psi},g'\vert\qq {\rm res}_{\iota_Q^P\circ\sigma}\. f = g' + {\rm res}_\sigma\. g''
\eqno £7.11.11;$$
in particular, it is clear that for suitable $x,y\in \Bbb Z$ we have
$$\eqalign{\vert f_{\sigma_R,\psi},f\vert &= \vert g_{\iota_R^Q,\psi},g'\vert
+ x\vert g_{\sigma_R^Q,\psi},g_{\sigma_R^Q,\psi}\vert\cr
\vert f_{\iota_R,\psi},f\vert &= \vert g_{\iota_R^Q,\psi},g'\vert
+ y\vert g_{\iota_R^Q,\psi},g_{\iota_R^Q,\psi}\vert\cr}
\eqno £7.11.12.$$
But, we already know that (cf.~£3.6.1 and~£5.5.1)
$$\eqalign{\vert g_{\sigma_R^Q,\psi},g_{\sigma_R^Q,\psi}\vert &= 
\big\vert\bar N_{Q\times P} \big(\Delta_{\sigma_R^Q,\psi}(R)\big)\big\vert\cr
\big\vert g_{\iota_R^Q,\psi},g_{\iota_R^Q,\psi}\vert &=\vert \bar N_{Q\times P}\big(\Delta_{\iota_R^Q,\psi}(R)\big)\big\vert\cr}
\eqno £7.11.13;$$
moreover, note that we have group isomorphims
$$\eqalign{N_{Q\times P} \big(\Delta_{\sigma_R^Q,\psi}(R)\big)\big/C_Q\big(\sigma (R)\big)&\times C_P\big(\psi(R)\big)\cr
 &\cong {}^{\sigma^{\!-1}}\!\F_Q\big(\sigma (R)\big) \cap{}^{\psi^*}\F_P\big(\psi(R)\big)\cr
N_{Q\times P}\big(\Delta_{\iota_R^Q,\psi}(R)\big)\big/C_Q (R)&\times C_P\big(\psi(R)\big)\cr
&\cong \F_Q (R)\cap {}^{\psi^*}\F_P\big(\psi (R)\big)\cr}
\eqno £7.11.14.$$
Hence, since $N_Q (R)\not= R$ and ${}^{\psi^*}\!\tilde\F_P\big(\psi (R)\big)$ contains the image of the quotient
$N_Q (R)/C_Q (R)$ (cf.~£7.11.7), $p$ divides the integers
$$\vert g_{\sigma_R^Q,\psi},g_{\sigma_R^Q,\psi}\vert
\big/\big\vert C_P\big(\psi(R)\big)\big\vert\qq \big\vert g_{\iota_R^Q,\psi},g_{\iota_R^Q,\psi}\vert
\big/\big\vert C_P\big(\psi(R)\big)\big\vert
\eqno £7.11.15.$$

\smallskip
Consequently, it follows from congruence~£7.11.5 and from the bottom equality in~£7.11.12 that 
 $$\vert g_{\iota_R^Q,\psi},g'\vert \big/\big\vert C_P\big(\psi (R)\big)\big\vert\not\equiv 0 \bmod{p} 
\eqno £7.11.16;$$
then, it follows from the top equality in~£7.11.12 that 
 $$\vert f_{\sigma_R,\psi},f\vert \big/\big\vert C_P\big(\psi (R)\big)\big\vert\not\equiv 0 \bmod{p} 
\eqno £7.11.17\phantom{.}$$
and, in particular, that $\vert f_{\sigma_R,\psi},f\vert = \vert f_{\psi,\sigma_R},f\vert$ is not zero; finally, equality~£7.10.1 yields
$$\vert f_{\iota_R^P,\sigma_R},f\vert\vert f_{\iota_R^Q,\psi},f\vert 
= \vert f_{\psi,\sigma_R},f\vert\vert f_{\iota_R^Q,\psi},f\vert 
\eqno £7.11.18\phantom{.}$$ 
which proves our claim since $\psi$ belongs to ${\rm Hom}^f (R,P)\,.$
\eject

\smallskip
On the other hand, since $\ell (f)$ is prime to~$p\,,$ it follows from Lemma~£7.5 that there is $\sigma\in {\rm Aut}(P)$ such that 
$$\vert f_{{\rm id}_P,\sigma},f\vert/\vert Z(P)\vert \not\equiv 0 \bmod{p}
\eqno £7.11.19\phantom{.}$$
 and therefore, according to equality~£7.10.1, for any $\tau\in {\rm Aut}(P)$ we still have $\vert f_{{\rm id}_P,\tau},f\vert =  \vert f_{\sigma,\tau},f\vert$ and, since $f_{\sigma,\sigma} = f_{{\rm id}_P,{\rm id}_P}\,,$ we get
$$\vert f_{{\rm id}_P,{\rm id}_P},f\vert/\vert Z(P)\vert = \vert f_{{\rm id}_P,\sigma},f\vert/\vert Z(P)\vert  \not\equiv 0 \bmod{p}
\eqno £7.11.20;$$
thus,  ${\rm id}_P$ belongs to ${\rm Hom}^f (P,P)$ and then equality~£7.10.1 implies that 
$${\rm Aut}^f (P) = {\rm Hom}^f (P,P)
\eqno £7.11.21.$$
Hence, it follows from our claim that $\F^f (P)$ contains ${\rm Hom}^f (P,P)\,.$

\smallskip
From now on, in order to prove that any $\varphi\in {\rm Hom}^f (Q,P)$ belongs to~$\F^f(P,Q)\,,$ we argue 
by induction on $\vert P\colon Q\vert$ and may assume that $Q\not= P\,.$ First of all, we assume  that 
$\big\vert N_P\big(\varphi (Q)\big)\big\vert$ is maximal and then it follows from congruence~£7.11.2 that 
$\F_P\big(\varphi (Q)\big)$ is a Sylow $p\-$subgroup of ${\rm Aut}^f\big(\varphi (Q)\big)$ and that we have
 $$\vert f_{\iota_{{\rm id}_{\varphi (Q)}}^P,\iota_{{\rm id}_{\varphi (Q)}}^P},f \vert \big/\big\vert C_P\big(\varphi (Q)\big)\big\vert\not\equiv 0 \bmod{p} 
\eqno £7.11.22.$$
Note that it suffices to prove that $\tau\circ\varphi$ belongs to $\F^f (P,Q)$ for some element $\tau$
in ${\rm Aut}^f \big(\varphi (Q)\big)\,;$ indeed, in this case $\iota_{{\rm id}_{\varphi (Q)}}^P$ belongs to 
${\rm Hom}^f \big(\varphi (Q),P\big)$ and it suffices to apply our claim above to $\varphi (Q)$ to get that 
$${\rm res}_{\iota_{\varphi (Q)}^P\circ\tau}\. f =  {\rm res}_{\iota_{\varphi (Q)}^P}\.f
\eqno £7.11.23,$$
so that $\varphi$ belongs to $\F^f (P,Q)$ too. Thus, we may assume that
$${}^\varphi \F_P (Q)\i \F_P\big(\varphi (Q)\big)
\eqno £7.11.24.$$

\smallskip
 We know that $f = f' - f''$ where $f'$ and $f''$ are the isomorphism classes of functors $\frak f_{\Omega'}$ and 
$\frak f_{\Omega''}$  for some finite $P\times P\-$sets $\Omega'$ and $\Omega''\,,$ and that we have (cf.~£3.6.1 and~£5.5.1)
$$\vert f_{\iota_Q^P,\varphi},f\vert = \vert \Omega'^{\Delta_{\iota_Q^P,\varphi}(Q)}\vert
- \vert \Omega''^{\Delta_{\iota_Q^P,\varphi}(Q)}\vert
\eqno £7.11.25;$$
 moreover, since  $f^{^\circ} = f$ and  $f$ belongs to $\B_{P\times P}^{^{\rm pp}}\,,$
 we may assume that $\Omega'^\circ\cong \Omega'$ and $\Omega''^\circ\cong \Omega''\,,$ and that $\{1\}\times P$ acts freely
 on these $P\times P\-$sets. Then, on the one hand, the group $N = N_{P\times P}\big(\Delta_{\iota_Q^P,\varphi} (Q)\big)$ acts on the quotient sets 
$$\Omega'^{\Delta_{\iota_Q^P,\varphi}(Q)}\!\big/ C_P\big(\varphi (Q)\big)\qq
\Omega''^{\Delta_{\iota_Q^P,\varphi}(Q)}\!\big/ C_P\big(\varphi (Q)\big)
\eqno £7.11.26;$$ 
on the other hand, since $f_{\iota_{{\rm id}_{\varphi (Q)}}^P,\iota_{{\rm id}_{\varphi (Q)}}^P} = f_{\varphi,\varphi}\,,$
it follows from equalities~£7.11.1 and congruence~£7.11.22 that we have
$$\vert f_{\iota_Q^P,\varphi},f\vert \big/\big\vert C_P\big(\varphi (Q)\big)\big\vert\not\equiv 0 \bmod{p} 
\eqno £7.11.27.$$
\eject

\smallskip
Consequently, we still have
$$\Big\vert\Big( \Omega'^{\Delta_{\iota_Q^P,\varphi}(Q)}\!\big/ C_P \big(\varphi (Q)\big)\Big)^N\Big\vert \not\equiv  
\Big\vert\Big( \Omega''^{\Delta_{\iota_Q^P,\varphi}(Q)}\!\big/ C_P \big(\varphi (Q)\big)\Big)^N\Big\vert  \bmod{p} 
\eqno £7.11.28\phantom{.}$$
and, in particular, the difference between both members is not zero; moreover, if $C_P \big(\varphi(Q)\big)\.\omega$ is a
fixed point of $N$ in the quotient set
$$\Big(\Omega'^{\Delta_{\iota_Q^P,\varphi}(Q)}\bigsqcup \Omega'^{\Delta_{\iota_Q^P,\varphi}(Q)}
\Big)\Big/ C_P\big(\varphi (Q)\big)
\eqno £7.11.29\phantom{.}$$
where $\omega\in \Omega'^{\Delta_{\iota_Q^P,\varphi}(Q)}\bigsqcup \Omega'^{\Delta_{\iota_Q^P,\varphi}(Q)}\,,$ then 
the stabilizer of~$\omega$ in~$N$ has the form  $\Delta_{\iota_T^P,\eta} (T)$  for some 
subgroup $T$ of $N_P (Q)$ containing $Q$ and a suitable   injective group homomorphism $\eta\,\colon T\to P\,.$ extending~$\varphi\,,$ and we necessarily have
$$N \cong C_P \big(\varphi(Q)\big)\rtimes\Delta_{\iota_T^P,\eta} (T)
\eqno £7.11.30.$$

\smallskip
In conclusion,   for some  subgroup $T$ of $N_P (Q)$ containing $Q$ and a suitable   injective group homomorphism $\eta\,\colon T\to P\,.$ 
extending~$\varphi\,,$ we still get
$$\big\vert\Omega'^{\Delta_{\iota_T^P,\eta}(T)}\big/C_P\big(\eta (T)\big)\big\vert\not\equiv 
\big\vert\Omega''^{\Delta_{\iota_T^P,\eta}(T)}\big/C_P\big(\eta (T)\big)\big\vert\bmod{p} 
\eqno £7.11.31;$$
in particular, we have 
$$\vert f_{\iota_T^P,\eta},f\vert = \vert \Omega'^{\Delta_{\iota_T^P,\eta}(T)}\vert
- \vert \Omega''^{\Delta_{\iota_T^P,\eta}(T)}\vert \not= 0
\eqno £7.11.32,$$ 
so that $\eta$  belongs to ${\rm Hom}^f (T,P)\,.$ On the other hand, it is easily checked that inclusion~£7.11.24 
implies that $p$ divides $\big\vert N/C_P \big(\varphi(Q)\big)\big\vert\,,$
which forces $Q\not= T\,;$ hence,  by the induction hypothesis,  $\eta$ belongs to~$\F^f(P,T)\,;$
thus, since $\F^f$ is divisible and $\eta$ extends $\varphi\,,$ $\varphi$ belongs to $\F^f(P,Q)\,.$

\smallskip
Finally, for any $\psi\in {\rm Hom}^f (Q,P)\,,$ from equality~£7.10.1 we obtain
$$0\not= \vert f_{\iota_Q^P,\varphi},f\vert= \vert f_{\psi,\varphi,},f\vert= \vert f_{\iota_{\psi (Q)}^P,\varphi\circ \psi^*},f\vert
\eqno £7.11.33\phantom{.}$$
and therefore $\varphi\circ \psi^*$ also belongs to  ${\rm Hom}^f \big(\psi (Q),P\big)\,;$ but, the cardinal
$$\big\vert N_P\big(\varphi (Q)\big)\big\vert = \Big\vert N_P\Big((\varphi\circ \psi^*)\big)\psi (Q)\big)\Big)\Big\vert
\eqno £7.11.34\phantom{.}$$
 is still maximal; consequently, the argument above proves that
$\varphi\circ \psi^*$ actually belongs to $\F^f \big(\psi (Q),P\big)\,;$ once again, the divisibility of $\F^f$
implies that $\psi$ belongs to $\F^(P,Q)\,.$ We are done.
\eject

 \bigskip
 \noindent
 {\bf £8. The case of a finite group} 
\bigskip
£8.1. Let $G$ be a finite group admitting $P$ as a Sylow $p\-$subgroup and consider
the {\it Frobenius $P\-$category\/} $\F = \F_G\,.$ On the one hand, we have introduced above the
{\it Hecke $\Bbb Z\-$algebra\/}~$\H_\F$ of $\F\,.$ On the other hand, the 
$\Bbb Z\-$algebra
$$\dot\H_G= {\rm End}_{\Bbb Z G}\big({\rm Ind}_{N_G (P)}^G (\Bbb Z)\big)
\eqno £8.1.1\phantom{.}$$
is usually called the {\it Hecke $\Bbb Z\-$algebra\/} of $G\,.$ In this section, our purpose is to relate $\H_\F$ with the 
{\it extended  Hecke $\Bbb Z\-$algebra\/}
of $G\,,$ namely with the $\Bbb Z\-$algebra
$$\hat\H_G = {\rm End}_{\Bbb Z G}\big({\rm Ind}_P^G (\Bbb Z)\big)
\eqno £8.1.2.$$
Actually, the $\Bbb ZG\-$module ${\rm Ind}_P^G (\Bbb Z)$ is isomorphic to the {\it permutation\/}
$\Bbb ZG\-$mo-dule with $\Bbb Z\-$basis $G/P$ and therefore the corresponding images of $1\otimes 1$ in ${\rm Ind}_P^G (\Bbb Z)$
determines a canonical $\Bbb Z\-$module isomorphism
$$\hat\H_G\cong \big(\Bbb Z (G/P)\big)^P\cong \Bbb Z (P\backslash G/P)
\eqno £8.1.3;$$
more precisely, in the group $\Bbb Q\-$algebra $\Bbb Q G$ setting $h_D = {1\over\vert P\vert}\.\sum_{x\in D} x$
for any $D\in P\backslash G/P\,,$ it is easily checked that $\bigoplus_{D\in P\backslash G/P} \Bbb Z\.h_d$ is a $\Bbb Z\-$subalgebra
of~$\Bbb Q G$ {\it canonically\/} isomorphic to $\hat\H_G\,.$

\medskip
£8.2.  In order to describe this relationship, let us consider the 
{\it transporter category\/} $\T = \T_G$ of $G$ [17,~17.2], namely the category where
the objects are all the subgroups of $P$ and, for any pair of subgroups $Q$ and $R$
of $P\,,$ we set
$$\T (Q,R) = \{x\in G \mid R\i Q^x\}
\eqno £8.2.1,$$
the composition of $\T\-$morphisms being induced by the product in $G\,;$
note that we have a canonical functor
$$\frak c : \T\too \F
\eqno £8.2.2.$$
Recall that the {\it exterior\/} quotient $\tilde \T$ of $\T$ is the quotient category defined by the inner
automorphisms of the objects [4,~1.3], and that the {\it additive cover\/} $\ad (\tilde\T)$ of $\tilde\T$
is the category where the objects are the finite sequences $\bigoplus_{i\in I} Q_i$ of $\tilde\T\-$objects
$Q_i$ and the morphisms are the pairs formed by a map between the sets of indices and by a sequence
of $\tilde\T\-$morphisms between the corresponding $\tilde\T\-$objects [4,~A2.7.3].
\eject

\bigskip
\noindent
{\bf Lemma~£8.3.} {\it The additive cover $\ad (\tilde\T)$ of the exterior quotient
 $\,\tilde\T$ of $\,\T$ admits pull-backs.\/}
 
\medskip
\noindent
{\bf Proof:}  Let $Q\,,$ $R$ and $T$ be three subgroups of $P$ and consider two $\tilde\T\-$mor-phisms 
$$\tilde x^{^Q}\,\colon R\too Q\qq \tilde y^{^Q}\,\colon T\too Q
\eqno £8.3.1\phantom{.}$$ 
where $\tilde x^{^Q}$ and $\tilde y^{^Q}$ are the respective classes of t$x\in \T (Q,R)$ and $y\in \T (Q,T)\,;$
choose a set of representatives $W\i Q$ for ${}^x R\backslash Q/ {}^y T$ and set $U = \bigoplus_{w\in W} U_w$ 
where $U_w = {}^x R\cap  {}^{wy} T$ for any $w\in W\,.$  Then, since $x^{-1}$ and $(wy)^{-1}$ determine 
$\tilde\T\-$morphisms 
$$\widetilde{x^{-1}}^{^R} : U_w \too R\qq  \widetilde{(wy)^{-1}}^{\!^T} : U_w \too T
\eqno £8.3.2\phantom{.}$$
for any $w\in W\,,$ we get two obvious $\ad (\tilde\T)\-$morphisms
$$a : U\too R\qq b : U\too T
\eqno £8.3.3\phantom{.}$$
and we claim that the $\ad (\tilde\T)\-$diagram
$$\matrix{&&\hskip-10pt Q \hskip-10pt&&\cr
&\hskip-5pt{\tilde x^{^Q} \atop }\hskip-10pt\nearrow&
&\nwarrow\hskip-4pt{\tilde y^{^Q} \atop }\hskip-10pt\cr
R\hskip-10pt&&&&\hskip-15pt T\cr
&{\atop a}\hskip-4pt\nwarrow&&\nearrow\hskip-4pt{\atop b}&\cr
&&\hskip-13pt U \hskip-13pt&&\cr}
\eqno £8.3.4\phantom{.}$$
is a {\it pull-back.\/}

\smallskip
 Indeed, since $W\i Q\,,$ this diagram is clearly commutative; moreover,
let $V$ be an $\ad (\tilde\T)\-$object and 
$$c : V\too R\qq d : V\too T
\eqno £8.3.5\phantom{.}$$
 two $\ad (\tilde\T)\-$morphisms fulfilling $\tilde x^{^Q}\circ c 
 = \tilde y^{^Q}\circ d\,;$ in order to prove that there is a unique  
 $\ad (\tilde\T)\-$morphism $h\,\colon V\to U$ fulfilling $c\circ h = a$ and 
 $d\circ h = b\,,$  it is quite clear that we may assume that $V$ is a subgroup
 of $P$ and that we have $c = \tilde r^{^R}$ and $d = \tilde t^{^T}$ for suitable 
 $r\in \T (R,V)$  and $t\in \T (T,V)\,.$

\smallskip
 In this case, we have $\widetilde{xr}^{^Q}
 = \widetilde{yt}^{^Q}$ and therefore there is $w\in Q$ such that $xr = wyt\,;$
 actually, up to modifying the choices of $x$ in $\tilde x^{^Q}$ and of $y$ 
 in~$\tilde y^{^Q}\,,$ we may assume that $w$ belongs to $W$ and then we have
 $${}^{xr} V =  {}^{wyt} V \i U_w
 \eqno £8.3.6,$$
  so that we can define $h\,\colon V\to U$ as the $\ad (\tilde\T)\-$morphism formed by
 the map to~$W$ from the set with an element corresponding to $V$ which maps
 this element on $w\,,$ and by the $\tilde\T\-$morphism
 $$\widetilde{xr}^{^{U_w}} = \widetilde{wyt}^{^{U_w}}  : V\too U_w
 \eqno £8.3.7.$$
 \eject
 \noindent
 Finally, any $\ad (\tilde\T)\-$morphism $h'\,\colon V\to U$  fulfilling 
 $c\circ h' = a$ and  $d\circ h' = b$ determines some $w'\in W$ and a 
 $\tilde\T\-$morphism
 $$\widetilde{z}^{^{U_{\!w'}}}  : V\too U_{w'}
 \eqno £8.3.8\phantom{.}$$
 fulfiiling $\tilde{r}^{^R} = \widetilde{x^{-1}z}^{^R}$ and $\tilde{ t}^{\,^T} = \widetilde{(w'y)^{-1}z}^{\!\!\!^T}\,;$ once again, we may assume that 
 $xr = z = w'yt$ which forces $w' = w$ and $h' = h\,.$ This proves our claim.

\smallskip
Now, if $Q = \bigoplus_{i\in I} Q_i\,,$ $R = \bigoplus_{j\in J} R_j$ and 
$T = \bigoplus_{\ell\in L} T_\ell$ are three $\ad(\tilde\T)\-$objects and we consider
two $\ad(\tilde\T)\-$morphisms
$$R\too Q\qq T\too Q
\eqno £8.3.9\phantom{.}$$
given by  two maps $f\,\colon J\to I$ and $g\,\colon L\to I\,,$ and by 
$\tilde\T\-$morphisms
$$\tilde x_j^{^Q}\,\colon R_j\too Q_{f(j)}\qq 
\tilde y_\ell^{^Q}\,\colon T_\ell\too Q_{g(\ell)}
\eqno £8.3.10,$$ 
it is easily checked that, considering the {\it pull-backs\/}
$$\matrix{&&Q_{f(j)} &&\cr
&{\tilde x_j^{^Q} \atop }\hskip-10pt\nearrow&
&\hskip-4pt\nwarrow\hskip-4pt{\tilde y_\ell^{^Q} \atop }\cr
R_j\hskip-20pt  &&&&\hskip-20pt T_\ell\cr
&\hskip-5pt{\atop a_{(j,\ell)}}\hskip-8pt\nwarrow&&\nearrow\hskip-6pt{\atop b_{(j,\ell)}}&\cr
&&U_{(j,\ell)}&&\cr}
\eqno £8.3.11\phantom{.}$$
for any $(j,\ell)\in J\times_I L\,,$ the ``sum''  $U = \bigoplus_{(j,\ell)\in J\times_I L} U_{(j,\ell)}$ together with the  $\ad(\tilde\T)\-$morphisms
$$a : U\too R\qq b : U\too T
\eqno £8.3.12\phantom{.}$$
given by the projections $J\times_I L\to J$ and $J\times_I L\to L$
and by the $\tilde\T\-$morphisms
$$a_{(j,\ell)} : U_{(j,\ell)}\too R_j\qq b_{(j,\ell)} : U_{(j,\ell)}\too T_\ell
\eqno £8.3.13\phantom{.}$$
define a {\it pull-back\/} of the  $\ad(\tilde\T)\-$morphisms~£8.3.9. We are done.

\medskip
£8.4. Let $\frak H_\T$ be the category where the objects $\frak t_{a,a'}$ are the
triples formed  by an $\ad(\tilde\T)\-$object $Q$ and a pair of 
$\ad(\tilde\T)\-$morphisms $a$ and $a'$ from $Q$ to~$P$ --- to be more explicit, we also say that $\frak t_{a,a'}$ 
is an object {\it over\/} $Q$  --- and where, for a pair of $\ad(\tilde\T)\-$objects $Q$ and $R\,,$ and a  pair of pairs of 
$\ad(\tilde\T)\-$morphisms 
$$\matrix{P&&P\cr
{\atop a}\hskip-5pt\nwarrow\hskip-10pt&
&\hskip-10pt\nearrow\hskip-5pt{\atop a'}\cr
&Q}\qq \matrix{P&&P\cr
{\atop b}\hskip-5pt\nwarrow\hskip-10pt&
&\hskip-10pt\nearrow\hskip-5pt{\atop b'}\cr
&R}
\eqno £8.4.1,$$ 
a morphism from $\frak t_{a,a'}$ to $\frak t_{b,b'}$ is an $\ad(\tilde\T)\-$morphism $h\,\colon R\to Q$ fulfilling
$a\circ h = b$ and $a'\circ h = b'\,,$ the composition being induced by the composition in~$\ad(\tilde\T)\,.$~Then, 
denoting by $\nabla\,\colon P\oplus P\to P$ the $\ad(\tilde\T)\-$morphism determined by the identity ${\rm id}_P\,,$ we define
$$\eqalign{a \nabla b = \nabla\circ (a\oplus b) &: Q\oplus R\too P\oplus P\too P\cr
\frak t_{a,a'}\oplus \frak t_{b,b'} &= \frak t_{a \nabla b\,,a' \nabla b'}\cr}
\eqno £8.4.2\phantom{.}$$
\eject
\noindent
and this {\it direct sum\/} is clearly {\it associative\/}; coherently, we say that $\frak t_{a,a'}$ is {\it indecomposable\/} 
whenever $Q$ is a subgroup of $P$ and, in this case, that it is {\it maximal\/} if any  $\frak H_\T\-$morphism from 
$\frak t_{a,a'}$ to an {\it indecomposable\/} $\frak H_\T\-$object is an isomorphism. Moreover, from any pair of 
$\ad(\tilde\T)\-$morphisms
$$h : \frak t_{a,a'}\too \frak t_{\bar a,\bar a'}\qq 
\ell : \frak t_{b,b'}\too \frak t_{\bar b,\bar b'}
\eqno £8.4.3,$$
we clearly obtain a new $\ad(\tilde\T)\-$morphism
$$h\oplus \ell : \frak t_{a,a'}\oplus \frak t_{b,b'}\too \frak t_{\bar a,\bar a'}\oplus \frak t_{\bar b,\bar b'}
\eqno £8.4.4.$$

\medskip
£8.5. On the other hand, choosing a {\it pull-back\/}
$$\matrix{&&\hskip-10pt P \hskip-10pt&&\cr
&{a' \atop }\hskip-10pt\nearrow&
&\hskip-5pt\nwarrow\hskip-4pt{b \atop }\cr
Q\hskip-10pt&&&&\hskip-20pt R\cr
&{\atop \hat b}\hskip-4pt\nwarrow&&\nearrow\hskip-4pt{\atop \hat a'}&\cr
&&\hskip-13pt Q\times_P R \hskip-13pt&&\cr}
\eqno £8.5.1,$$
we define the {\it tensor product\/} as the following $\frak H_\T\-$object over $Q\times_P R$ 
$$\frak t_{a,a'}\otimes \frak t_{b,b'} = \frak t_{a\circ \hat b\,,b'\circ \hat a'}
\eqno £8.5.2;$$
as above, from any pair of $\ad(\tilde\T)\-$morphisms
$$h : \frak t_{a,a'}\too \frak t_{\bar a,\bar a'}\qq 
\ell : \frak t_{b,b'}\too \frak t_{\bar b,\bar b'}
\eqno £8.5.3,$$
we clearly obtain a new $\ad(\tilde\T)\-$morphism
$$h\otimes \ell : \frak t_{a,a'}\otimes \frak t_{b,b'}\too \frak t_{\bar a,\bar a'}\otimes \frak t_{\bar b,\bar b'}
\eqno £8.5.4.$$

\medskip
£8.6. Moreover, for a third $\frak H_\T\-$object $\frak t_{c,c'}\,,$ we have a {\it natural\/}
isomorphism
$$(\frak t_{a,a'}\otimes \frak t_{b,b'})\otimes \frak t_{c,c'}\cong 
\frak t_{a,a'}\otimes (\frak t_{b,b'}\otimes \frak t_{c,c'})
\eqno £8.6.1\phantom{.}$$
and it is easily checked that the corresponding {\it pentagonal $\frak H_\T\-$diagram\/} is commutative;
similarly, the  {\it tensor product\/} is {\it distributive\/} with respect to the {\it direct sum\/} in the sense that we have natural isomorphisms
$$\eqalign{(\frak t_{a,a'}\oplus \frak t_{b,b'})\otimes \frak t_{c,c'}&\cong 
(\frak t_{a,a'}\oplus  \frak t_{c,c'})\otimes (\frak t_{b,b'}\oplus \frak t_{c,c'})\cr
\frak t_{c,c'}\otimes (\frak t_{a,a'}\oplus \frak t_{b,b'})&\cong 
(\frak t_{c,c'}\oplus  \frak t_{a,a'})\otimes (\frak t_{c,c'}\oplus \frak t_{b,b'})\cr}
\eqno £8.6.2.$$
Note that the functor $\frak c\,\colon \T\to \F$ above induces new functors
$$ \tilde\frak c : \tilde\T\too \tilde\F\qq
\ad (\tilde \frak c) : \ad(\tilde\T)\too \ad( \tilde\F)
\eqno £8.6.3,$$
and finally a {\it surjective\/} functor 
$$\frak H_\frak c : \frak H_\T\too \frak H_\F
\eqno £8.6.4\phantom{.}$$ 
which is compatible with the {\it direct sum\/} and maps the  {\it tensor product\/} on the composition of functors;
moreover, this functor preserves {\it indecomposable\/} and {\it maximal indecomposable\/} objects.

\bigskip
\noindent
{\bf Proposition~£8.7.} {\it With the notation above, for any indecomposable 
$\frak H_\T\-$ob-ject $\frak t_{a,a'}$ over $Q$ there exists a unique isomorphism class of maximal  
indecomposable $\frak H_\T\-$objects $\frak t_{\hat a,\hat a'}\,,$ over a suitable subgroup $\hat Q$ of $P\,,$
such that  $\frak H_\T (\frak t_{a,a'},\frak t_{\hat a,\hat a'})$ is not empty, and then we have
$$\frak H_\T (\frak t_{a,a'},\frak t_{\hat a,\hat a'})= \{\tilde r^{^{\hat Q}}\}
\eqno £8.7.1.$$
In particular, for any $\frak H_\T\-$object $\frak t_{b,b'}$ over $R$ admitting an $\frak H_\T\-$morphism $\tilde s^{^R}\,\colon 
\frak t_{a,a'}\to \frak t_{b,b'}\,,$ there exists a unique $\tilde t^{^{\hat Q}}\in
\frak H_\T (\frak t_{b,b'},\frak t_{\hat a,\hat a'})$ fulfilling 
$\widetilde{ts}^{^{\hat Q}}= \tilde r^{^{\hat Q}}\,.$\/}

\medskip
\noindent
{\bf Proof:} We argue by induction on $\vert P\,\colon Q\vert$ and may assume that $Q\not= P\,;$ respectively denote by $x,x'\in \T (P,Q)$ representatives of the $\tilde\T\-$morphisms  $a,a'$ and set
$$N = N_{P^x\cap P^{x'}} (Q)
\eqno £8.7.2;$$
then, it is clear that $x$ and $x'$ belong to $\T (P,N)$ and that, setting 
$c = \tilde x^{^P}$~and $c' = \tilde x'^{^P}\,,$ the inclusion $Q\i N$ determines
an $\frak H_\T\-$morphism $h\,\colon \frak t_{a,a'}\to \frak t_{c,c'}\,;$
moreover, $\frak t_{c,c'}$ admits an $\frak H_\T\-$morphism 
$\hat h\,\colon \frak t_{c,c'}\to \frak t_{\hat c,\hat c}$ to a maximal  
inde-composable $\frak H_\T\-$object $\frak t_{\hat c,\hat c}$ over $\hat N\,;$ if $\hat c$ and 
$\hat c'$ belong to $\tilde\T (P,\hat N)\,,$ denote by~$\tilde r^{^{\hat N}}$ the corresponding composition of natural maps
$$\frak t_{a,a'}\buildrel h\over\too \frak t_{c,c'}
\buildrel \hat h\over\too \frak t_{\hat c,\hat c'}
\eqno £8.7.3.$$

\smallskip
Now, if  $\tilde s^{^R}\,\colon \frak t_{a,a'}\to \frak t_{b,b'}$ is an  $\frak H_\T\-$morphism where 
$b = \tilde y^{^P}$ and $b' = \tilde y'^{^P}$ belong to $\tilde\T (P,R)\,,$ there are $u,u'\in P$ such that 
$$x = uy s\qq x' = u'y' s
\eqno £8.7.4;$$
set $M = N_R (Q^{s^{-1}})$ and respectively denote by $d\,\colon M\to P$
and $d'\,\colon M\to P$ the composition of $c$ and $c'$ with the 
$\tilde\T\-$morphism $\hat h'$ determined by the inclusion $M\i R\,;$ it is clear that
$\tilde s^{^R}$ coincides with the composition
$$\frak t_{a,a'}\buildrel h'\over\too\frak t_{d,d'}\buildrel \hat h'\over\too
\frak t_{b,b'}
\eqno £8.7.5\phantom{.}$$
where $h'$ is also determined by $s\,.$ In particular, since 
$R\i P^{^y}\cap P^{^{y'}}\,,$ we have
$$M^s = N_{R^s}(Q)\i N_{P^{ys}\cap P^{y's}} (Q)\ = N
\eqno £8.7.6\phantom{.}$$
and thus $s^{-1}$ induces an $\frak H_\T\-$morphism
$\ell\,\colon \frak t_{d,d'}\to \frak t_{c,c'}$ fulfilling~$\ell\circ h' = h\,.$
\eject

\smallskip
At this point, it follows from the  induction hypothesis that  there is a unique 
$\tilde t^{^{\hat N}}\in\frak H_\T (\frak t_{b,b'},\frak t_{\hat c,\hat c'})$ fulfilling 
$\tilde t^{^{\hat N}}\circ \hat h' = \hat h\circ \ell$ and therefore we get
$$\tilde t^{^{\hat N}}\circ \tilde s^{^R} = \tilde t^{^{\hat N}}\circ \hat h'\circ h' 
=  \hat h\circ \ell \circ h'  = \hat h\circ h = \tilde r^{^{\hat N}}
\eqno £8.7.7;$$
this already proves that $\frak t_{\hat c,\hat c'}$  is the unique isomorphism class of maximal  indecomposable $\frak H_\T\-$objects such that 
$\frak H_\T (\frak t_{a,a'},\frak t_{\hat c,\hat c'})$ is not empty.
On the other hand, if $\tilde t'^{^{\hat N}}\,\colon \frak t_{b,b'}\to 
\frak t_{\hat c,\hat c'}$ is another $\frak H_\T\-$morphism fulfilling 
$\tilde t'^{^{\hat N}}\circ \tilde s^{^R} = \tilde r^{^{\hat N}}\,,$
we have 
$$\tilde t'^{^{\hat N}}\circ \hat h'\circ h' = \tilde t'^{^{\hat N}}\circ \tilde s^{^R} = \tilde r^{^{\hat N}} = \hat h\circ h  = \hat h\circ \ell \circ h' 
\eqno £8.7.8\phantom{.}$$
which clearly implies $\tilde t'^{^{\hat N}}\circ \hat h' =  \hat h\circ \ell$ and then
the uniqueness of $\tilde t^{^{\hat N}}$ forces the equality $\tilde t'^{^{\hat N}} =\tilde t^{^{\hat N}}\,.$ We are done.

\medskip
£8.8. Let us denote by $\H_\T$ the free $\Bbb Z\-$module over the set of isomorphism classes $t_{a,a'}$ of {\it indecomposable\/} $\frak H_\T\-$objects
$\frak t_{a,a'}$ --- called {\it indecomposable ele-ments\/} of $\H_\T\,;$ we say that
the {\it indecomposable element\/} $t_{a,a'}$ is {\it maximal\/} whenever the
corresponding {\it indecomposable\/} $\frak H_\T\-$object $\frak t_{a,a'}$ is so.
It is clear that any $\frak H_\T\-$object determines a {\it positive\/}  element of $\H_\T$ and that the {\it direct sum\/} of 
$\H_\T\-$objects corresponds to the sum in $\H_\T\,;$
similarly, the {\it tensor\/} product of $\frak H_\T\-$objects induces an {\it associative and distributive product\/} in $\H_\T\,,$ so that $\H_\T$ becomes a $\Bbb Z\-$algebra. It is clear that the functor $\frak H_\frak c\,\colon \frak H_\T\to \frak H_\F$
above induces a {\it surjective\/} $\Bbb Z\-$algebra homomorphism 
$$\H_\frak c : \H_\T\too \H_\F
\eqno £8.8.1.$$

\bigskip
\noindent
{\bf Proposition~£8.9.} {\it  For any indecomposable element  $t$ in $\H_\T $ there is a  unique {\it maximal\/} indecomposable element $\hat t$ fulfilling 
$\big\vert \H_\frak c (t),\H_\frak c (\hat t)\big\vert\not= 0\,.$ Denote by
$e_\T\,\colon \H_\T\to \H_\T$ the $\Bbb Z\-$module endomorphism mapping  any indecomposable element  $t\in \H_\T$ on  $\displaystyle{\ell (\H_\frak c (t))\over\ell (\H_\frak c (\hat t))}\.\hat t\,.$ If $\H_\frak c (t)\not\in \N_\F$ then 
$$e_\F \big(\H_\frak c (t)\big) = \H_\frak c \big(e_\T (t)\big)
\eqno £8.9.1\phantom{.}$$
 and, for any pair of elements $r$ and $s$ of~$\H_\T\,,$ we have 
$$e_\T (rs) = e_\T \big(e_\T (r)e_\T (s)\big)
\eqno £8.9.2.$$
In particular, ${\rm Ker}(e_\T)$ is a two-sided ideal of $\H_\T\,.$\/}

\medskip
\noindent
{\bf Proof:} It follows from Proposition~£8.7 that, if $t$ is the isomorphism class of
an $\frak H_\T\-$object $\frak t\,,$ there is a unique isomorphism class $\hat t$ of {\it maximal indecomposable\/} 
$\H_\T\-$objects $\hat \frak t$ admitting  an $\frak H_\T\-$morphism $h\,\colon \frak t\to \hat\frak t$  or, equivalently, fulfilling 
$\big\vert \H_\frak c (t),\H_\frak c (\hat t)\big\vert\not= 0\,.$ Moreover, since 
$\H_\frak c (\hat t)$ is a {\it maximal indecomposable\/} element in $\H_\F\,,$  if 
$\H_\frak c (t)$ does not belong to $\N_\F$ then equality~£8.9.1
follows from the very definition  of~$e_\F$ (cf.~Proposition~£5.12).
\eject

 \smallskip
 Consider a second {\it indecomposable\/} element  $s$ in $\H_\T\,,$ coming from an {\it indecomposable\/} $\frak H_\T\-$object $\frak s\,,$ and the corresponding isomorphism class $\hat s$ of {\it maximal  indecomposable\/} $\frak H_\T\-$objects  $\hat\frak s$ admitting an $\frak H_\T\-$morphism  $\ell\,\colon \frak s\to 
\hat\frak s\,;$ then, we also have an $\frak H_\T\-$morphism (cf.~£8.5)
$$h\otimes \ell : \frak t\otimes \frak s\too \hat\frak t\otimes \hat\frak s
\eqno £8.9.3.$$
More explicitly, assume that $\frak t\,,$ $\hat\frak t\,,$
$\frak s$ and $\hat\frak s$ respectively come from subgroups $Q\,,$ $\hat Q\,,$
$R$ and $\hat R$ of $P$ and from $x$ and  $x'$ in $\T (P,Q)\,,$ $\hat x$ and
$\hat x'$ in $\T (P,\hat Q)\,,$  $y$ and  $y'$ in $\T (P,R)\,,$ and $\hat y$ and
$\hat y'$ in $\T (P,\hat R)\,;$ actually,  we may assume that $Q\i \hat Q\,,$ 
that $R\i \hat R\,,$ that $\hat x = x\,,$ $\hat x' = x'\,,$ $\hat y = y$ and
$\hat y' = y'\,,$ and that $h$ and $\ell$ are both determined by the trivial element.

\smallskip
Then, we still may assume that the $\ad (\tilde\T)\-$morphism
$$\matrix{h\otimes \ell : &Q\times_P R&\too &\hat Q\times_P\hat R\cr
&\Vert &&\Vert\cr
&\bigoplus_{w\in W} U_w&&\bigoplus_{\hat w\in \hat W} \hat U_{\hat w}\cr}
\eqno £8.9.4\phantom{.}$$
comes from the canonical map 
$${}^{x'}\! Q\backslash P /\,{}^y \!R\too {}^{x'}\!\hat Q\backslash P /\,{}^y\! \hat  R
\eqno £8.9.5\phantom{.}$$
 and  from the family of $\tilde\T\-$morphisms
$$\matrix{\widetilde{{}^{x'}\! x''}^{^{\hat U_{\hat w}}} : 
&{}^{x'}\! Q\cap \,{}^{wy} \!R&\too  & {}^{x'}\!\hat Q\cap \,{}^{\hat w y}\! \hat  R\cr
\Vert&\Vert&\phantom{\big\uparrow}&\Vert\cr
(h\otimes \ell)_w &U_w&&\hat U_{\hat w}\cr}
\eqno £8.9.6\phantom{.}$$
where, borrowing notation  from~Lemma~£8.3 above, we have chosen  sets of representatives $\hat W$ for  
${}^{x'}\!\hat Q\backslash P /\,{}^y\! \hat R\,,$ $X$ for $\hat Q/Q$ and $Y$ for $\hat R/R\,,$ and set
$$w =  ({}^{x'}\! x'')^{-1}\.\hat w\.\,{}^y y''
\eqno £8.9.7$$ 
for  $\hat w$ running over $\hat W\,,$ $x''$ over $X$ and $y''$ over $Y\,,$ so that $W = ({}^{x'}\!X)^{-1}\.\hat W\.\,{}^y Y$ 
is a set of representatives for ${}^{x'}\! Q\backslash P /\,{}^y \!R\,.$

\smallskip
 Consequently, for any $\hat w\in \hat W\,,$ any $x''\in X$ and any $y''\in Y\,,$ setting
  $w =  ({}^{x'}\! x'')^{-1}\.\hat w\.\,{}^y y''$ and
$$\frak r_w = \frak t_{\widetilde{x x'^{-1}}^{^P}\widetilde{,y' y^{-1}w^{-1}}^{\!\!\!^P}}\qq \frak r_{\hat w} 
= \frak t_{\widetilde{x x'^{-1}}^{^P}\widetilde{,y' y^{-1}\hat w^{-1}}^{\!\!\!^P}}
\eqno £8.9.8,$$ 
we have the $\frak H_\T\-$morphism $(h\otimes \ell)_w\,\colon\frak r_w\to 
\frak r_{\hat w}\,;$ finally, the $\frak H_\T\-$morphism~£8.9.3 becomes
$$h\otimes \ell = \bigoplus_{w\in W} (h\otimes \ell)_w : \frak s\otimes \frak t
= \bigoplus_{w\in W} \frak r_w \too \hat \frak s\otimes \hat \frak t
= \bigoplus_{\hat w\in \hat W} \frak r_{\hat w}
\eqno £8.9.9.$$
Moreover, it follows from Proposition~£8.7 that, for any $\hat w\in \hat W\,,$ 
if $\hat\frak r_{\hat w}$ is a {\it maximal indecomposable\/} $\frak H_\T\-$object
admitting an $\frak H_\T\-$morphism from $\frak r_{\hat w}\,,$ it is also  a {\it maximal indecomposable\/} $\frak H_\T\-$object admitting an $\frak H_\T\-$morphism from~$\frak r_w\,.$
 
\smallskip
Thus, denoting by $\hat r_{\hat w}$ the isomorphism class of 
$\hat\frak r_{\hat w}\,,$ by the very definition of $e_\T$ we get
 $$e_\T (s t) = \sum_{w\in W} {\vert P \colon U_w\vert\over\ell \big(\H_\frak c 
 (\hat r_{\hat w})\big)}\.\hat r_{\hat w}\qq 
 e_\T(\hat s\hat t) = \sum_{\hat w\in \hat W} {\vert P \colon \hat U_{\hat w}\vert\over\ell \big(\H_\frak c  (\hat r_{\hat w})\big)}\.\hat r_{\hat w}
 \eqno £8.9.10;$$
hence, we still get
$$\eqalign{e_\T (st) &= \sum_{\hat w\in \hat W}
\Big(\sum_{w\in ({}^{x'}\!X)^{-1}\.\hat w\.\,{}^y Y} 
{\vert P \colon U_w\vert\over\ell \big(\H_\frak c  (\hat r_{\hat w})\big)}\Big)\.\hat r_{\hat w}\cr 
 e_\T \big(e_\T ( s)e_\T ( t)\big) 
&= \vert \hat Q\colon Q\vert\vert \hat R\colon R\vert 
\sum_{\hat w\in \hat W} {\vert P \colon \hat U_{\hat w}\vert 
\over \ell \big(\H_\frak c  (\hat r_{\hat w})\big)}\.\hat r_{\hat w}\cr}
\eqno £8.9.11;$$
moreover, for any $\hat w\in \hat W\,,$ it is clear that
$$\eqalign{\sum_{w\in ({}^{x'}\!X)^{-1}\.\hat w\.\,{}^y Y} \vert P \colon U_w\vert 
&= \sum_{w\in ({}^{x'}\!X)^{-1}\.\hat w\.\,{}^y Y} {\vert P\vert\vert {}^{x'}Q\.w\.\,{}^y R\vert\over \vert Q\vert \vert R\vert}\cr
& = {\vert P\vert\vert {}^{x'}\hat Q\.\hat w\.\,{}^y \hat R\vert\over \vert Q\vert \vert R\vert} 
=  \vert \hat Q\colon Q\vert\vert \hat R\colon R\vert \vert P \colon \hat U_{\hat w}\vert\cr}
\eqno £8.9.12;$$
finally, we obtain $ e_\T \big(e_\T ( s)e_\T ( t)\big) = e_\T (st)\,.$
We are done.

\medskip
£8.10. On the other hand, we have a $\Bbb Z\-$module homomorphism (cf.~£8.1)
$$d_G : \H_\T\too \hat \H_G \i \Bbb Q G
\eqno £8.10.1\phantom{.}$$
mapping the isomorphism class $t_{a,a'}$ of an {\it indecomposable\/} 
$\frak H_\T\-$object $\frak t_{a,a'}\,,$ where $a = \tilde x^{^Q}$ and $a' = \tilde x'^{^Q}\,,$ on 
$$d_G (t_{a,a'}) =  {\vert P\vert\over \vert P^{x \,x'^{-1}}\cap P\vert}
\. h_P\,x \,x'^{-1}h_P = h_{Px x'^{-1}P} 
\eqno £8.10.2;$$ 
moreover, for an {\it indecomposable\/} $\frak H_\T\-$object $\frak t_{b,b'}$ admitting an
$\frak H_\T\-$morphism $h\,\colon \frak t_{b,b'}\to \frak t_{a,a'}\,,$
it is easily checked that, denoting by $t_{b,b'}$ the isomorphism class of 
$\frak t_{b,b'}\,,$ we have~$d_G (t_{b,b'}) = d_G  (t_{a,a'})\,.$

\bigskip
\noindent
{\bf Proposition~£8.11.} {\it The $\Bbb Z\-$module homomorphism $d_G$
induces a $\Bbb Z\-$algebra isomorphism
$$\H_T/{\rm Ker}(e_\T)\cong \hat \H_G
\eqno £8.11.1\phantom{.}$$\/}
\eject

\par
\noindent
{\bf Proof:} Since $(e_\T)^2 = e_\T$ and for any {\it indecomposable\/} element
$t$ of $\H_\T$ we have $d_G(t) = d_G\big(e_\T(t)\big)\,,$ it is clear that $d_G$
induces a $\Bbb Z\-$module homomorphism
$$\H_T/{\rm Ker}(e_\T)\too  \hat \H_G
\eqno £8.11.2;$$
moreover, since any $x\in G$ determines a $\tilde \T\-$morphism 
$\tilde x^{^P}\colon 1\to P\,,$ considering the {\it indecomposable\/} element
$t_{\tilde x^{^P},\tilde 1^{^P}}$ of $\H_\T\,,$ we have
$$d_G (t_{\tilde x^{^P},\tilde 1^{^P}}) =  {\vert P\vert\over \vert P^{x }\cap P\vert}\. h_P\, x\, h_P = h_{PxP}
\eqno £8.11.3\phantom{.}$$
and therefore the  $\Bbb Z\-$module homomorphism~£8.11.2 is surjective.

\smallskip
On the other hand, for any subgroup $Q$ of $P$ and any $x,x'\in \T (P,Q)\,,$
 it is clear that both $\tilde\T\-$morphisms 
 $$\tilde x^{^P} : Q\too P\qq \tilde x'^{^P} : Q\too P
 \eqno £8.11.4\phantom{.}$$
can be factorized throughout 
$$\tilde x'^{^{P^{\,xx'^{-1}}\cap P}} :  Q\too P^{\,xx'^{-1}}\cap P
\eqno £8.11.5;$$
it easily follows that the  $\Bbb Z\-$module homomorphism~£8.11.2 is injective too.

\smallskip
It remains to prove that $d_G$ is also a $\Bbb Z\-$algebra homomorphism;
let $Q$ and $R$ be subgroups of $P\,,$ $x$ and $x'$ elements of $\T (P,Q)$
and $y$ and $y'$  elements of~$\T (P,R)\,;$ as above, choose a set of 
representatives $W\i P$ for ${}^{x'}\! Q\backslash P /\,{}^y \!R\,;$ then, 
for any $w\in W\,,$  we have the commutative $\tilde\T\-$diagram
$$\matrix{&&\hskip-10pt P \hskip-10pt&&\cr
&{\tilde x'^{^P} \atop }\hskip-10pt\nearrow&
&\hskip-5pt\nwarrow\hskip-4pt{\tilde y^{^P} \atop }\cr
Q\hskip-10pt&&&&\hskip-35pt R\cr
&\hskip-15pt{\atop \widetilde{x'^{-1}}^{^Q}}\hskip-4pt\nwarrow&&\nearrow\hskip-4pt{\atop \widetilde{(wy)^{-1}}^{^R}}\hskip-15pt&\cr
&&\hskip-13pt {}^{x'}\! Q\cap \,{}^{wy} \!R \hskip-13pt&&\cr}
\eqno £8.11.6\phantom{.}$$
and therefore  we clearly have (cf.~£8.5)
$$\frak t_{\tilde x^{^P},\tilde x'^{^P}}\otimes \frak t_{\tilde y^{^P},\tilde y'^{^P}}
= \bigoplus_{w\in W} \frak  t_{\widetilde{x x'^{-1}}^{^P}\widetilde{,y' y^{-1}w^{-1}}^{\!\!\!^P}}
\eqno £8.11.7.$$ 
Thus,  respectively denoting by $t_{\tilde x^{^P},\tilde x'^{^P}}\,,$  $t_{\tilde y^{^P},\tilde y'^{^P}}$ and 
$t_{\widetilde{x x'^{-1}}^{^P}\widetilde{,y' y^{-1}w^{-1}}^{\!\!\!^P}}$  the isomorphism classes of the 
$\frak H_\T\-$objects  $\frak t_{\tilde x^{^P},\tilde x'^{^P}}\,,$ $\frak t_{\tilde y^{^P},\tilde y'^{^P}}$
and $\frak t_{\widetilde{x x'^{-1}}^{^P}\widetilde{,y' y^{-1}w^{-1}}^{\!\!\!^P}}\,,$
in the $\Bbb Z\-$algebra $\H_\T$ we get
$$t_{\tilde x^{^P},\tilde x'^{^P}}\, t_{\tilde y^{^P},\tilde y'^{^P}}
= \sum_{w\in W}  t_{\widetilde{x x'^{-1}}^{^P}\!,\,\widetilde{y' y^{-1}w^{-1}}^{\!\!\!^P}}
\eqno £8.11.8.$$ 
and therefore we have
$$d_G (t_{\tilde x^{^P},\tilde x'^{^P}}\, t_{\tilde y^{^P},\tilde y'^{^P}}) = 
\sum_{w\in W}  h_{P\,x x'^{-1} w y y'^{-1} P}
\eqno £8.11.9.$$
\eject

\smallskip
On the other hand, we also have
$$\eqalign{d_G (t_{\tilde x^{^P},\tilde x'^{^P}})\, 
d_G ( t_{\tilde y^{^P},\tilde y'^{^P}}) &= h_{P x x'^{-1}P} \,h_{P y y'^{-1}P}\cr
&= {1\over \vert P\vert}\. \sum_{z''\in  P x x'^{-1}Py y'^{-1}P} z''\cr}
\eqno £8.11.10;$$
moreover, we still have $P = \bigsqcup_{w\in W}{}^{x'}\! Q\.w\.\,{}^y R$ and, for any
$u\in Q$ and $v\in R\,,$ we get
$$x \,x'^{-1}({}^{x'}\!u \,w\, {}^yv)\,y \,y'^{-1} =  
{}^x u\,(x \,x'^{-1}w\,y \,y'^{-1}) {}^{y'}\! v
\eqno £8.11.11\phantom{.}$$
where ${}^x u$ and $ {}^{y'}\! v$ still belong to $P\,,$ so that 
$$h_P\,{}^x u = h_P\qq h_P\,{}^{y'}\! v = h_P
\eqno £8.11.12.$$
In conclusion, we obtain
$$\eqalign{d_G (t_{\tilde x^{^P},\tilde x'^{^P}})\, d_G ( t_{\tilde y^{^P},\tilde y'^{^P}}) 
&= {1\over \vert P\vert}\.\sum_{w\in W}\,\sum_{z''\in  P x x'^{-1}wy y'^{-1}P} z''\cr
&= \sum_{w\in W}  h_{P\,x x'^{-1} w y y'^{-1} P} \cr
&= d_G (t_{\tilde x^{^P},\tilde x'^{^P}}\, t_{\tilde y^{^P},\tilde y'^{^P}})\cr}
\eqno £8.11.13.$$
We are done.

\bigskip
\noindent
{\bf Knowledgment.\/}  This paper would not have emerged  without our know-ledge of the work  
{\it Saturated fusion systems as idempotents in the double Burnside ring\/} by Kari Ragnarsson and Radu Stancu [6]. Our effort to understand
the meaning of the {\it Frobenius condition\/} they discovered has led us to the {\it Hecke algebra\/} above.

\bigskip
\noindent
{\bf References}
\bigskip
\noindent
\smallskip\noindent
[1]\phantom{.} Jon Alperin,  {\it Sylow intersections and fusion\/}, Journal of
Algebra, 1(1964), 110-113.
\smallskip\noindent
[2]\phantom{.} James Green, {\it Functors on categories of finite group
representations\/}, Journal of Pure and Applied Algebra, 37(1985), 265-298.
\smallskip\noindent
[3]\phantom{.} Llu\'\i s Puig, {\it Structure locale dans les groupes
finis\/}, Bull.~Soc.~Math.~France, M\'emoire N$^o\,$47(1976)
\smallskip\noindent
[4]\phantom{.} Llu\'\i s Puig, {\it ``Frobenius categories versus Brauer blocks''\/}, Progress in Math.
274(2009), Birkh\"auser, Basel.
\smallskip
\noindent
[5]\phantom{.} Kari Ragnarsson, {\it Classifying spectra of saturated fusion systems\/}, Algebraic and Geometric Topology 6 (2006), 195-252
\smallskip
\noindent
[6]\phantom{.} Kari Ragnarsson and Radu Stancu, {\it Saturated fusion systems as idempotents in the double Burnside ring\/}, preprint

\vfill
\eject
\bigskip

\end